\documentclass[]{article}

\PassOptionsToPackage{%
	maxalphanames	= 4,%
	maxcitenames	= 4,%
	style			= alphabetic%
}{biblatex}

\usepackage[doi=false,isbn=false,url=false,
	hyperref=auto,
	sorting=nyt,
	maxnames=10,
	backend=biber,
	texencoding=auto,
	giveninits=true,
	block=space,
]{biblatex}

\DefineBibliographyStrings{english}{%
	andothers = {et al}
} 
 
\DeclareFieldFormat{eprint:arxiv}{%
	\href{https://arxiv.org/abs/#1}{\emph{arXiv:\allowbreak #1}}}

\DeclareFieldFormat{eprint:mrnumber}{%
	\href{http://www.ams.org/mathscinet-getitem?mr=MR#1}{MR\allowbreak #1}}

\DeclareFieldFormat{eprint:jstor}{%
	\href{http://www.jstor.org/stable/#1}{JSTOR\allowbreak #1}}

\DeclareFieldFormat{eprint:customeprint}{%
	\em Available at \upshape \ttfamily \href{https://www.#1}{#1}}

\DeclareFieldFormat{eprint:online}{%
	Available \href{https://#1}{online}}

\DeclareFieldFormat{eprint:inprep}{%
	\emph{In preparation}}

\DeclareFieldFormat{eprint:manual}{%
	\emph{#1}}

\DeclareFieldFormat{eprint:onarxiv}{%
	\emph{#1}}

\DeclareFieldFormat{eprint:toappear}{%
	\mbox{\emph{To appear}}}

\DeclareFieldFormat{eprint:accepted}{%
	\emph{Accepted at {#1}; to appear}}

\DeclareSourcemap{
	\maps[datatype=bibtex]{
		\map{
			\step[fieldsource=mrnumber,		fieldtarget=eprint, final]
			\step[fieldset=eprinttype,		fieldvalue=mrnumber]
		}
		\map{
			\step[fieldsource=arxiv,		fieldtarget=eprint, final]
			\step[fieldset=eprinttype,		fieldvalue=arxiv]
		}
		\map{
			\step[fieldsource=jstor,		fieldtarget=eprint, final]
			\step[fieldset=eprinttype,		fieldvalue=jstor]
		}
		\map{
			\step[fieldsource=customeprint,	fieldtarget=eprint, final]
			\step[fieldset=eprinttype,		fieldvalue=customeprint]
		}
		\map{
			\step[fieldsource=online,		fieldtarget=eprint, final]
			\step[fieldset=eprinttype,		fieldvalue=online]
		}
		\map{
			\step[fieldsource=inprep,		fieldtarget=eprint, final]
			\step[fieldset=eprinttype,		fieldvalue=inprep]
		}
		\map{
			\step[fieldsource=manual,		fieldtarget=eprint, final]
			\step[fieldset=eprinttype,		fieldvalue=manual]
		}
		\map{
			\step[fieldsource=onarxiv,		fieldtarget=eprint, final]
			\step[fieldset=eprinttype,		fieldvalue=onarxiv]
		}
		\map{
			\step[fieldsource=toappear,		fieldtarget=eprint, final]
			\step[fieldset=eprinttype,		fieldvalue=toappear]
		}
		\map{
			\step[fieldsource=accepted,		fieldtarget=eprint, final]
			\step[fieldset=eprinttype,		fieldvalue=accepted]
		}
	}
}

\DeclareFieldFormat{doi}{%
	\href{https://doi.org/#1}{DOI}%
}

\DeclareFieldFormat{url}{%
	\href{https://#1}{\texttt{#1}}%
}

\DeclareFieldFormat{comments}{%
	\emph{#1}%
} 
 
\DeclareFieldFormat
	[article,inbook,incollection,inproceedings,patent,thesis,unpublished]
	{title}{{#1\isdot}}
\DeclareFieldFormat
	[article,book,inbook,incollection,inproceedings,patent,thesis,unpublished, online]
	{date}{\mkbibparens{#1}}
\DeclareFieldFormat
	[article,book,inbook,incollection,inproceedings,patent,thesis,unpublished]
	{volume}{\mkbibbold{#1}}
\DeclareFieldFormat
	{pages}{\mkbibparens{#1}}

\DeclareFieldFormat{edition}{%
	\ifinteger{#1}
	{\mkbibordedition{#1}~\bibstring{edition}\addcomma}
	{#1~\bibstring{edition}\addcomma}}

\renewbibmacro*{volume+number+eid}{%
	\printfield{volume}%
\setunit*{\adddot}%
	\printfield{number}%
\setunit{\space}%
	\printfield{eid}%
}

\DeclareBibliographyDriver{article}{%
	\usebibmacro{bibindex}%
	\usebibmacro{begentry}%
	\usebibmacro{author}%
\setunit{\addspace}%
	\usebibmacro{date}%
\newunit\newblock
	\usebibmacro{title}%
\newunit\newblock
	\printfield{journaltitle}\addperiod%
\setunit*{\addspace}%
	\usebibmacro{volume+number+eid}%
\setunit*{\addspace}\newblock
	\printfield{pages}
\setunit*{\addspace}%\newblock
	\printfield{comments}%
	\usebibmacro{eprint}%
\setunit*{\addspace}%
	\printfield{doi}%
	\usebibmacro{finentry}%
}

\DeclareBibliographyDriver{online}{%
	\usebibmacro{bibindex}%
	\usebibmacro{begentry}%
	\usebibmacro{author}%
\setunit{\addspace}%
	\usebibmacro{date}%
\newunit\newblock
	\usebibmacro{title}%
\newunit\newblock
	\usebibmacro{eprint}%
	\usebibmacro{finentry}%
}

\DeclareBibliographyDriver{book}{%
	\usebibmacro{bibindex}%
	\usebibmacro{begentry}%
	\usebibmacro{author}%
\setunit{\addspace}\newblock
	\usebibmacro{date}%
\newunit\newblock
	\usebibmacro{title}%
\setunit*{\addspace}\newblock
	\printfield{edition}%
\newunit\newblock
	\printlist{publisher}%\addcomma
\setunit*{\addspace}%
	\printfield{volume}%
\setunit*{\addspace}\newblock%
	\printfield{comments}%
	\usebibmacro{eprint}%
\setunit*{\addspace}
	\printfield{doi}
	\usebibmacro{finentry}%
}

\DeclareBibliographyDriver{thesis}{%
	\usebibmacro{bibindex}%
	\usebibmacro{begentry}%
	\usebibmacro{author}%
\setunit{\addspace}\newblock
	\usebibmacro{date}%
\newunit\newblock
	\usebibmacro{title}.%
\setunit*{\addspace}\newblock
	\space Thesis, \printlist{institution}
	\usebibmacro{eprint}%
	\printfield{comments}%
\setunit*{\addspace}
	\printfield{doi}%
	\usebibmacro{finentry}%
}

\DeclareBibliographyDriver{inproceedings}{%
	\usebibmacro{bibindex}%
	\usebibmacro{begentry}%
	\usebibmacro{author}%
\setunit{\addspace}%
	\usebibmacro{date}%
\newunit\newblock
	\usebibmacro{title}%
\newunit\newblock
	\printfield{booktitle}\addcomma%
\setunit*{\addspace}%
	\printfield{series}\addcomma%
\setunit*{\addspace}\newblock
	\usebibmacro{editor}\addcomma
\setunit*{\addspace}\newblock
	\printlist{publisher}
\setunit*{\addspace}%
	\printfield{volume}%
\setunit*{\addspace}%
	\printfield{pages}
\setunit*{\addspace}%\newblock
	\usebibmacro{eprint}%
	\printfield{comments}%
\setunit*{\addnnspace}
	\printfield{doi}
	\usebibmacro{finentry}%
}

\DeclareBibliographyDriver{inbook}{%
	\usebibmacro{bibindex}%
	\usebibmacro{begentry}%
	\usebibmacro{author}%
\setunit{\addspace}%
	\usebibmacro{date}%
\newunit\newblock
	\usebibmacro{title}%
\newunit\newblock
	\printfield{booktitle}\addcomma%
\setunit*{\addspace}%
	\printfield{series}\addcomma%
\setunit*{\addspace}\newblock
	\usebibmacro{editor}\addcomma
\setunit*{\addspace}\newblock
	\printlist{publisher}%\addcomma%
\setunit*{\addspace}%
	\printfield{volume}%
\setunit*{\addspace}%
	\printfield{pages}
\setunit*{\addspace}%\newblock
	\usebibmacro{eprint}%
	\printfield{comments}%
\setunit*{\addspace}
	\printfield{doi}
	\usebibmacro{finentry}%
}

\DeclareBibliographyDriver{incollection}{%
	\usebibmacro{bibindex}%
	\usebibmacro{begentry}%
	\usebibmacro{author}%
\setunit{\addspace}%
	\usebibmacro{date}%
\newunit\newblock
	\usebibmacro{title}%
\newunit\newblock
	\printfield{booktitle}\addcomma%
\setunit*{\addspace}%
	\printfield{series}\addcomma%
\setunit*{\addspace}\newblock
	\printlist{publisher}%\addcomma%
\setunit*{\addspace}%
	\printfield{volume}%
\setunit*{\addspace}%
	\printfield{pages}
\setunit*{\addspace}%\newblock
	\usebibmacro{eprint}%
	\printfield{comments}%
\setunit*{\addspace}
	\printfield{doi}
	\usebibmacro{finentry}%
}

\AtEveryBibitem{%
	\clearfield{day}%
	\clearfield{month}%
} 
 
\DeclareNameAlias{sortname}{given-family}
\DeclareNameAlias{default}{given-family}

\usepackage{geometry}
\usepackage[hidelinks]{hyperref}
\hypersetup{
	colorlinks,
	linkcolor={blue},
	citecolor={ForestGreen},
	urlcolor={blue}
}

\usepackage{graphicx}
\graphicspath{{Figures/}}

\usepackage[%
	labelfont	= {bf, sf},
	labelformat	= simple,
	labelsep	= period,
	singlelinecheck	= false%
]{caption}

\usepackage[%
	labelformat	= simple%
]{subcaption}

\usepackage{fancyhdr}
\usepackage{fmtcount}
\usepackage{lastpage}
\fancypagestyle{plain}{
	\fancyhf{}
	\fancyfoot[R]{\sffamily Page~\decimal{page} of~\pageref*{LastPage}}
	
}
\newcommand{\blfootnote}[1]{%
	\begingroup
	\renewcommand\thefootnote{}\footnote{\sffamily#1}%
	\addtocounter{footnote}{-1}%
	\endgroup
}
\usepackage{titlesec}
\titleformat{\section}
{\sffamily \Large \bfseries \boldmath}{\thesection}{1em}{}
\titleformat{\subsection}
{\sffamily \large \bfseries \boldmath}{\thesubsection}{1em}{}
\titleformat{\subsubsection}
{\sffamily \normalsize \bfseries \boldmath}{\thesubsubsection}{1em}{}
\titleformat{\paragraph}[runin]
{\sffamily \normalsize \bfseries \boldmath}{\theparagraph}{1em}{}
\titleformat{\subparagraph}[runin]
{\sffamily \normalsize \bfseries \boldmath}{\thesubparagraph}{1em}{}
\usepackage[utf8]{inputenc}

\usepackage[UKenglish]{babel}
\usepackage{amsmath, amsthm, amssymb, mathrsfs, bm, mathtools}
\usepackage{wasysym}
\allowdisplaybreaks

\usepackage{titlesec}
\usepackage{xifthen}
\usepackage{xspace}

\usepackage{footnotebackref}
\usepackage{enumitem}

\setlist[itemize,enumerate]{%
	itemsep	= 0pt,%
	topsep	= \smallskipamount%
}

\newcommand{\romannumbering}{%
	\renewcommand{\labelenumi}{\upshape(\roman{enumi})}%
	\renewcommand{\theenumi}{\upshape(\roman{enumi})}}

\usepackage{scrextend}

\setlist[description]{%
	topsep = \smallskipamount,		% space before start / after end of list
	itemsep = \smallskipamount,		% space between items
	font = {\mdseries\itshape},		% set the label font
	leftmargin = \parindent
}

\newlist{ranlist}{enumerate}{3}
\setlist[ranlist,1]{
	label=(\roman*) }
\setlist[ranlist,2]{
	label=(\textit{\alph*}),
	ref=(\roman{ranlisti}.\textit{\alph*})	}
\setlist[ranlist,3]{
	label=\arabic*.,
	ref=(\roman{ranlisti}.\textit{\alph{ranlistii}}.\arabic*)	}
\usepackage[dvipsnames,hyperref]{xcolor}
\definecolor{darkorange}{rgb}{1.0, 0.55, 0.0}

\def\IfAmpersandUseAlign#1#2&#3\EndIfAmpersandUseAlign%
{%
	\if\relax\detokenize{#3}\relax
	\begin{equation*}%
	#1%
	\end{equation*}%
	\else
	\begin{align*}%
	#1%
	\end{align*}%
	\fi
}
\def\[#1\]%
{%
	\IfAmpersandUseAlign{#1}#1&\EndIfAmpersandUseAlign
}

\newcommand{\tvt}[1]{\| #1 \|_{\textsf{\textup{TV}}}}
\newcommand{\tv}[1]{
	\mathchoice
	{\bigl\| #1 \bigr\rVert_{\textup{\textsf{TV}}}}%
	{\lVert #1 \|_{\textup{\textsf{TV}}}}%
	{\| #1 \|_{\textup{\textsf{TV}}}}%
	{\| #1 \|_{\textup{\textsf{TV}}}}
}

\newcommand{\ONE}     {\bm1}
\newcommand{\one}  [1]{\bm1\{ #1 \}}

\let\originalexp\exp
\let\exp\relax
\DeclareRobustCommand{\exp} [1]{\originalexp(#1)}

\newcommand{\expb}[2][]{
	\ifthenelse{\equal{}{#1}}
	{\originalexp\bigl( #2 \bigr)}
	{\originalexp_{#1}\bigl( #2 \bigr)}
}

\newcommand{\logn}[1][]{
	\ifthenelse{\equal{}{#1}}
	{\log n}
	{(\log n)^{#1}}
}
\newcommand{\logk}[1][]{
	\ifthenelse{\equal{}{#1}}
	{\log k}
	{(\log k)^{#1}}
}

\newcommand{\Quad}[1]{
	\mathchoice
	{\quad\text{#1}\quad}%	\displaystyle
	{\text{ #1 }}%			\textsyle
	{\text{ #1 }}%			\scriptstyle
	{\text{ #1 }}%			\scriptscriptstyle
}

\newcommand{\Qand}{\Quad{and}}
\newcommand{\Qfor}{\Quad{for}}
\newcommand{\Qforall}{\Quad{for all}}
\newcommand{\Qwhere}{\Quad{where}}

\newcommand{\Qthen}{\Quad{then}}

\newcommand{\fnrestrict}[2]{{% we make the whole thing an ordinary symbol
		\left.\kern-\nulldelimiterspace % automatically resize the bar with \right
		#1 % the function
		\vphantom{\big|} % pretend it's a little taller at normal size
		\right|_{#2} % this is the delimiter
}}

\usepackage{calc}
\newlength{\halfplusheight}
\newcommand{\MAX}[1]{%
	\mathop{\raisebox{\halfplusheight}{\(\displaystyle\max_{#1}\)}}\:%
}
\newcommand{\SUP}[1]{%
	\mathop{\raisebox{\halfplusheight}{\(\displaystyle\sup_{#1}\)}}\:%
}

\newcommand{\LIM}[1]{%
	\mathop{\raisebox{\halfplusheight}{\(\displaystyle\lim_{#1}\)}}\:%
}
\newcommand{\LIMSUP}[1]{%
	\mathop{\raisebox{\halfplusheight}{\(\displaystyle\limsup_{#1}\)}}\:%
}
\newcommand{\LIMINF}[1]{%
	\mathop{\raisebox{\halfplusheight}{\(\displaystyle\liminf_{#1}\)}}\:%
}

\newcommand{\bcdot}{\ensuremath{\bm{\cdot}}}
\newcommand{\abs}  [1]{| #1 |}
\newcommand{\absb} [1]{\big| #1 \bigr|}

\newcommand{\rbr} [1]{ ( #1 ) }
\newcommand{\rbb} [1]{\bigl( #1 \bigr)}

\newcommand{\bra} [1]{ \{ #1 \} }
\newcommand{\brb} [1]{\bigl\{ #1 \bigr\}}

\newcommand{\floor}  [1]{\lfloor #1 \rfloor}

\newcommand{\ceil}[1]  {\lceil #1 \rceil}

\newcommand{\midb}{\bigm|}

\newcommand{\Oh}  [1]{\mathcal{O}( #1 )}
\newcommand{\Ohb} [1]{\mathcal{O}\bigl( #1 \bigr)}

\newcommand{\oh}  [1]{o( #1 )}

\newcommand{\Th}  [1]{\Theta( #1 )}

\newcommand{\pr}[2][]{
	\mathchoice
	{\ifthenelse{\isempty{#1}}
		{\mathbb{P}\bigl(#2\bigr)}
		{\mathbb{P}_{#1}\bigl(#2\bigr)}}%	\displaystyle
	{\ifthenelse{\isempty{#1}}
		{\mathbb{P}(#2)}
		{\mathbb{P}_{#1}(#2)}}%	\textsyle
	{\ifthenelse{\isempty{#1}}
		{\mathbb{P}(#2)}
		{\mathbb{P}_{#1}(#2)}}%	\scriptstyle
	{\ifthenelse{\isempty{#1}}
		{\mathbb{P}(#2)}
		{\mathbb{P}_{#1}(#2)}}%	\scriptscriptstyle
}
\newcommand{\prt}[2][]{
	\ifthenelse{\equal{}{#1}}
	{\mathbb{P}( #2 )}
	{\mathbb{P}_{#1}( #2 )}
}
\newcommand{\prb}[2][]{
	\ifthenelse{\equal{}{#1}}
	{\mathbb{P}\bigl( #2 \bigr)}
	{\mathbb{P}_{#1}\bigl( #2 \bigr)}
}
\newcommand{\prB}[2][]{
	\ifthenelse{\equal{}{#1}}
	{\mathbb{P}\Bigl( #2 \Bigr)}
	{\mathbb{P}_{#1}\Bigl( #2 \Bigr)}
}
\newcommand{\prbb}[2][]{
	\ifthenelse{\equal{}{#1}}
	{\mathbb{P}\biggl( #2 \biggr)}
	{\mathbb{P}_{#1}\biggl( #2 \biggr)}
}
\newcommand{\prBB}[2][]{
	\ifthenelse{\equal{}{#1}}
	{\mathbb{P}\Biggl( #2 \Biggr)}
	{\mathbb{P}_{#1}\Biggl( #2 \Biggr)}
}
\newcommand{\prs}[2][]{
	\ifthenelse{\equal{}{#1}}
	{\mathbb{P}\left( #2 \right)}
	{\mathbb{P}_{#1}\left( #2 \right)}
}

\newcommand{\ex}[2][]{
	\mathchoice
	{\ifthenelse{\isempty{#1}}
		{\mathbb{E}\bigl(#2\bigr)}
		{\mathbb{E}_{#1}\bigl(#2\bigr)}}%	\displaystyle
	{\ifthenelse{\isempty{#1}}
		{\mathbb{E}(#2)}
		{\mathbb{E}_{#1}(#2)}}%	\textsyle
	{\ifthenelse{\isempty{#1}}
		{\mathbb{E}(#2)}
		{\mathbb{E}_{#1}(#2)}}%	\scriptstyle
	{\ifthenelse{\isempty{#1}}
		{\mathbb{E}(#2)}
		{\mathbb{E}_{#1}(#2)}}%	\scriptscriptstyle
}
\newcommand{\ext}[2][]{
	\ifthenelse{\equal{}{#1}}
	{\mathbb{E}( #2 )}
	{\mathbb{E}_{#1}( #2 )}
}
\newcommand{\exb}[2][]{
	\ifthenelse{\equal{}{#1}}
	{\mathbb{E}\bigl( #2 \bigr)}
	{\mathbb{E}_{#1}\bigr( #2 \bigr)}
}
\newcommand{\exB}[2][]{
	\ifthenelse{\equal{}{#1}}
	{\mathbb{E}\Bigl( #2 \Bigr)}
	{\mathbb{E}_{#1}\Bigl( #2 \Bigr)}
}
\newcommand{\exbb}[2][]{
	\ifthenelse{\equal{}{#1}}
	{\mathbb{E}\biggl( #2 \biggr)}
	{\mathbb{E}_{#1}\biggl( #2 \biggr)}
}
\newcommand{\exBB}[2][]{
	\ifthenelse{\equal{}{#1}}
	{\mathbb{E}\Biggl( #2 \Biggr)}
	{\mathbb{E}_{#1}\Biggl( #2 \Biggr)}
}

\newcommand{\Varb}[2][]{
	\ifthenelse{\equal{}{#1}}
	{\mathbb{V}\textnormal{ar} \bigl(#2\bigr)}
	{\mathbb{V}\textnormal{ar}_{#1} \bigl(#2\bigr)}
}
\newcommand{\VAR}[2][]{
	\ifthenelse{\equal{}{#1}}
	{\textnormal{Var}(#2)}
	{\textnormal{Var}_{#1}(#2)}
}
\newcommand{\sumt}[2][]{
	\mathchoice
	{\ifthenelse{\isempty{#1}}
		{\textstyle \sum_{#2}      \displaystyle}
		{\textstyle \sum_{#2}^{#1} \displaystyle}}%	\displaystyle
	{\ifthenelse{\isempty{#1}}
		{\sum_{#2}}
		{\sum_{#2}^{#1}}}%			\textsyle
	{\ifthenelse{\isempty{#1}}
		{\sum_{#2}}
		{\sum_{#2}^{#1}}}%			\scriptstyle
	{\ifthenelse{\isempty{#1}}
		{\sum_{#2}}
		{\sum_{#2}^{#1}}}%			\scriptscriptstyle
}
\newcommand{\sumd}[2][]{
	\ifthenelse{\isempty{#1}}
		{\sum_{#2}}
		{\sum_{#2}^{#1}}
}

\newcommand{\intt}[2][]{
	\mathchoice
	{\ifthenelse{\isempty{#1}}
		{\textstyle \int_{#2}      \displaystyle}
		{\textstyle \int_{#2}^{#1} \displaystyle}}
	{\ifthenelse{\isempty{#1}}
		{\int_{#2}}
		{\int_{#2}^{#1}}}
	{\ifthenelse{\isempty{#1}}
		{\int_{#2}}
		{\int_{#2}^{#1}}}
	{\ifthenelse{\isempty{#1}}
		{\int_{#2}}
		{\int_{#2}^{#1}}}
}

\newcommand{\prodt}[2][]{
	\mathchoice
	{\ifthenelse{\isempty{#1}}
		{\textstyle \prod_{#2}      \displaystyle}
		{\textstyle \prod_{#2}^{#1} \displaystyle}}
	{\ifthenelse{\isempty{#1}}
		{\prod_{#2}}
		{\prod_{#2}^{#1}}}
	{\ifthenelse{\isempty{#1}}
		{\prod_{#2}}
		{\prod_{#2}^{#1}}}
	{\ifthenelse{\isempty{#1}}
		{\prod_{#2}}
		{\prod_{#2}^{#1}}}	
}
\newcommand{\prodd}[2][]{
	\ifthenelse{\isempty{#1}}
		{\prod_{#2}}
		{\prod_{#2}^{#1}}
}
\newcommand{\toinf}[1]{\ensuremath{#1\to\infty}\xspace}
\newcommand{\asinf}[1]{\text{as \(#1\to\infty\)}\xspace}

\newcommand{\tozero}[1]{\ensuremath{#1\to0}\xspace}
\newcommand{\aszero}[1]{\text{as \(#1\to0\)}\xspace}

\DeclareMathOperator{\Unif}{Unif}

\DeclareMathOperator{\Bern}{Bern}

\newcommand{\iid}{\textup{\textsf{iid}}\xspace}

\newcommand{\TV}{\textsf{TV}\xspace}

\newcommand{\tmix}{t^\mix}

\newcommand{\mix}{\textup{\textsf{mix}}}

\newcommand{\MC}[2]{(#1_#2)_{#2\ge0}}
\newcommand{\mbn}{\mathbb{N}}

\newcommand{\mbr}{\mathbb{R}}

\newcommand{\mbz}{\mathbb{Z}}

\newcommand{\mcc}{\mathcal{C}}

\newcommand{\mcl}{\mathcal{L}}

\newcommand{\mcp}{\mathcal{P}}

\newcommand{\mfc}{\mathfrak{C}}

\newcommand{\mfm}{\mathfrak{M}}

\newcommand{\mfp}{\mathfrak{P}}

\newcommand{\mfs}{\mathfrak{S}}

\newcommand{\cq}{\coloneqq}

\newenvironment{Proof}[1][\proofname]{%
	\proof[\upshape\bfseries\sffamily\boldmath#1]
}{\endproof}

\newcommand{\eps}{\varepsilon}

\makeatletter
\let\save@mathaccent\mathaccent
\newcommand*\if@single[3]{%
	\setbox0\hbox{${\mathaccent"0362{#1}}^H$}%
	\setbox2\hbox{${\mathaccent"0362{\kern0pt#1}}^H$}%
	\ifdim\ht0=\ht2 #3\else #2\fi
}
\newcommand*\rel@kern[1]{\kern#1\dimexpr\macc@kerna}
\newcommand*\widebar[1]{\@ifnextchar^{{\wide@bar{#1}{0}}}{\wide@bar{#1}{1}}}
\newcommand*\wide@bar[2]{\if@single{#1}{\wide@bar@{#1}{#2}{1}}{\wide@bar@{#1}{#2}{2}}}
\newcommand*\wide@bar@[3]{%
	\begingroup
	\def\mathaccent##1##2{%
		\let\mathaccent\save@mathaccent
		\if#32 \let\macc@nucleus\first@char \fi
		\setbox\z@\hbox{$\macc@style{\macc@nucleus}_{}$}%
		\setbox\tw@\hbox{$\macc@style{\macc@nucleus}{}_{}$}%
		\dimen@\wd\tw@
		\advance\dimen@-\wd\z@
		\divide\dimen@ 3
		\@tempdima\wd\tw@
		\advance\@tempdima-\scriptspace
		\divide\@tempdima 10
		\advance\dimen@-\@tempdima
		\ifdim\dimen@>\z@ \dimen@0pt\fi
		\rel@kern{0.6}\kern-\dimen@
		\if#31
		\overline{\rel@kern{-0.6}\kern\dimen@\macc@nucleus\rel@kern{0.4}\kern\dimen@}%
		\advance\dimen@0.4\dimexpr\macc@kerna
		\let\final@kern#2%
		\ifdim\dimen@<\z@ \let\final@kern1\fi
		\if\final@kern1 \kern-\dimen@\fi
		\else
		\overline{\rel@kern{-0.6}\kern\dimen@#1}%
		\fi
	}%
	\macc@depth\@ne
	\let\math@bgroup\@empty \let\math@egroup\macc@set@skewchar
	\mathsurround\z@ \frozen@everymath{\mathgroup\macc@group\relax}%
	\macc@set@skewchar\relax
	\let\mathaccentV\macc@nested@a
	\if#31
	\macc@nested@a\relax111{#1}%
	\else
	\def\gobble@till@marker##1\endmarker{}%
	\futurelet\first@char\gobble@till@marker#1\endmarker
	\ifcat\noexpand\first@char A\else
	\def\first@char{}%
	\fi
	\macc@nested@a\relax111{\first@char}%
	\fi
	\endgroup
}
\makeatother
\numberwithin{equation}{subsection}


\newcounter{parentnumber}
\newcommand{\qedtriangle}{\renewcommand{\qedsymbol}{\ensuremath{\triangle}}}

\usepackage[capitalise]{cleveref}

\newtheoremstyle{sfsl}
{1\baselineskip}		% Space above
{1\baselineskip}		% Space below
{\slshape}				% Theorem body font
{}						% Indent amount
{\bfseries\sffamily}	% Theorem head font
{.}						% Punctuation after theorem head
{0.5em}					% Space after theorem head
{\thmname{#1}\thmnumber{ #2}\thmnote{: \textnormal{\sffamily#3}}}
\newtheoremstyle{sfup}
{1\baselineskip}		% Space above
{1\baselineskip}		% Space below
{\upshape}				% Theorem body font
{}						% Indent amount
{\bfseries\sffamily}	% Theorem head font
{.}						% Punctuation after theorem head
{0.5em}					% Space after theorem head
{\thmname{#1}\thmnumber{ #2}\thmnote{: \textnormal{\sffamily#3}}}
\theoremstyle{sfsl}

\newtheorem*{thm*}{Theorem}
\newtheorem{thm} {Theorem}[section]
\crefname{thm}{Theorem}{Theorems}

\newtheorem*{introthm*}{Theorem}

\crefname{introthm}{Theorem}{Theorems}

\newtheorem*{cor*}{Corollary}
\newtheorem{cor} [thm]{Corollary}
\crefname{cor}{Corollary}{Corollaries}

\newtheorem*{introcor*}{Corollary}

\crefname{introcor}{Corollary}{Corollaries}

\newtheorem*{introconj*}{Conjecture}

\crefname{introconj}{Conjecture}{Conjectures}

\newtheorem*{introques*}{Question}

\crefname{introques}{Question}{Questions}

\newtheorem*{lem*}    {Lemma}
\newtheorem{lem} [thm]{Lemma}
\crefname{lem}{Lemma}{Lemmas}

\newtheorem*{introlem*}{Lemma}

\crefname{introlem}{Lemma}{Lemmas}

\newtheorem*{prop*}    {Proposition}
\newtheorem{prop} [thm]{Proposition}
\crefname{prop}{Proposition}{Propositions}

\newtheorem*{defn*}    {Definition}
\newtheorem{defn} [thm]{Definition}
\crefname{defn}{Definition}{Definitions}

\newtheorem*{introdefn*}{Definition}

\crefname{introdefn}{Definition}{Definitions}

\newtheorem*{nota*}{Notation}
\newtheorem{nota}[thm]{Notation}

\newtheorem*{term*}{Terminology}
\newtheorem{term}[thm]{Terminology}

\newtheorem{innercustomgeneric}{\customgenericname}
\providecommand{\customgenericname}{}
\newcommand{\newcustomtheorem}[2]{%
	\newenvironment{#1}[1]
	{%
		\renewcommand\customgenericname{#2}%
		\renewcommand\theinnercustomgeneric{##1}%
		\innercustomgeneric
	}
	{\endinnercustomgeneric}
}

\newcustomtheorem{customthm} {Theorem}
\newcustomtheorem{customprop}{Proposition}
\newcustomtheorem{customlem} {Lemma}
\newcustomtheorem{customcor} {Corollary}
\newcustomtheorem{customres}{\!\!}

\theoremstyle{sfup}

\newtheorem{exmT}[thm]{Example}
\newenvironment{exmt}
{\pushQED{\qed}\renewcommand{\qedsymbol}{$\triangle$}\exmT}
{\popQED\endexmT}

\newtheorem*{exm*} {Example}
\newtheorem*{exmTT}{Lemma}
\newenvironment{exmt*}
{\pushQED{\qed}\renewcommand{\qedsymbol}{$\triangle$}\exmTT}
{\popQED\endexmTT}

\newtheorem{rmk} [thm]{Remark}
\newtheorem{rmkT}[thm]{Remark}
\newenvironment{rmkt}
{\pushQED{\qed}\renewcommand{\qedsymbol}{$\triangle$}\rmkT}
{\popQED\endrmkT}

\newtheorem*{rmk*} {Remark}
\newtheorem*{rmkTT}{Remark}
\newenvironment{rmkt*}
{\pushQED{\qed}\renewcommand{\qedsymbol}{$\triangle$}\rmkTT}
{\popQED\endrmkTT}

\crefname{intrormk} {Remark}{Remarks}
\crefname{intrormkT}{Remark}{Remarks}

\newtheorem*{intrormk*} {Remark}
\newtheorem*{intrormkTT}{Remark}
\newenvironment{intrormkt*}
{\pushQED{\qed}\renewcommand{\qedsymbol}{\ensuremath{\triangle}}\intrormkTT}
{\popQED\endintrormkTT}

\newtheorem*{rmks*} {Remarks}
\newtheorem*{rmksTT}{Remarks}
\newenvironment{rmkst*}
{\pushQED{\qed}\renewcommand{\qedsymbol}{$\triangle$}\rmksTT}
{\popQED\endrmksTT}

\newtheorem*{note*} {Note}
\newtheorem*{noteTT}{Note}
\newenvironment{notet*}
{\pushQED{\qed}\renewcommand{\qedsymbol}{$\triangle$}\noteTT}
{\popQED\endnoteTT}

\newtheorem{algT}[thm]{Algorithm}
\newenvironment{algt}
	{\pushQED{\qed}\renewcommand{\qedsymbol}{$\triangle$}\algT}
	{\popQED\endalgT}

\crefname{thm}{Theorem}{Theorems}
\crefname{lem}{Lemma}{Lemmas}
\crefname{cor}{Corollary}{Corollaries}
\crefname{prop}{Proposition}{Propositions}
\crefname{algT}{Algorithm}{Algorithms}

\crefname{figure}{Figure}{Figures}
\numberwithin{figure}{section}

\crefformat{section}{\S#2#1#3}
\crefformat{subsection}{\S#2#1#3}
\crefformat{subsubsection}{\S#2#1#3}

\newenvironment{center-small}
	{\par\centering\smallskip}
	{\par\smallskip}

\newcommand{\cf}{coalescence--\allowbreak fragmentation\xspace}

\newcommand{\bsc}[1][]{\cite[#1]{BS:cutoff-conj-inv}\xspace}
\newcommand{\bst}[1][]{\textcite[#1]{BS:cutoff-conj-inv}\xspace}
\newcommand{\bsn}{\citeauthor{BS:cutoff-conj-inv}\xspace}

\newcommand{\bszc}[1][]{\cite[#1]{BSZ:k-cycle}\xspace}
\newcommand{\bszt}[1][]{\textcite[#1]{BSZ:k-cycle}\xspace}

\newcommand{\bklmc}[1][]{\cite[#1]{BKLM:interchange-rev}\xspace}
\newcommand{\bklmt}[1][]{\textcite[#1]{BKLM:interchange-rev}\xspace}

\renewcommand{\aa}{\ensuremath{a}}

\newcommand{\bb}{\ensuremath{b}}
\newcommand{\BB}{\ensuremath{B}}
\newcommand{\cc}{\ensuremath{c}}
\newcommand{\Eta}{H}
\newcommand{\GG}{\ensuremath{G}}
\newcommand{\ii}{\ensuremath{i}}
\newcommand{\II}{\ensuremath{I}}
\newcommand{\jj}{\ensuremath{j}}
\newcommand{\JJ}{\ensuremath{J}}
\newcommand{\kk}{\ensuremath{k}}
\newcommand{\KK}{\ensuremath{K}}
\newcommand{\mm}{\ensuremath{m}}
\newcommand{\MM}{\ensuremath{M}}
\newcommand{\nn}{\ensuremath{n}}

\newcommand{\pp}{\ensuremath{\lambda}}
\newcommand{\PP}{\ensuremath{P}}
\newcommand{\qq}{\ensuremath{\mu}}
\newcommand{\QQ}{\ensuremath{Q}}
\newcommand{\rr}{\ensuremath{r}}

\renewcommand{\ss}{\ensuremath{s}}
\renewcommand{\SS}{\ensuremath{S}}
\renewcommand{\tt}{\ensuremath{t}}
\newcommand{\TT}{\ensuremath{T}}
\newcommand{\uu}{\ensuremath{u}}
\newcommand{\UU}{\ensuremath{U}}
\newcommand{\vv}{\ensuremath{v}}

\newcommand{\ww}{\ensuremath{w}}

\newcommand{\Mu}{\ensuremath{M}}
\newcommand{\Nu}{\ensuremath{N}}

\newcommand{\id}{\textup{\textsf{id}}\xspace}
\newcommand{\whp}{\textsf{whp}\xspace}

\newcommand{\wrt}{\textsf{wrt}\xspace}
\newcommand{\uar}{\textsf{uar}\xspace}

\newcommand{\CM}{\textsf{CM}\xspace}

\newcommand{\CS}{\textsf{CS}\xspace}
\newcommand{\CSs}{\textsf{CS}s\xspace}
\newcommand{\IP}{\textsf{IP}\xspace}

\newcommand{\LLN}{\textsf{LLN}\xspace}
\newcommand{\LLNs}{\textsf{LLN}s\xspace}
\newcommand{\PM}{\textsf{PM}\xspace}
\newcommand{\PMs}{\textsf{PM}s\xspace}
\newcommand{\RT}{\textsf{RT}\xspace}

\newcommand{\RW}{\textsf{RW}\xspace}
\newcommand{\RWs}{\textsf{RW}s\xspace}

\DeclareMathOperator{\PD}{PD}

\newcommand{\oci} [1]{\mathopen( #1 \mathclose]}

\fancypagestyle{body}{
	\fancyhf{}
	\fancyhead[L]{\sffamily Cutoff for Rewiring Dynamics on Perfect Matchings}
	\fancyhead[R]{\sffamily Sam Olesker-Taylor}
	\fancyfoot[L]{\sffamily \nouppercase \leftmark}
	\fancyfoot[R]{\sffamily Page~\decimal{page} of~\pageref*{LastPage}}
	
}
\title{\sffamily%
	Cutoff for Rewiring Dynamics on Perfect Matchings
}
\author{\sffamily%
	Sam Olesker-Taylor%\quad%
}
\date{}
\begin{document}

\thispagestyle{plain}

\maketitle

\renewcommand{\abstractname}{\sffamily Abstract}

\begin{abstract}
We establish cutoff for a natural random walk (\textit{\RW}) on the set of perfect matchings (\textit{\PM}s), based on `rewiring'.
An $\nn$-\PM is a pairing of $2\nn$ objects.
The \textit{$\kk$-\PM \RW}
	selects $\kk$ pairs uniformly at random,
	disassociates the corresponding $2\kk$ objects,
then
	chooses a new pairing on these $2\kk$ objects uniformly at random.
The equilibrium distribution is uniform over all $\nn$-\PMs.

The $2$-\PM \RW was first introduced by Diaconis and Holmes~\cite{DH:matchings-trees,DH:rw-matchings}, seen as a \RW on phylogenetic trees.
They established cutoff in this case.
We establish cutoff for the $\kk$-\PM \RW whenever $2 \le \kk \ll \nn$.
If $\kk \gg 1$, then the mixing time is $\tfrac\nn\kk \log \nn$ to leading order.

Diaconis and Holmes~\cite{DH:rw-matchings} relate the $2$-\PM \RW to the random transpositions card shuffle.
Ceccherini-Silberstein, Scarabotti and Tolli~\cite{CsST:gelfand-applications,CsST:harmonic-analysis-finite-groups} establish the same result using representation theory.
We are the first to handle $\kk > 2$.
We relate the \PM \RW to conjugacy-invariant \RWs on the permutation group by introducing a `cycle structure' for \PMs, then build on work of Berestycki, Schramm, Şengül and Zeitouni~\cite{S:random-trans,BSZ:k-cycle,BS:cutoff-conj-inv} on such \RWs.
\end{abstract}

%\begin{abstract}
%	%
%We establish cutoff for a natural random walk (RW) on the set of perfect matchings (PMs), based on 'rewiring'. An $n$-PM is a pairing of $2n$ objects. The $k$-PM RW selects $k$ pairs uniformly at random, disassociates the corresponding $2k$ objects, then chooses a new pairing on these $2k$ objects uniformly at random. The equilibrium distribution is uniform over all $n$-PM.
%
%The $2$-PM RW was first introduced by Diaconis and Holmes (1998, 2002), seen as a RW on phylogenetic trees. They established cutoff in this case. We establish cutoff for the $k$-PM RW whenever $2 \le k \ll n$. If $k \gg 1$, then the mixing time is $\tfrac nk \log n$ to leading order.
%
%Diaconis and Holmes (2002) relate the $2$-PM RW to the random transpositions card shuffle. Ceccherini-Silberstein, Scarabotti and Tolli (2007, 2008) establish the same result using representation theory. We are the first to handle $k > 2$. We relate the PM RW to conjugacy-invariant RWs on the permutation group by introducing a 'cycle structure' for PMs, then build on work of Berestycki, Schramm, Şengül and Zeitouni (2005, 2011, 2019) on such RWs.
%	%
%\end{abstract}

\small
\begin{quote}
\begin{description}
	\item [Keywords:]
	mixing time,
	cutoff,
	random walks,
	perfect matchings,
	coalescence--fragmentation,
	random transpositions,
	conjugacy-invariant random walks
	
	\item [MSC 2020 subject classifications:]
	60B15;
	60C05;
	60J10, 60J90
\end{description}
\end{quote}
\normalsize

%05C12: Distance in graphs
%05C48: Expander graphs
%05C80:	Random graphs
%05C81:	Random walks on graphs

%20C15: Ordinary representations and characters
%20C30: Representations of finite symmetric groups

%42A61: Probilistic methods in Fourier analysis

%43A30: Fourier and Fourier-Stieltjes transforms on non-Abelian groups and on semigroups, etc
%43A65: Representations of groups, semigroups, etc
%43A75: Analysis on specific compact groups
%43A90: Spherical functions (in Abstract harmonic analysis)

%60B15: Probability measures on groups or semigroups, Frouier transforms, factorization
%60C05: Combinatorial probability
%60G50: Sums of independent random variables; random walks
%60J10: Markov chains (discrete-time Markov processes on discrete state spaces)
%60J20: Applications of Markov chains and discrete-time Markov processes on general state spaces
%60J90: Coalescent processes
%60J27:	Continuous-time Markov processes on discrete state spaces
%60K35: Interacting random processes; statistical mechanics type models; percolation theory
%60K37:	Processes in random environments

\blfootnote{{}\par%
	Sam Olesker-Taylor,\quad%
	\href{mailto:oleskertaylor.sam@gmail.com}{oleskertaylor.sam@gmail.com}%
\\
	Department of Mathematical Sciences, University of Bath, UK%
\\
	Statistical Laboratory, DPMMS, University of Cambridge, UK%
\\
	Research supported by EPSRC grant EP/N004566/1
}

\vspace{1cm}

\vfill
%\newpage
\sffamily
\setcounter{tocdepth}{1}
\tableofcontents
\normalfont

\vspace*{\bigskipamount}

\romannumbering

\newpage
\section{Introduction}
\label{sec:intro}

\subsection{Model Set-Up}
\label{sec:intro:set-up}

%Let $\nn \in \mbn$.
We analyse a random walk (\textit{\RW}) on the set of perfect matchings (\textit{\PM}s) on $2\nn$ objects, for $\nn \in \mbn$.
We represent a \PM $\eta$ on $2\nn$ objects by a collection of unordered pairs:
\[
	\mfm_\nn
\cq
	\brb{
		\eta = \cup_{\ell=1}^\nn \bra{ \bra{\eta_{2\ell-1}, \eta_{2\ell}} }
	\midb
		\eta_\ii \in [2\nn] \ \forall \: \ii \in [2n], \;
		\cup_{\ii=1}^{2\nn} \bra{\eta_i} = [2\nn]
	}
\Qfor
	\nn \in \mbn;
\]
here, $[\mm] \cq \bra{1, ..., \mm}$ for $\mm \in \mbn$.
We refer to an element of $\mfm_\nn$ as an \textit{$\nn$-perfect matching}.
% (\textit{$\nn$-\PM}).
Note the double-braces: the union is a set of $\nn$ pairs.
That is, an $\nn$-\PM is a collection of $\nn$ disjoint pairs.

\begin{defn}[Perfect Matching Random Walk]
Let $\kk \in [2, \nn] \cap \mbn$.
The \textit{$\kk$-perfect matching random walk} (\textit{$\kk$-\PM \RW}) on $\mfm_\nn$ has discrete-time dynamics, a step of which is described as follows:
\begin{itemize}[noitemsep]
	\item 
	choose $\kk$ matched pairs, say $\cup_{\ell=1}^\kk \bra{ \bra{i_\ell, j_\ell} }$;
	
	\item 
	disassociate the pairs to give $2\kk$ unpaired elements, ie $\cup_{\ell=1}^\kk \bra{i_\ell, j_\ell}$;
	
	\item 
	uniformly re-pair these $2\kk$ elements.
\end{itemize}
That is,
	$\kk$ matched pairs are chosen,
	the matches are broken
and
	a new matching on these $2\kk$ elements is chosen uniformly.
We denote this process by $\MM_{\nn, \kk} \cq (\MM_{\tt, \nn, \kk})_{\tt\ge0}$.
\end{defn}

%\begin{rmkt}[Use of Sub- and Superscripts]\color{blue}
%	%
%To give fair warning to the reader, this paper is fairly notationally heavy.
%The reader should particularly note that superscripts, as well as powers, are used as indices. Typically, superscripts index the size of the underlying \emph{space}, whilst subscripts indicate \emph{time}.
%Eg, the superscript-$\nn$ in $\mfm_\nn$ indicates that the \PMs are on $2\nn$ objects, whilst the subscript-$\tt$ in $\MM_{\tt, \nn, \kk}$ indicates that the process $\MM_{\nn, \kk}$ has been run for time $\tt$. The $\kk$-superscript in $\MM_{\nn, \kk}$ indicates that the size of the \PM re-pairing is of size $2\kk$.
%	%
%\end{rmkt}

The \PM \RW was first introduced by \textcite{DH:matchings-trees,DH:rw-matchings} in the case $\kk = 2$.
They originally introduced a \RW on \textit{phylogenetic trees}, a biological concept object, along with a bijection between these leaf-labelled trees and \PMs.
Our extension to larger $\kk$ allows more general \RWs on phylogenetic trees to be studied.
More on the biological and other motivations, including randomised algorithms and coding theory, can be found in \cite{DH:matchings-trees,DH:rw-matchings}.
The $2$-\PM \RW was later studied in the representation theory community, where it is known as the \textit{party model}.
%\emph{Cutoff} is studied in these~papers.

We give a full discussion on related work in \S\ref{sec:intro:previous-work}.
Prior to this, we give
	a few brief remarks below,
then give
	precise mixing definitions in \S\ref{sec:intro:mix-def}
and
	state the main theorem in \S\ref{sec:intro:thm}.

\medskip

The dynamics of the \PM \RW are clearly transitive for the space $\mfm_\nn$.
Thus, we may assume~that
\[
	\MM_{0, \nn, \kk}
=
	\brb{ \bra{1,2}, \bra{3,4} , ..., \bra{2\nn-1, 2\nn} },
\]
without loss of generality.
We refer to this \PM as the `identity' matching and denote
\[
	\id_\ell
\cq
	\brb{ \bra{1,2}, \bra{3,4} , ..., \bra{2\ell-1, 2\ell} }
\in
	\mfm_\ell
\Qfor
	\ell \in \mbn.
\]
The dynamics are irreducible.
Thus, another consequence of the transitivity is that the unique invariant distribution of the dynamics, which we denote $\pi_{\mfm_\nn}$, is uniform on $\mfm_\nn$, ie $\pi_{\mfm_\nn} = \Unif(\mfm_\nn)$.

%\begin{rmkt}
%\label{rmk:intro:set-up:re-pairing}
	%
The re-pairing process involves choosing a new matching on a $2\kk$-size subset of $[2\nn]$ and leaving the remainder fixed.
By transitivity, it suffices to be able to sample a \PM with at most $\kk$ `non-fixed points' \wrt $\id_\nn$, ie pairs $\bra{ \eta_{2\ell-1}, \eta_{2\ell} } \ne \bra{ 2\ell-1, 2\ell }$ for $\ell \in [\nn]$.
Indeed, given a general \PM, one first applies a permutation to the labels to send it to the identity matching $\id_\nn$, then replaces this with the sampled \PM. Finally, the inverse of the original permutation is applied.
	%
%\end{rmkt}

%Define $\mfm'$ to be the set of \PMs on $[2\nn]$ with at least $\kk$ non-fixed points:
%\[
%	\mfm'
%&
%\cq
%	\brb{ \eta \in \mfm_\nn \midb \abs{\bra{ \ell \in [\nn] \mid \bra{ \eta_{2\ell-1}, \eta_{2\ell} } \ne \bra{ 2\ell-1, 2\ell } } } \le \kk }
%%\\&
%%=
%%	\brb{ \eta \in \mfm_\nn \midb \abs{\bra{ \ell \in [n] \mid \bra{ \eta_{2\ell-1}, \eta_{2\ell} } = \bra{ 2\ell-1, 2\ell } } } \ge n - k }.
%\]

\medskip

We may drop subscripts, defaulting to $\nn$ or $(\nn, \kk)$, as appropriate.
Eg, $\mfm = \mfm_\nn$ but $\MM = \MM_{\nn, \kk}$.
%Herein, we drop the $(\nn, \kk)$-subscript from $\MM_{\nn, \kk}$, writing $\MM \cq \MM_{\nn, \kk}$.
%We typically only consider \PMs on either $2\nn$ or $2\kk$ objects.
%Nevertheless, we tend to keep the subscripts for clarity;
%we sometimes abbreviate $\mfm \cq \mfm_\nn$.
We also use some abbreviations for frequently-occurring words or phrases.
There are some usual ones:
	``\RW'', ``\TV'', ``\uar'', ``\wrt'' and ``\whp''
abbreviating
	``random walk'', ``total variation'', ``uniformly at random'', ``with respect to'' and ``with high probability'',
respectively.
We also use
	``\CS'' for ``cycle structure''
and
	``\PM'' for ``perfect matching'';
these are not standard, but together they appear close to 300 times throughout the paper, so we feel their abbreviation is legitimate.
%We rarely repeat these.

We do introduce other abbreviations throughout the paper, but only use these `locally'---shortly after their definition. The reader is not expected to remember such abbreviations for more than a couple of paragraphs.
Contrastingly, the reader should remember \PM and \CS throughout~the~paper.

\subsection{Mixing and Cutoff Definitions}
\label{sec:intro:mix-def}

Let $X = (X_\tt)_{\tt \ge 0}$ be an ergodic Markov chain on a finite state space $\Omega$.
Write $\pi$ for its unique invariant distribution.
We are interested in the distance between the law of $X_\tt$ and $\pi$.
%We measure \textit{distance} through \textit{total variation}.

\begin{defn}[Total Variation Distance]
	Let $\mu$ and $\nu$ be probability distributions on $\Omega$.
	The \textit{total variation} (\textit{\textsf{TV}}) distance between $\mu$ and $\nu$ is defined to be
	\[
		\tvt{\mu - \nu}
	\cq
		\SUP{A \subseteq \Omega}
		\absb{ \mu(A) - \pi(A) }.
	\]
	This is known to be equivalent to half the $\ell_1$ distance; see, eg, \cite[Proposition~4.2]{LPW:markov-mixing}.
\end{defn}

The \textit{mixing time} is the time $t$ at which the law of $X_\tt$ is close to $\pi$ in \TV.

\begin{defn}[Mixing Time]
	Define the \textit{mixing time} $\tmix(\cdot)$ by
	\[
		\tmix(\eps)
	\cq
		\inf\brb{
			\tt \ge 0
		\mid
			d(\tt) \cq \MAX{x \in \Omega} \, \tv{ \pr[x]{ X_\tt \in \cdot } - \pi } \le \eps
		}
	\Qfor
		\eps \in (0, 1).
	\]
\end{defn}

We are interested in a sequence $(X_\nn)_{\nn \in \mbn} = ( (X_{\tt, \nn})_{\tt \ge 0} )_{\nn \in \mbn}$ of finite, ergodic Markov chains.
Write
	$\Omega_\nn$ for the state space,
	$\pi_\nn$ for the unique invariant distribution,
	$d_\nn(\cdot)$ for the worst-case \TV distance
and
	$\tmix_\nn(\cdot)$ for the mixing time
of the $\nn$-th chain $X^\nn$.
We want to determine the asymptotic behaviour of $\tmix_\nn(\eps)$ \asinf \nn for each fixed $\eps \in (0, 1)$.

In some special cases, the leading order term of $\tmix_\nn(\eps)$ \asinf \nn does not depend on $\eps$.
This is known as \textit{cutoff}.
It is conjectured to hold for many natural sequences of Markov chains.

\begin{defn}[Cutoff]
	A sequence $(X_\nn)_{\nn \in \mbn}$ of finite, ergodic Markov chains exhibits \textit{cutoff} if
	\[
		\LIMSUP{\toinf \nn}
%		\frac{\tmix_\nn(\eps)}{\tmix_\nn(1-\eps)}
		\tmix_\nn(\eps) \big/ \tmix_\nn(1-\eps)
	=
		1
	\Qforall
		\eps \in (0,1).
	\]
	An equivalent definition is that there exists a sequence $(\tt_\nn)_{\nn \in \mbn}$ of times such that
	\[
		\LIMINF{\toinf \nn}
		d^\nn\rbb{ (1 - \delta) \tt_\nn }
	=
		1
	\Qand
		\LIMSUP{\toinf \nn}
		d^\nn\rbb{ (1 + \delta) \tt_\nn }
	=
		0
	\Qforall
		\delta \in (0,1).
	\]
	The sequence is then said to exhibit \textit{cutoff at time $(\tt_\nn)_{\nn \in \mbn}$}.
\end{defn}

\subsection{Main Theorem}
\label{sec:intro:thm}

Our main results is that the $\kk$-\PM \RW exhibits cutoff if $2 \le \kk \ll \nn$.
Further, we find the leading order of the cutoff time; if \toinf \kk \asinf \nn, then the leading order is given by $\tfrac \nn\kk \log \nn$.

\begin{thm}[Cutoff for the \PM \RW]
\label{res:intro:main}
Let $\nn \in \mbn$.
Let $\kk \in \mbn \setminus \bra{1}$ with \tozero{\kk/\nn} \asinf \nn.
%Let $\MM = (\MM_\tt)_{\tt \ge 0}$ be the $\kk$-\PM \RW on $\mfm$.
Let
\[
	\tt
\cq
	\frac{\nn \log \nn}{\kk - \kk/(2\kk - 1)}
=
	\frac\nn\kk \log \nn
\cdot
	\frac1{1 - 1/(2\kk-1)}.
\]
Then, the $\kk$-\PM \RW $\MM = \MM_{\nn, \kk}$ exhibits cutoff at time $\tt$.
In particular, if \toinf\kk \asinf \nn, then
\[
	\tt
=
	\frac\nn\kk \log \nn
\cdot
	\rbb{ 1 + \oh1 }.
\]
[Officially, this is for a sequence $(\kk_\nn)_{\nn\in\mbn}$ with $\kk_\nn/\nn \to 0$ \asinf \nn and \PM \RWs $(\MM_{\nn, \kk_\nn})_{\nn\in\mbn}$.]
\end{thm}

%\begin{thm}[Cutoff for the Perfect Matching Random Walk]
%\label{res:intro:main}
%	%
%Let $(\kk_\nn)_\ninn \in (\mbn \setminus \bra{1})^\mbn$ satisfy $\kk_\nn / \nn \to 0$ \asinf \nn.
%%	Let $\kk : \nn \to \kk_\nn : \mbn \to \mbn$ satisfy $\kk_\nn / \nn \to 0$ \asinf \nn.
%Let $\MM_{\nn, \kk_\nn} = (\MM_{\tt, \nn, \kk_\nn})_{\tt \ge 0}$ be the $\kk_\nn$-\PM \RW on $\mfm_\nn$ for each $\nn \in \mbn$.
%Define
%\[
%	\tt_\nn
%\cq
%	\frac{\nn \log \nn}{ \kk_\nn - \kk_\nn / (2\kk_\nn - 1) }
%=
%	\frac{\nn}{\kk_\nn} \log \nn
%\cdot
%	\frac{1}{ 1 - 1/(2\kk_\nn - 1) }
%\Qfor
%	\nn \in \mbn.
%\]
%Then, $(\MM_{\nn, \kk_\nn})_{\nn \in \mbn}$ exhibits cutoff at time $(\tt_\nn)_{\nn \in \mbn}$.
%In particular, if $\kk_\nn \to \infty$ \asinf \nn, then
%\[
%	\tt_\nn
%=
%	\frac \nn{\kk_\nn} \log \nn
%\cdot
%	\rbb{ 1 + o_{\toinf \nn}(1) }.
%\]
%	%
%\end{thm}

The upper bound on mixing is the primary focus of this paper.
There are five main steps, analysed in \S\ref{sec:red}--\S\ref{sec:graph}.
These are pulled together to conclude in \S\ref{sec:cutoff:upper}.
The lower bound in a straightforward coupon-collector argument.
This is done in \S\ref{sec:cutoff:lower}, with various details omitted.

\subsection{Related Previous Work}
\label{sec:intro:previous-work}

Establishing cutoff for $\kk = 2$ has received attention in the past, but we are the first to study $\kk > 2$. Not only this, but we handle any $\kk \ll \nn$.
The case $\kk = 2$ is known in representation-theoretic literature as the \textit{party model}. The chain is described by the \textit{Gelfand pair} $(\mfs_{2\nn}, \: \mfs_2 \wr \mfs_\nn)$.
A very readable introduction to this field can be found in the book of \textcite{CsST:harmonic-analysis-finite-groups}.
The $\kk$-\PM \RW with $\kk = 2$ is covered in \cite[\S 11]{CsST:harmonic-analysis-finite-groups}; see also \cite[\S 8]{CsST:gelfand-applications}.
It appears that this approach may extend beyond $\kk = 2$, but it would involve complicated estimation of
\(
	(2\kk-1)!!
=
	(2\kk-1) \cdot (2\kk-3) \cdots 3 \cdot 1
\)
eigenvalues.
The complexity of the transition matrix thus rapidly gets out of hand as $\kk$ grows. It is perhaps only feasible for small $\kk$.

\textcite{DH:rw-matchings} are able to avoid the theory of Gelfand pairs by directly relating the transition matrix for the \PM \RW to that used in the random transpositions (\textit{\RT}) card shuffle; see \cite[Proposition~1]{DH:rw-matchings}.
This only applies for $\kk = 2$.
It becomes a mixture of transition matrices of different card shuffles for larger $\kk$.
Such matrices are not necessarily jointly diagonalisable, rendering their approach much trickier for larger $\kk$.
Our approach can be seen as relating the \PM \RW to another type of card shuffle---namely, conjugacy-invariant \RWs on the permutation group.
There are a significant number of challenges in our own approach which arise only when $\kk > 2$.

A Markov chain on the permutation group $\mfs_\nn$ is a \textit{conjugacy-invariant} \RW if some conjugacy class $\Gamma \subseteq \mfs_\nn$ generates the walk:
	a step involves choosing $\sigma \sim \Unif(\Gamma)$ and composing the current location (permutation) with $\sigma$.
A by-now standard approach to analysing card shuffles corresponding to conjugacy-invariant \RWs, such as the \RT shuffle or $\kk$-cycle \RW, is to project
	from the \RW to its conjugacy class,
then
	from the conjugacy class to the corresponding integer partition
and finally
	use a variant of a coupling due to \textcite{S:random-trans};
see, eg, \cite{BSZ:k-cycle,BS:cutoff-conj-inv,B:random-trans-coupling}.
A similar projection for the \PM \RW can be defined; see, eg, \cite{BKLM:interchange-rev,GUW:quantum-heisenberg}.
We use this.

\medskip

The previous work to which our approach is most related is that of \bst, which builds in part on work of \bszt.
Cutoff for the $\kk$-cycle \RW for any fixed $\kk$, not depending on $\nn$, is established in \bszc.
They strongly believe that their argument can be extended to consider general conjugacy classes $\Gamma$ with bounded support:
	$\kk = \abs \Gamma$, independent of $\nn$.
This is extended to allow any $\kk = \abs \Gamma \ll \nn$ in \bsc, similarly to what we allow.

Parts of our argument are very similar to those of \bsc.
We feel that it is important to detail which ideas are our own and which parts are adjustments or extensions of \bsc.
We do so in
%We detail this in
\S\ref{sec:outline-compare-bs}.

\medskip

Cutoff for the $\kk$-cycle \RW with $2 \le \kk \ll \nn$ has also been established independently by \textcite{H:cutoff-k-cycle}.
He introduces and uses an asymptotic estimation of the characters of $\mfs_\nn$ evaluated at cycles.
\citeauthor{H:cutoff-k-cycle} tentatively suggests that his method can be extended to some conjugacy classes. However, he believes that some new ideas are needed to obtain the full result of \bsc.

The \emph{limit profile} was determined for the $\kk$-cycle case by the current author and Nestoridi \cite{NOt:limit-profiles}, building on work of \textcite{H:cutoff-k-cycle,T:limit-profile}.
Extending this argument to general conjugacy classes, even of bounded support, appears to be very technically challenging.

\medskip

\begin{figure}
	\centering
%\noindent%
\begin{minipage}{0.9\textwidth}
	\includegraphics[width = 1\textwidth]{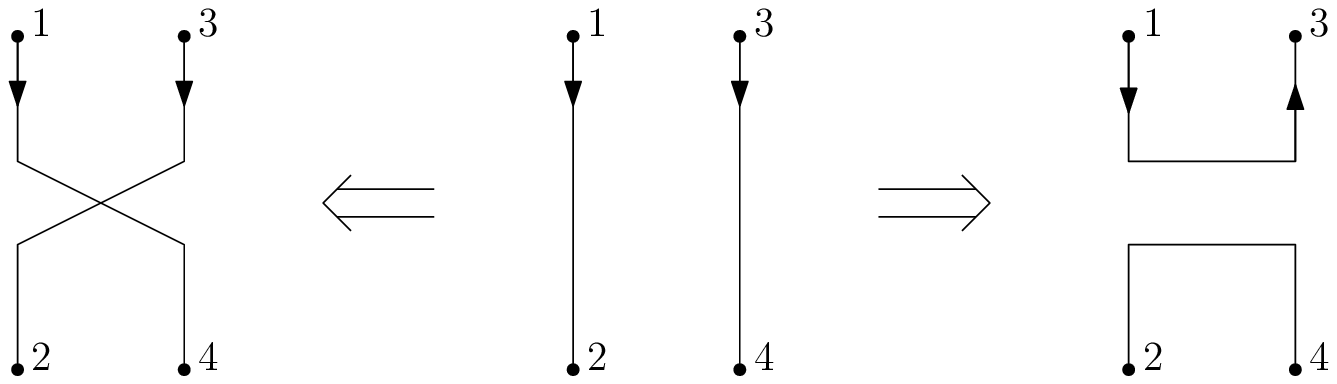}
	
	\caption{%
		The rematchings have the following correspondences:
		\begin{alignat*}{2}
			&\bcdot\ \text{the `cross' (left)}&
		\quad
			\brb{ \bra{1,2},\: \bra{3,4} }
		\ &\longrightarrow{} \
			\brb{ \bra{1,4},\: \bra{2,3} }
		\Quad{to a}
			\text{`transposition';}
		\\
			&\bcdot\ \text{the `bar' (right)}&
		\quad
			\brb{ \bra{1,2},\: \bra{3,4} }
		\ &\longrightarrow{} \
			\brb{ \bra{1,3},\: \bra{2,4} }
		\Quad{to a}
			\text{`reversal'.}
		\qedhere
		\end{alignat*}
		\vspace*{-7mm}
	}
	\label{fig:trans_rev}
\end{minipage}
\end{figure}

The other work most related to ours is by \bklmt on an interchange process with reversals.
It is related to stochastic representations of quantum spin systems, namely anti/ferromagnetic Heisenberg models.
Very roughly, the ferromagnetic model has interactions between spins which behave like `transpositions'; antiferromagnetic models additionally have `reversals'.
\cref{fig:trans_rev} shows the two possible $\PM$ rematchings,
	known as `transpositions' and `reversals',
when $\kk = 2$.
%these are known as `transpositions' and `reversals'.
This shows the corresponds between the $2$-\PM \RW and Heisenberg models.

The $2$-\PM \RW is briefly discussed by \textcite{CLR:aldous-spectral-gap}. This is the paper in which Aldous's famous spectral gap conjecture is proved.
The conjecture---now a theorem---regards the interchange processes (\textit{\IP}), which is a generalisation of the \RT shuffle.
%, and the \RW on a graph.
They show that the spectral gap $\lambda_1^\IP$ of the \IP equals the spectral gap $\lambda_1^\RW$ of the \RW.
The $2$-\PM \RW is a projection of the \IP, in a precise sense, which implies that its spectral gap \(\lambda_1^\text{$2$-\PM}\) satisfies
\[
	\lambda_1^\RW
=
	\lambda_1^\IP
\le
	\lambda_1^\text{$2$-\PM}.
\]
\textcite{DH:rw-matchings} completely characterise the spectrum of the $2$-\PM \RW.
In particular, their work shows that this inequality is not tight when $\nn = 2$, ie when considering \PMs on $4$ objects.
We defer the reader to \cite[\S 4.2.2]{CLR:aldous-spectral-gap} for further details.

\medskip

The history of the $2$-\PM \RW goes back to \textcite{DH:rw-matchings}, being studied later by \textcite{CsST:gelfand-applications,CsST:harmonic-analysis-finite-groups,CLR:aldous-spectral-gap}, as detailed above.
The $\kk$-\PM \RW appeared recently in work of
Avena, Güldaş, van der Hofstad and den Hollander~%
\cite{AGHH:dynamic-cm,AGHH:dynamic-cm2},
%\cite{AGHH:dynamic-cm2},
albeit in a slightly different set-up:
	they use the $\kk$-\PM \RW to drive a dynamic graph model via the configuration model.
%	 (\textit{\CM}).

The \textit{configuration model} (\textit{\CM}),
	 introduced in different forms by \textcite{BC:cm,B:cm,B:random-graphs},
randomly samples a graph with a given degree sequence as follows.
\begin{itemize}[noitemsep]
	\item 
	Let $d = (d_1, ..., d_m) \in \mbn^m$.
	Assume that $n \cq \tfrac12(d_1 + \cdots d_m) \in \mbn$.
	
	\item 
	Place $m$ vertices and attach $2n$ `half edges' to the vertices:
		$d_v$ to vertex $v$ for each $v \in V$.
	
	\item 
	Uniformly pair the half edges to create a graph on $m$ vertices with $n$ edges.
\end{itemize}
The resulting graph has $m$ vertices, $n$ edges and degree sequence $d$.
An excellent introduction to the configuration model, with multiple explanatory figures,
is given by van der Hofstad~\cite[\S 7]{H:random-graphs-book}.
%can be found in \cite[\S 7]{H:random-graphs-book}.
%See Figure~\ref{fig:cm} for an example.

The above viewpoint is of the \CM as a projection of a \PM.
The $\kk$-\PM \RW induces a dynamic random graph process:
	simply select $\kk$ edges,
	cut them to produce $2\kk$ half-edges
and
	randomly re-pair the half-edges.
This drives a dynamic \CM by keeping the half-edges attached to the same vertices throughout.
A non-backtracking \RW is placed on this in
\cite{AGHH:dynamic-cm,AGHH:dynamic-cm2}
%\cite{AGHH:dynamic-cm2}
and its mixing properties are studied.
Focus is on properties of the walk. The graph is not studied in detail.
In particular, the question of the mixing time of the dynamic random graph is left open.
%It is from their work that this current project arose.

\medskip

The \CM need not be a simple graph.
It is known that
	the probability of being simple is bounded
	under some regularity conditions for vertex degrees;
see \cite[\S 7.4]{H:random-graphs-book} for precise details.
Importantly, the law of the \CM \emph{conditioned on being simple} is uniform over all simple graphs, with the appropriate degree sequence.
Analogously, the dynamic \CM does not consist only of simple graphs.

The \textit{switch chain} is defined to be the dynamic \CM, but where transitions are rejected if they give rise to a non-simple graph.
The invariant distribution of this chain is uniform over all simple graphs with the appropriate degree sequence.
The purpose of the switch chain is to draw from such graphs uniformly at random.
Control on the mixing time is naturally required for such sampling.

Analysis of the switch chain has a long and rich history; far too much to discuss in totality here.
Recent overviews can be found in the introductions of the recent papers \cite{AK:switch-chain:conf,AK:switch-chain:jour}, \cite{EGMMSS:switch-chain:arxiv} or \cite{TY:switch-chain}.
It was introduced by \textcite{KTV:switch-chain:conf,KTV:switch-chain:jour} in the late 90s, making it over 20 years old.
Even so, it is still an extremely active area of research.
%	the paper \cite{GG:switch-chain}, released in 2020, cites seven related papers, released in either 2018 or 2019 alone.

\subsection{Acknowledgements}

This is a single-author project, but I would be remiss not to acknowledge the input of others.

The initial question arose out of a research visit of mine to the EURANDOM group in the Netherlands during my PhD, in early 2019.
I discussed this question extensively with Güldaş, as well as Avena, van der Hofstad and den Hollander.
We also discussed their work
\cite{AGHH:dynamic-cm,AGHH:dynamic-cm2}
%\cite{AGHH:dynamic-cm2}
on the dynamic \CM.
I gratefully acknowledge their insights and comments at the start of this project, as well as the hospitality of the EURANDOM group more widely.
I met with Berestycki a few months later to discuss his work \bsc and how it might adapted to my set-up.

I subsequently started building on these ideas at the end of my PhD. I had a number of very fruitful discussions with my then PhD supervisor Perla Sousi around this time, in late 2020. I gratefully acknowledge her ideas, comments and contributions.

The helpful comments provided by the anonymous reviewer significantly improved the presentation and clarity of this paper.
%particularly in a lightening of the technical weight.
They have my sincere thanks.
%---and the readers' too, no doubt.

%\section[Outline of Approach and Comparison with Berestycki and Şengül \cite{BS:cutoff-conj-inv}]{Outline of Approach and Comparison with \cite{BS:cutoff-conj-inv}}
\section{Outline of Approach and Comparison with \cite{BS:cutoff-conj-inv}}
\label{sec:outline-compare-bs}

This section first outlines the underlying approach.
There are five main steps for the upper bound on mixing, which we detail below.
The lower bound is much more straightforward, via a coupon-collector argument.
Establishing the upper bound is the primary focus of the article.

To close the section, we compare and contrast the methods used in the current article with those of developed by \bst.
Related comments are made throughout the paper.
% when the relevant concepts are introduced.
%There is significant repetition between the current section and those comments.
We feel that it is important to be transparent regarding the similarities and differences between our work and theirs,
so we include this summary to gather all relevant remarks are together.

\vspace*{-\smallskipamount}

\paragraph*{\S\ref{sec:red}: Projecting to Cycle Structure and Partitions.}

The first step involves projecting the \PM \RW to its \textit{cycle structure}---a concept that we introduce below, akin to that for permutations---and then further to its corresponding partition.
This idea has become a standard approach when analysis conjugacy-invariant \RWs on groups by now, being mentioned at least as early as \cite{DH:matchings-trees}; it is used by \cite{DH:matchings-trees,S:random-trans,BSZ:k-cycle,B:random-trans-coupling,BS:cutoff-conj-inv} and surely many more.
The current article extends the $2$-\PM \RW of \cite{DH:matchings-trees} to the general $\kk$-\PM \RW via a decomposition of a $\kk$-\PM into $\kk-1$ $2$-\PMs.

\vspace*{-\smallskipamount}

\paragraph*{\S\ref{sec:decomp}: Decomposing into Swaps.}

It is natural to break down permutations into products of transpositions.
%\(
%	(a_1, ..., a_\ell) = (a_1, a_2) (a_2, a_3) \cdots (a_{\ell-1}, a_\ell).
%\)
Doing so permits analysis of the above partition walk via an adjustment of a coupling due to \textcite{S:random-trans}; this was done first by \textcite{BSZ:k-cycle} and used in \bsc.
Constructing an analogous decomposition of a $\kk$-\PM into a sequence of $\kk-1$ $2$-\PMs (`swaps') is one of the most fundamental parts of the paper, permitting analysis~of~the~partition~walk.

%It is easy to decompose an $\ell$-cycle into $\ell-1$ $2$-cycles (transpositions):
The decomposition is easy for permutations:
	\(
		(a_1, ..., a_\ell) = (a_1, a_2) (a_2, a_3) \cdots (a_{\ell-1}, a_\ell).
	\)
Unfortunately, such a natural idea does not apply for \PMs:
	after pairs labelled $(a_1, a_2)$ are swapped, there is no way of identifying a specific one of the new pairs with $a_1$ and the other with $a_2$;
	see \cref{fig:destroy_label}.
One of the main inventions of the whole paper is an algorithmic approach for drawing a uniform `single-cycle' $\ell$-\PM, roughly corresponding to an $\ell$-cycle in permutation language, from $\ell-1$ swaps.
%	 See \cref{alg:decomp:gen:gen-cycle} for the details.
We are able to use the principle behind this approach in a variety of other scenarios.

\vspace*{-\smallskipamount}

\paragraph*{\S\ref{sec:cf}: Analysis of Partition Walk.}

A partition can be seen as a tiling of $(0,1]$ by rescaling.
A single \textit{step} of the partition walk when $\kk = 2$ involves choosing markers $u, v \in \bra{1/\nn, ..., \nn/\nn}$ uniformly:
	if $u$ and $v$ are in different blocks, then the blocks are merged;
	if they are in the same block, then the block is split according to some simple rule.
When $\kk > 2$, a single \textit{round} is broken into $\kk-1$ \emph{steps} via the above decomposition into swaps. The markers are no longer uniform:
	one corresponds to a marker from a previous step and the other is uniform on what is yet to be chosen.

If $\kk \ll \sqrt \nn$, then this sampling \emph{without} replacement can be well-approximated by sampling \emph{with} replacement, somewhat decoupling the $\kk-1$ steps.
This idea originated in \bszc where $\kk \asymp 1$.
\bst realised that if $\kk \ll \nn$, then each individual draw is still a uniform draw from a collection of at least $\nn-\kk \approx \nn$ objects.
So, \emph{marginally} the two processes are similar, even conditional on what has come before.
They show that this is sufficient.
We use the same idea, combined with our new decomposition algorithm discussed immediately above.

A variant of the coupling of \textcite{S:random-trans} is introduced in \bszc and used almost unchanged in \bsc.
The variation from \cite{S:random-trans} is only fairly minor, but is crucial to make sure that blocks in the partition do not become unmanageably small.
Adjusting the coupling from conjugacy-invariant \RWs to \PM \RWs is not trivial.
%A very similar argument for $2$-\PM \RW using \citeauthor{S:random-trans}'s original coupling is given in \bklmc.
The first marker in one swap is the same as the second marker in the previous swap in \cite{BSZ:k-cycle,BS:cutoff-conj-inv}.
Such a statement cannot hold for the \PM \RW, because of the lack of identifiability of the previous matched discussed above.
%This difficulty is very similar to that faced in \S\ref{sec:decomp} and
It makes defining and controlling the coupling of \cite{S:random-trans} in our case more challenging compared with in \cite{BSZ:k-cycle,BS:cutoff-conj-inv}.
The subtlety does not arise in \bklmc where $\kk = 2$, as there swaps are completely independent.
%It thus does not arise in \bklmc, where a related argument is used for $2$-\PM \RWs.
%There is also a concept of reversing a cycle rather than splitting it; see \cref{fig:reject_split}.
%This does not cause any major issues, however.

\vspace*{-\smallskipamount}

\paragraph*{\S\ref{sec:3coup}: Path Coupling Structure.}

We use a path-coupling argument to couple two partition walks.
The structure of this argument originated in \bsc.
%\footnote{%
%	There, they refer to it as ``discrete Ricci curvature'',
%		introduced by \textcite{O:ricci-curvature}.
%	They remark that it is the same concept as the famous ``path coupling'',
%		introduced by \textcite{BD:path-coupling}}
The application to conjugacy-invariant \RWs was new in \bsc,
	in particular inspecting the relative distance of the two walks after order $\beta \nn / \kk$ steps and letting $\beta \to \infty$.
This is markedly different to the usual inspection after just a single step.
The justification for this time $\beta \nn / \kk$ is outlined in the next part.
Their general path-coupling structure requires only very minor adjustment to apply in our set-up.

\vspace*{-\smallskipamount}

\paragraph*{\S\ref{sec:graph}: Auxiliary Graph Process.}

%In short, originality of ideas in this section is not claimed.

The mixing time of random transpositions ($2$-cycles) on $\nn$ cards is order $\nn \log \nn$.
However, \textcite{S:random-trans} showed that the approximate structure of the \emph{large} cycles relax to uniformity in time order $\nn$.
His proof goes via an auxiliary graph process:
	$i$ and $j$ are connected at time $t$ if transposition $(i,j)$ has been applied by this time.
This gives \emph{precisely} the usual Erd\H{o}s--R\'enyi graph.
If $\gamma n$ transpositions are applied with $\gamma > \tfrac12$, then the graph has a giant component \whp. This giant component is key in analysing the structure of the large cycles.

The same idea is used in both \cite{BSZ:k-cycle,BS:cutoff-conj-inv}.
	The former restrict to \CSs with support $\kk \asymp 1$, whilst the latter allow any $\kk \ll \nn$.
%	---although they comment that they believe it can be extended, at least to any $\kk \ll \sqrt \nn$.
Both generalise \citeauthor{S:random-trans}'s construction to a hyper-graph:
	hyper-edge $\bra{a_1, ..., a_\ell}$ is added if the cycle $(a_1, ..., a_\ell)$ is applied as part of the cycle decomposition of the element of $\Gamma$.
The analysis in \bsc is far more tricky than in \cite{BSZ:k-cycle}.

We use exactly the same ideas, once we have the correct viewpoint relating an $\nn$-\PM to an $\nn$-permutation.
In particular, we \emph{do not} consider a graph on $2\nn$ objects and connect objects $x$ and $y$ if they are matched at some point.
There is a slight difference in our set-up:
	\bsc considers a \emph{fixed} \CS $\Gamma$, while our \CS changes from step to step.
Some non-trivial adjustments are needed.

\vspace*{-\smallskipamount}

\paragraph*{\S\ref{sec:cutoff:upper}: Combining Results.}

All the above is developed for the \emph{upper} bound on mixing.
From these, particular the path-coupling bounds, concluding an upper bound on the mixing time is not difficult.
%These results are pulled together in an analogous way to in \bsc, implementing the path coupling bounds.

\vspace*{-\smallskipamount}

\paragraph*{\S\ref{sec:cutoff:lower}: Lower Bound.}

The \emph{lower} bound is a standard coupon-collector argument, as in \cite{BSZ:k-cycle,BS:cutoff-conj-inv}:
	the number of fixed points is used as a distinguishing statistic.
The idea is not new to \cite{BSZ:k-cycle,BS:cutoff-conj-inv}; rather, it has been in a variety of related papers in the past.
It is so standard that it is deferred to the appendix of \bsc and omitted completely from \bszc.
%We omit the majority of the details.

\vspace*{-\smallskipamount}

\paragraph*{Comparison with \bsc.}

%To close this section, we compare and contrast the methods used in the current article with those of \bsc.
%Related comments are made throughout the paper.
%% when the relevant concepts are introduced.
%%There is significant repetition between the current section and those comments.
%It is important to be transparent regarding the similarities and differences between our work and theirs,
%so we include this summary to gather all relevant remarks are together.

%To summarise, the current article is strongly inspired by the methods \bsc.
Our article is strongly inspired by the methods of \bsc.
However, a significant number of new ideas are required.
It is not even clear a priori that the \PM \RW can be related so closely to a conjugacy-invariant \RW.
The natural approach of viewing an $\nn$-\PM as a permutation on its $2\nn$ objects does not allow this: the corresponding Cayley graph is not generated by a union conjugacy class.
It can be viewed as a \emph{Gelfand pair} $(\mfs_{2\nn}, \mfs_2 \wr \mfs_\nn)$,
%as discussed before,
but this moves far from the probabilistic approach of \bsc towards the representation-theoretic of \cite{CsST:gelfand-applications,CsST:harmonic-analysis-finite-groups}.

Key is to introduce the `cycle structure' of a \PM via cycle lengths in an induced graph.
We use this to relate an $\nn$-\PM to an $\nn$-permutation.
We then adjust the techniques developed for conjugacy-invariant \RWs to this \PM viewpoint.
Some of these adjustments are trivial, but many are far from easy.
Multiple subtleties arise for the \PM \RW which are not present for conjugacy-invariant \RWs.
%---even just drawing a \PM \uar via a sequence of swaps is difficult.
%---not least of which is just simply decomposing a uniform $\kk$-\PM into a correlated sequence of swaps.

We feel that the merit and contribution of this paper is not in the technical proficiency of the argument, but rather in developing the correct viewpoint.
The underlying ideas can be found in \bsc; being able to utilise them for \PMs is the challenge.
For example, the lack of consistency in the labelling of previously interacted with objects, discussed at length throughout the paper, is a constant source of difficulty:
	in decomposing a \PM into a sequence of swaps in \S\ref{sec:decomp};
	in the tiling of \S\ref{sec:cf};
	in the construction of the auxiliary graph process in \S\ref{sec:graph}.
Overcoming such challenges is paramount.

\section{Reductions and Realisations}
\label{sec:red}

\subsection{Cycle Definitions and Reduction to Cycle Structure}
\label{sec:red:cs}

%All definitions of `cycles' will be \wrt the appropriate identity matching.

Let $\eta \in \mfm$ be a \PM.
Consider the graph $\rbr{ [2\nn], \eta }$, which has $2\nn$ vertices and edge-set given by the pairs in the \PM $\eta$.
Every vertex is of degree $1$:
	it is paired with one other vertex.
Now consider the union of this graph with $\rbr{ [2\nn], \id }$, as a multigraph, written $\rbr{ [2\nn], \eta \cup \id }$.
	Every vertex is now of degree $2$ and each vertex is contained in a unique cycle of even length---this is counting an isolated double-edge between two vertices as a cycle of length $2$.
We define the \textit{cycle structure} of $\eta$ as the vector giving the number of $2$-cycles, $4$-cycles and so on; see \cref{def:red:cs} for the formal definition.

We sometimes consider arbitrary \PMs, not specifying the number of underlying objects:
write
\[
	\mfm_\infty
\cq
	\cup_{\nn=1}^\infty \mfm_\nn.
\]
This is a disjoint union and $\eta \in \mfm_\nn$ implies that $\abs \eta = n$.
Here, $\eta$ is a set whose elements are disjoint, unordered pairs; so $\abs \eta = n$ means that there are $\nn$ disjoint, unordered pairs in $\eta$.
We emphasise that every element of $\mfm_\infty$ is a \emph{finite} set, ie corresponds to a \PM on $2\nn$ objects~for~some~$\nn \in \mbn$.

\begin{defn}[Cycle Structure]
\label{def:red:cs}
The \textit{cycle structure} (\textit{\CS}) of a \PM $\eta \in \mfm_\infty$ is the vector
\[
	\mcc(\eta)
\cq
	\rbb{ \mcc_1(\eta), \mcc_2(\eta), ... }
\in
	\mbn_0^\mbn,
\]
where $\mcc_\ell(\eta)$ is the number of $2\ell$-cycles in the multigraph $\rbr{ [2\abs \eta], \eta \cup \id_{\abs \eta} }$ for $\ell \in \mbn$.
Write
\[
	\mfc_\nn
\cq
	\brb{ \cc \in \mbn_0^\mbn \midb \sumt[\infty]{\ell=1} \ell \cc_\ell = \nn }
\Qfor
	\nn \in \mbn;
\]
this is the set of possible \CSs for an element of $\mfm_\nn$, ie an $\nn$-\PM.
Write
\[
	\mfc_\infty
\cq
	\cup_{\nn=1}^\infty \mfc_\nn.
\]
\end{defn}

Some examples are given in \cref{fig:exm}.
These figures correspond to \cref{exm:decomp:dist-supp} in which other statistics---namely, \textit{support} and \textit{swap distance}---are analysed.
We do not repeat the image here; rather, the reader can look ahead to Page \pageref{fig:exm} to see the graphs corresponding to different \PMs $\eta$.

%\(
%	\mfc_\infty
%\ne
%	\mbn_0^\mbn,
%\)
%as $\sumt[\infty]{\ell=1} \ell \cc_\ell < \infty$ for all $\cc \in \mfc_\infty$.

\begin{term}[Cycle Terminology]
We drop the trailing $0$s when writing out a \CS, eg writing $(1, 2)$ rather than $(1, 2, 0, 0, ...)$.
A non/fixed point is assumed to be \wrt the appropriate identity---$\id_{\abs \eta}$ in the above case---if the \PM to which it is to be compared is omitted.

We refer to $\mcc_\ell(\eta)$ as the number of \textit{$\ell$-cycles} in the \PM $\eta \in \mfm_\infty$ for $\ell \in \mbn$, not \textit{$2\ell$-cycles}.
Then, $1$-cycles in the \PM, ie $2$-cycles in the multigraph, are precisely fixed points.
Thus, there are $\mcc_1(\eta)$ fixed points in the \PM $\eta \in \mfm_\infty$;
analogously, we say that there are $\cc_1$ fixed points in the~\CS~$\cc \in \mfc_\infty$.
\end{term}

This choice of terminology, specifically the $\ell$-cycle vs $2\ell$-cycle distinction, may seem somewhat peculiar at first. However, we shall soon see that it makes the definitions for $\nn$-\PMs analogous to the standard definitions for $\nn$-permutations, ie permutations on $\nn$ objects.

We use this terminology to set up the following notation.
Recall that $\mfm = \mfm_\nn$ and $\mfc = \mfc_\nn$.

\begin{nota}
Write $\mfm' \subseteq \mfm$, respectively $\mfc' \subseteq \mfc$, for those $\eta \in \mfm$, respectively $\cc \in \mfc$, with at most $\kk$ non-fixed points \wrt the identity $\id$.
Write $\mfm_\infty(\cc) \subseteq \mfm_\infty$ for those $\eta \in \mfm_\infty$ with \CS $\mcc(\eta) = \cc \in \mfc_\infty$.
In mathematical notation, make the following definitions for $\cc \in \mfc_\infty$:
\begin{gather*}
	\mfm_\infty(\cc) \cq \brb{ \eta \in \mfm_\infty \midb \mcc(\eta) = \cc },
\quad
	\mfc' \cq \brb{ \cc \in \mfc_\nn \midb \cc_1 \ge \nn - \kk };
\\
	\mfm'
\cq
	\brb{ \eta \in \mfm_\nn \midb \mcc(\eta) \in \mfc' }
=
	\brb{ \eta \in \mfm_\nn \midb \mcc_1(\eta) \ge \nn - \kk }
=
	\cup_{\cc \in \mfc'} \mfm_\infty(\cc).
\end{gather*}
We emphasise that this is an important definition which the reader should commit to memory:
	adding a prime ($\prime$) to $\mfm$ or $\mfc$, giving $\mfm'$ or $\mfc'$, indicates that there are at most $\kk$ non-fixed points.
%Note that $\mfm_{\nn, \nn} = \mfm_\nn$ and $\mfc_{\nn, \nn} = \mfc_\nn$.
	%
\end{nota}

%\begin{lem}[Uniformity Given Cycle Structure]
%\label{res:red:cs:cycle-unif}
%	Assume that $\MM_{\nn, \kk}_0 = \id$.
%	Let $\cc \in \mfc_\nn$.
%	Then
%	\[
%		\mcl\rbb{ \MM_{\tt, \nn, \kk} \midb \mcc(\MM_{\tt, \nn, \kk}) = \cc }
%	=
%		\Unif\rbb{ \mfm_\infty(\cc) },
%	\]
%	under the assumption that $\bra{\mcc(\MM_{\tt, \nn, \kk}) = \cc}$ is a non-null event.
%\end{lem}

The purpose of introducing this \emph{cycle structure} is that the law of the \PM \RW given its \CS is uniform over all \PMs with this given \CS.
A completely analogous projection is often used when studying conjugacy-invariant \RWs on the permutation group.
This means that projecting from $\mfm$ to $\mfc$ does not decrease the \TV distance from equilibrium.
We use $\mcl(\cdot)$ to denote the law of a random variable.
Abbreviate $C_{\tt} \cq \mcc(\MM_{\tt}) \in \mfc$.
Denote by $\pi_{\mfc}$ the invariant distribution
of $C \cq (C_\tt)_{\tt\ge0}$.

\begin{lem}[\TV-Preserving Projection to Cycle Structure]
\label{res:red:cs:proj}
%	For any $\nn \in \mbn$ and $\kk \in [2, \nn] \cap \mbn$,
	The projection of the \PM \RW from the perfect-matching space $\mfm$ to the cycle-structure space $\mfc$ is \TV-preserving:
	\[
		\tv{ \mcl(\MM_{\tt}) - \pi_{\mfm} }
	=
		\tv{ \mcl(C_{\tt}) - \pi_{\mfc} }.
	\]
\end{lem}

\begin{Proof}
The pairs at each round of the \PM \RW are chosen uniformly and independently between rounds.
The uniformity of the \PM \RW given its \CS is thus an immediate consequence of this symmetry.
The \TV-preservation claim follows immediately from this, eg by a trivial coupling.
\end{Proof}

\subsection{From Cycle Structures to Integer Partitions}
\label{sec:red:part}

We have reduced from \PMs to \CSs.
We now explain how to reduce further: from \CSs to partitions.

Recall that all the cycles in the graph $([2\nn], \eta \cup \id)$ have even lengths and are disjoint. They thus form an integer partition of $[2\nn] = \bra{1, ..., 2\nn}$.
The terminology we used divided these lengths by $2$---eg a fixed point, or $1$-cycle, of a \PM or \CS corresponded to a $2$-cycle in the graph.
These halved values are all integers and form an integer partition of $[\nn] = \bra{1, ..., \nn}$.

We can further divide these values by $\nn$ to get a partition of $[0, 1]$ with block lengths in $\bra{0, 1/\nn, ..., 1} = [0, 1] \cap (\tfrac1\nn \mbz)$.
We refer to this latter situation as a $\tfrac1\nn$-integer partition of $[0, 1]$.
We tend to drop $\tfrac1n$-prefactor, including it only when there may be ambiguity.
Write
\[
	\mfp_\nn
\cq
	\brb{
		(x_1, ..., x_\nn) \in \bra{0, 1/\nn, ..., (\nn-1)/\nn, 1}^\mbn
	\midb
		x_1 \ge \cdots \ge x_\nn, \; \sumt[\infty]{\ii=1} x_\ii = 1
	}.
\]
We are not always concerned about the non-increasing order of the \textit{blocks} in the partition; in this case, we write a partition $x \in \mfp_\nn$ as an unordered multiset $[x_1, ..., x_\nn]$.
Abbreviate $\mfp \cq \mfp_\nn$.

We also define the limiting case, which we refer to as a \textit{continuous partition} of $[0, 1]$.
Write
\[
	\mfp_\infty
\cq
	\brb{
		(x_1, x_2, ...) \in [0,1]^\mbn
	\midb
		x_1 \ge x_2 \ge \cdots, \; \sumt[\infty]{\ii=1} x_\ii = 1
	}.
\]
Contrary to our previous notation, $\mfp_\infty \ne \cup_{\nn=1}^\infty \mfp_\nn$:
	all partitions have finitely many blocks
	in the latter.
Further, the entries $x_\ii$ of $x \in \mfp_\infty$ need not be rationals.

\begin{defn}[Coalescence--Fragmentation Chain]
\label{def:red:part:cf-def}
	Define $\mcp(\eta) \in \mfp$ to be the integer partition corresponding to the \PM $\eta \in \mfm$, with blocks in non-increasing order of size.
	Abbreviate $\PP_{\tt} \cq \mcp(\MM_{\tt})$, where $(\MM_{\tt})_{t\ge0} \in \mfm^{\mbn_0}$ is the $\kk$-\PM \RW.
	Write $\pi_{\mfp}$ for the invariant distribution of $\PP \cq (\PP_{\tt})_{\tt\ge0}$.
	The chain $(\PP_{\tt})_{\tt\ge0}$ is a \textit{\cf chain}.
%	It is the only type of such processes that we study here. We thus refer to it as \textit{the \cf chain}.
\end{defn}

We now describe the evolution of the \cf chain when $\kk = 2$.
We then comment on how it differs from the corresponding chain for the random transpositions shuffle.
%Recall that we denote the \cf chain by $\MC P t$ and the partition space
%\[
%	\mfp_\nn
%=
%	\brb{ (x_1, ..., x^\nn) \in \bra{0, 1/\nn, ..., (\nn-1)/\nn, 1}^\mbn \midb x_1 \ge \cdots \ge x^\nn \text{ and } \sumt[\infty]{1} x_i = 1 }.
%\]
%Sometimes we are not concerned about the non-increasing order of the \textit{blocks} in the partition, and we write a partition $x \in \mfp_\nn$ as an unordered multiset $[x_1, ..., x^\nn]$.

\begin{lem}[Evolution of Coalescence--Fragmentation Chain]
Suppose that $\kk = 2$.
Suppose that the \cf chain is at the integer partition $\pp = [\pp_1, ..., \pp_\nn] \in \mfp$.
%Consider first $\kk = 2$.
Suppose that the two pairs/matches chosen are indexed by $i$ and $j$ are chosen, respectively.
There are two cases.

\begin{itemize}
	\item 
	If the two pairs are in different blocks, say $i \in \pp_1$ and $j \in \pp_2$, then these two blocks \textit{merge}.
	The resulting partition is given by $[\pp_1 \cup \pp_2, \pp_3, ..., \pp_n]$.
	
	\item 
	If the two pairs are in the same block, say $i,j \in \pp_1$, then with probability $\tfrac12$ nothing changes and with probability $\tfrac12$ the block \textit{splits}, say into $\pp_1^-$ and $\pp_1^+$.
	If there is a split, then location of the split is uniform.
	The resulting partition is given by $[\pp_1^-, \pp_1^+, \pp_2, ..., \pp_n]$ when there~is~a~split.
\end{itemize}
\end{lem}

The evolution is similar for $k > 2$:
	one generates the $\kk$-rematching via $\kk-1$ swaps
	in a certain manner, described in \cref{alg:decomp:gen:gen-cycle} below,
and
	applies the above cases to each of the $\kk-1$ swaps.

%\begin{Proof}
%	%
%This is immediate from the construction when $\kk = 2$.
%The case $\kk > 2$ follows analogously using \cref{alg:decomp:gen:gen-cycle}:
%	there,
%	the sampling of a $\kk$-\PM is broken down into $\kk-1$ swaps, ie $2$-\PMs.
%	%
%\end{Proof}

\begin{figure}
	\centering
%\noindent%
\begin{minipage}{0.9\textwidth}	
	\includegraphics[width = 1\textwidth]{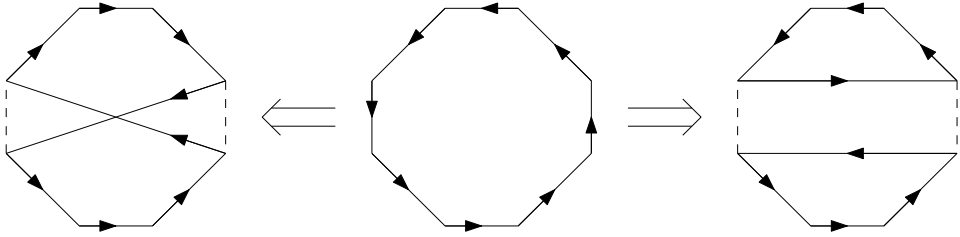}
	
	\caption{%
			The pair of vertical edges is picked and rematched into either a cross (left) or a bar (right). The bar splits the cycle in two, but the cross does not.
		}
	\label{fig:reject_split}
\end{minipage}
\end{figure}

\begin{rmkt}[Comparison with Random Transpositions]
%\label{rmk:red:part:evolution}
	%
The chain corresponding to the random transpositions shuffle is extremely similar.
The only difference is that there is a split \emph{every} time the cards are in the same block there.
The split occurs with probability $\tfrac12$ for our chain corresponding to \PMs.
This is the only difference.
We demonstrate this in \cref{fig:reject_split}.
Splits are rejected half the time in the antiferromagnetic Heisenberg model, studied in \bklmc, too; see \cref{fig:trans_rev}.
\end{rmkt}

The next lemma shows why we introduced the \cf chain.
Analogous results are used in \cite{S:random-trans,B:random-trans-coupling,BSZ:k-cycle,BS:cutoff-conj-inv} for conjugacy-invariant \RWs and in \cite{DH:rw-matchings,BKLM:interchange-rev} for $2$-\PM \RW.
%, amongst others.
Related claims are proved via representation theory, eg in \cite{CsST:gelfand-applications,CsST:harmonic-analysis-finite-groups,H:cutoff-k-cycle}.

\begin{lem}[\TV-Preserving Projection to Coalescence--Fragmentation Chain]
\label{res:red:part:proj}
%	For any $\nn \in \mbn$ and $\kk \in [2, \nn] \cap \mbn$,
	The projection of the $\kk$-\PM \RW from the perfect-matching space $\mfm$ to the integer-partition space $\mfp$ is \TV-preserving:
	\[
		\tv{ \mcl(\MM_{\tt}) - \pi_{\mfm} }
	=
		\tv{ \mcl(C_{\tt}) - \pi_{\mfc} }
	=
		\tv{ \mcl(P_{\tt}) - \pi_{\mfp} }.
	\]
\end{lem}

\begin{Proof}
The first equality is precisely \cref{res:red:cs:proj}.
The second equality follows from arguments analogous to those used in there.
The cycles partition $\bra{1, ..., \nn}$ into blocks and the integer partition records how many blocks of each size there are.
It does not, however, record \emph{where} the blocks are located.
Eg, the partitions of $\bra{1, 2, 3, 4, 5, 6}$ given by $\bra{ \bra{1,2}, \: \bra{3,4,5,6} }$ and $\bra{ \bra{1,2,3,4}, \: \bra{4,5} }$ are different partitions, yet have the same block sizes.
%The same symmetry as in the argument for \cref{res:red:cs:proj} shows that
By symmetry,
the law of the \CS given its integer partition is uniform over all \CSs with this given integer partition.
\end{Proof}

\begin{rmkt}[Limiting Distribution]
The limiting invariant distribution \asinf \nn is known to be the so-called \textit{Poisson--Dirichlet} distribution with parameter $\theta = \tfrac12$, denoted $\PD(\theta)$. More precisely, the joint law of the rescaled cycle sizes converge in distribution to $\PD(\tfrac12)$; see \bklmc[Theorem~1.1].
This is revisited in more detail in \S\ref{sec:cf:tiling}; see, in particular, \cref{res:cf:tiling:pd} and the surrounding discussion.
A related Poisson--Dirichlet limit is investigated by \textcite{P:partitions-pd}.
%A slightly more general situation is established by \textcite{P:partitions-pd}.
%As we see below, \emph{splits} are accepted with probability $\tfrac12$; \textcite{P:partitions-pd} deals with a general probability $\theta \in [0,1]$.
	%
\end{rmkt}

\section{Decomposing a Perfect Matching into a Sequence of Swaps}
\label{sec:decomp}

\subsection{Generating a $\kk$-\PM via $\kk-1$ Swaps}
\label{sec:decomp:gen}

Let $\eta \in \mfm$ be an $\nn$-\PM.
Consider a single step of the $\kk$-\PM \RW, starting from $\eta$:
	$\kk$ pairs are chosen uniformly at random (\textit{\uar});
	the \PM restricted to these $2\kk$ objects is resampled and the remaining $2(\nn-\kk)$ objects are left alone.
The resampled object is, up to a permutation of the labels, a $\kk$-\PM.
Being able to sample a $\kk$-\PM uniformly is then sufficient in order to run the dynamics.
%Recall \cref{rmk:intro:set-up:re-pairing} for more details.
We now describe a way to sample a $\kk$-\PM \uar choosing only $2$ matches at a time.

Choose an arbitrary cycle structure $\cc \in \mfc'$; this has at most $\kk$ non-fixed points and corresponds to a $\kk$-rematching inside an $\nn$-\PM.
We show how to draw $\eta$ \uar conditional on having \CS $\cc$, ie on $\mcc(\eta) = \cc$.
Recall that the \CS $\cc$ implicitly partitions $[\kk]$ and, by symmetry, the partition is uniform amongst all partitions with appropriately sized parts.
The relative matching inside different blocks of the partition is independent.
Thus, it suffices to be able to sample a single cycle of arbitrary length, ie an $\ell$-\PM with one $\ell$-cycle for any $\ell \in \mbn$.
This is analogous to sampling uniformly a permutation given its \CS.
\PMs with a single $\ell$-cycle are elements of $\widebar \mfm_\ell \cq \mfc_\infty(\delta_\ell)$ where $\delta_\ell(m) \cq \ONE\bra{m = \ell}$ for $\ell, m \in \mbn$.
How to sample such a \PM is described~in~\cref{alg:decomp:gen:gen-cycle}~below.

It is well-known and easy to generate an $\ell$-cycle permutation via $\ell-1$ transpositions, ie $2$-cycles. This does not generalise to \PMs, however.
This was a highly non-trivial obstacle for us.

\begin{figure}
	\centering
%\noindent%
\begin{minipage}{0.9\textwidth}
	\includegraphics[width = 1\textwidth]{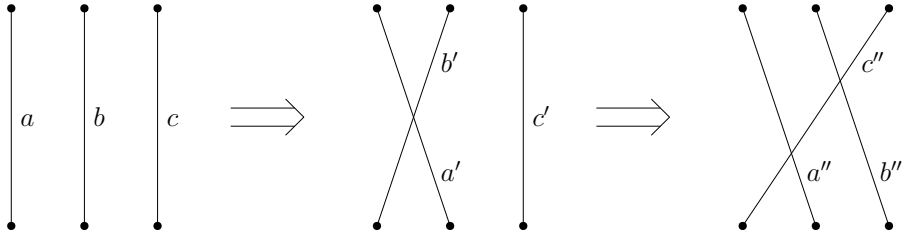}
	
	\caption{%
		The first step interacts with only $\bra{a, b}$; we can thus associate $c' = c$.
		The second step interacts with $\bra{b', c' = c}$; there is no natural way to say whether $b'$ corresponds to $a$ or to $b$.
		Our algorithm
%				does not need any such association:
			only needs the equality $a \cup b = a' \cup b'$ (as sets):
			it chooses $U \sim \Unif(\bra{a',b'})$ does a swap with $\bra{U, c}$.
	}
	\label{fig:destroy_label}
\end{minipage}
\vspace*{-1.5mm}
\end{figure}

%\begin{rmkt}[Generating Cyclic Permutations and \PMs]
%\label{rmk:decomp:gen:difficulty}
	%
%We now present a way of sampling an $\ell$-cycle permutation via $\ell-1$ transpositions, ie $2$-cycles.
%\cref{alg:decomp:gen:gen-cycle} translates this method from permutations to \PMs.
We first describe the usual way to sample an $\ell$-cycle permutation \uar.
%, which is not amenable to \PM sampling.
Generation of an $\ell$-cycle permutation is trivial for $\ell = 1$ as there is only one $1$-cycle.
Assume now that $\ell > 1$.
We use an inductive construction.
Choose independently
	$a_1 \sim \Unif([\nn])$,
then
	$a_2 \sim \Unif([\nn] \setminus \bra{a_1})$,
then
	$a_3 \sim \Unif([\nn] \setminus \bra{a_1, a_2})$
and
	so on.
Then a uniformly random $\ell$-cycle can be generated via
\[
	\sigma
\cq
	(a_1, a_2) (a_2, a_3) \cdots (a_{\ell-1}, a_\ell)
=
	(a_1, ..., a_\ell).
\]
This does not generalise well to \PMs because there is no concept of ``using the second element of the previous transposition as the first element of the current transposition'': which of the new matches is the `second' one?
See \cref{fig:destroy_label} for a visualisation of this difficulty.
This difficulty and highly related ones will rear its ugly head repeatedly throughout our arguments.

We can adjust this method for sampling an $\ell$-cycle permutation in a subtle way, which is then amenable to \PMs.
Let $b_\ii \sim \Unif(\bra{a_1, ..., a_\ii})$ independently for each $\ii \in [\ell]$ and set
\[
	\sigma
\cq
	(b_1, a_2) (b_2, a_3) \cdots (b_{\ell-1}, a_\ell).
\]
It is easy to check that $\sigma$ is still a uniform $\ell$-cycle.
We are ``using a uniformly random previously-used element as the first element of the current transposition''.
We can translate this into the realm of \PMs
	since we do know the set of previously interacted with pairs,
	so can choose one \uar.
%	we are not able to label the previously interacted with pairs individually,
%	but we do know the set of such pairs;
%	we can thus choose one uniformly at random.
	%
%\end{rmkt}

\begin{algt}[Generating a Uniform Cycle via Swaps]
\label{alg:decomp:gen:gen-cycle}
	Initialise $\eta_0 = \bra{\eta_{0, 1}, ..., \eta_{0, \ell}} \cq \id_\ell$.
	Choose $\ii \sim \Unif([\ell])$ and set $\SS_0 \cq \bra{\ii}$.
	Perform the following steps sequentially for $\ss = 1, ..., \ell-1$.
	
	\begin{itemize}
		\item 
		Choose $\ii \sim \Unif(\SS_{\ss-1})$ and $\jj \sim \Unif([\ell] \setminus \SS_{\ss-1})$ independently.
		
		\item 
		Choose uniformly a new matching on $\eta_{\ss-1, \ii} \cup \eta_{\ss-1, \jj}$ conditional on not being equal to $\bra{ \eta_{\ss-1, \ii}, \: \eta_{\ss-1, \jj} }$.
		Denote this new matching $\bra{ x, \: y }$.
		Note that $x \cup y = \eta_{\ss-1, \ii} \cup \eta_{\ss-1, \jj}$.
		
		\item 
		Set $\eta_{\ss, m} \cq \eta_{\ss-1, m}$ for $m \notin \bra{\ii,\jj}$,
		set $\eta_{\ss, \ii} \cq x$
		and
		set $\eta_{\ss, \jj} \cq y$.
%		(The roles of $x$ and $y$ can be interchanged.)
%		
%		\item 
		Set $\SS_\ss \cq \SS_{\ss-1} \cup \bra{\ii}$.
	\end{itemize}
	Output $\Eta \cq \eta_{\ell-1} = \bra{ \eta_{\ell-1, 1}, ..., \eta_{\ell-1, \ell} }$.
\end{algt}

\begin{lem}[Generating a Uniform Cycle via Swaps]
\label{res:decomp:gen:gen-cycle}
	Let $\ell \in \mbn$.
	Let $\Eta$ denote the (random) output of \cref{alg:decomp:gen:gen-cycle}.
	Recall that
	\(
		\widebar \mfm_\ell
	=
		\mfc_\infty(\delta_\ell)
	\)
	is the set of single-cycle $\ell$-\PMs.
	Then,
	\[
		\Eta \sim \Unif\rbb{ \widebar \mfm_\ell }.
	\]
\end{lem}

\begin{Proof}
We use induction on $\ell$.
The bases cases $\ell \in \bra{1,2}$ are trivial.
The case $\ell = 3$ is easy to check by hand.
Now assume that the claim holds for $\ell$. We establish it for $\ell+1$.

The algorithm consists of $\ell$ steps.
We break it into two parts:
	the first $\ell-1$ steps and the final step.
We use the notation from \cref{alg:decomp:gen:gen-cycle}.
Note that $\abs{\SS_\ss} = \ss+1$.
Let $\ii$ be the (random) unique element of $[\ell+1] \setminus \SS_{\ell-1}$.
By symmetry, $\ii \sim \Unif([\ell+1])$.
The inductive hypothesis implies~that
\[
	\fnrestrict{\eta_{\ell-1}}{[\ell+1] \setminus \bra{i}} \sim \Unif([\ell+1] \setminus \bra{i})
\Quad{given}
	\ii,
\Quad{or equivalently}
	\SS_{\ell-1}.
\]
That is, the restriction of $\eta_{\ell-1}$ to $[\ell+1] \setminus \bra{\ii} = \SS_{\ell-1}$ is a uniformly random $\ell$-\PM on its support, given $\ii$.
The nature of \cref{alg:decomp:gen:gen-cycle} means that the remaining matched pair is untouched:
\(
	\eta_{0, \ii}
%=
%	\eta_{1, \ii}
=
	\cdots
=
	\eta_{\ell-1, \ii}.
\)
The final step comprises a swap with the $\ii$-th pair and the $\jj$-th, where $\jj \sim \Unif([\ell+1] \setminus \bra{\ii})$ and is independent of $\ii$.
%It is straightforward to check that
This leads to a uniform, single-cycle $(\ell+1)$-\PM, as desired.
%\footnote{%
%	In fact, this does not even use the uniformity of $\ii$.
%	The proof is analogous to that used to show that the uniform distribution is invariant for the dynamics of the top-to-random card shuffle%
%}
	%
\end{Proof}

We now use this to generate a uniform $\kk$-rematching in the space of $\nn$-\PMs, ie an $\nn$-\PM with at most $\kk$ non-fixed points. 
We break down a \PM into its individual cycles according to its \CS.

\begin{algt}[Generating a Uniform $k$-Matching via Its Cycle Decomposition]
\label{alg:decomp:gen:gen-match}
%	Let $\kk, \nn \in \mbn$ with $\kk \in [2, \nn]$.
	Let $\kk \in [2, \nn] \cap \mbn$.
%	Define $\Eta$ according to the following steps and output it.
	\begin{itemize}
		\item 
		Draw $C \sim \mcc\rbr{ \Unif(\mfm') }$.
		
		\item 
		Draw partition $\PP = (\PP_1, ..., \PP_\SS)$ uniformly, conditional on having block lengths given by $C$.
%		In particular, $\SS = \sumt[\infty]{\ell=1} C_\ell$.
		
		\item 
%		Let $\SS \cq \sumt[\infty]{\ell=1} C_\ell$.
%		; this is the number of blocks in the partition.
		Draw $\Eta_\ss \sim \Unif(\widebar \mfm_{\abs{\PP_\ss}})$,
		which is a single $\abs{\PP_\ss}$-cycle,
		independently for each $\ss \in [\SS]$.
		
		\item 
		Combine to create $\Eta$:
			place the $\ss$-th cycle $\Eta_\ss$
			in the $\ss$-th block $\PP_\ss$
			for each $\ss \in [\SS]$.
%	\qedhere
	\end{itemize}
	Output $\Eta$.
\end{algt}

It is immediate from \cref{res:decomp:gen:gen-cycle} and the cycle decomposition that \cref{alg:decomp:gen:gen-match} gives rise to a uniform element of $\mfm'$, ie $\nn$-\PM with at most $\kk$ non-fixed points.

\begin{cor}[Generating a Uniform $k$-Matching via Its Cycle Decomposition]
\label{res:decomp:gen:gen-match}
	Let $\kk, \nn \in \mbn$ with $2 \le k \le n$.
%	Let $\kk \in [2, \nn] \cap \mbn$.
	Let $\Eta$ denote the (random) output of \cref{alg:decomp:gen:gen-match}.
	Then,
	\[
		\Eta
	\sim
		\Unif\rbb{ \mfm' }.
	\]
\end{cor}

%\begin{Proof}
%	%
%This follows immediately from the results we have already proved.
%	%
%\end{Proof}

\subsection{Support Size and Distance from Identity for a Uniform $\kk$-Rematch}
\label{sec:decomp:dist-supp}

We now know how to sample a $\kk$-\PM uniformly at random given its \CS.
A priori, one may assume that we must now calculate the law \CS of a uniform $\kk$-\PM.
It turns out that our proof does not require this, however.
The only information we need is the \emph{support} of the \CS---namely, the number of pairs interacted with. Eg, the support of the \CS $(0, 2)$ is $4$ and of $(0, 0, 1)$ is $3$; see \cref{def:decomp:dist-supp:supp}.

We take inspiration from the work of \bst on conjugacy-invariant \RWs.
They show that the mixing time is inversely proportional to the \textit{support} of the \CS used.

\begin{defn}[Support]
\label{def:decomp:dist-supp:supp}
The (\textit{size of the}) \textit{support} of
	a \CS $\cc \in \mfc_\infty$
and
	a \PM $\eta \in \mfm_\infty$
is
\[
	\# \cc
\cq
	\sumt[\infty]{\ell=2}
	\ell \cc_\ell
\Qand
	\# \eta
\cq
	\# \mcc(\eta),
\]
respectively.
This is the number of non-fixed points:
\[
	\text{if}
\quad
	\mcc(\eta) = \cc \in \mfc^\kk,
\Qthen
	\# \eta
=
	\# \cc
=
	\kk - \cc_1.
\]

We can view a $\kk$-\PM as an $\nn$-\PM by padding the end with $\nn-\kk$ fixed points.
More formally, view an element $\eta \in \mfm_\kk$ as an element $\eta' \in \mfm' \subseteq \mfm_\nn$ by setting $\eta'_\ii \cq \ii$ for $\ii \in \bra{2\kk+1, ..., \nn}$.
Let $\cc \cq \mcc(\eta)$ and $\cc' \cq \mcc(\eta')$ denote the \CSs.
Then $\cc'_1 = \cc_1 + (\nn - \kk)$ and $\cc'_\ell = \cc_\ell$ for $\ell \ge 2$.
Thus,
\[
	\# \cc
=
	\sumt[\infty]{\ell=2}
	\ell \cc_\ell
=
	\sumt[\infty]{\ell=2}
	\ell \cc'_\ell
=
	\# \cc'
\Qand
	\# \cc
=
	\kk - \cc_1
=
	\nn - \cc'_1
=
	\# \cc'.
\]
\end{defn}

\bst use $\abs \cdot$ to denote the support.
This already has an established meaning of ``size'' or ``cardinality'' for sets, which \PMs are.
We use $\#$ to avoid this notational~clash.

We also define the \textit{swap distance}.
This is just the number of swaps required to reach the identity.

\begin{defn}[Swap Distance]
\label{def:decomp:dist-supp:dist}
For
	a \CS $\cc \in \mfc_\infty$
and
	a \PM $\eta \in \mfm_\infty$,
define
\[
	d(\cc)
\cq
	\sumt[\infty]{\ell=2}
	(\ell-1) \cc_\ell
\Qand
	d(\eta)
\cq
	d\rbb{ \mcc(\eta) }.
\]
We refer to $d$ as the \textit{swap distance from the identity}, often referred to as just \textit{distance} for brevity.
Equivalently, $d(\eta)$ is the minimal number of swaps required to reach $\id_{\abs \eta}$. Indeed, each $\ell$-cycle needs precisely $\ell-1$ swaps to resolve its disparities compared with the identity.

Analogously, define $d(\eta, \eta')$ to be the minimal number of swaps required to move from $\eta$ to $\eta'$ for $\eta, \eta' \in \mfm_\infty$ with $\abs{\eta} = \abs{\eta'}$.
This is the usual distance in the graph which has \PMs as vertices and  edges between \PMs which differ by a single swap.
This graph is transitive.
Given $(\eta, \eta') \in \mfm \times \mfm$, relabel the objects so that these become $(\tilde \eta, \tilde \eta') \in \mfm \times \mfm$ with $\tilde \eta' = \id$. Then,~$d(\eta, \eta') = d(\tilde \eta)$.
\end{defn}

\begin{figure}
\centering
\begin{subfigure}[b]{0.45\textwidth}
	\includegraphics[width=\textwidth]{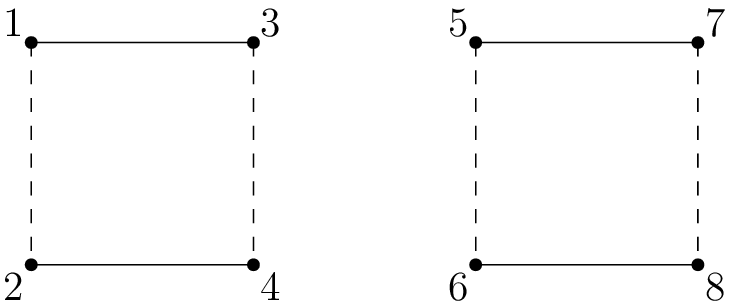}
%	\smallskip
	\caption{Graph corresponding to \cref{exm:decomp:dist-supp}	\ref{exm:a}}
	\label{fig:exm:a}
\end{subfigure}
\hfill
\begin{subfigure}[b]{0.45\textwidth}
	\includegraphics[width=\textwidth]{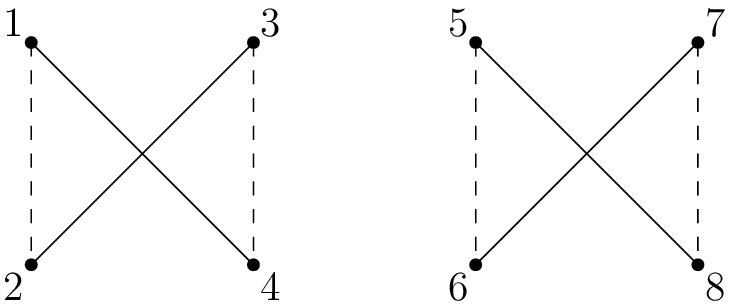}
%	\smallskip
	\caption{Graph corresponding to \cref{exm:decomp:dist-supp}	\ref{exm:b}}
	\label{fig:exm:b}
\end{subfigure}

\bigskip

\begin{subfigure}[b]{0.45\textwidth}
	\includegraphics[width=\textwidth]{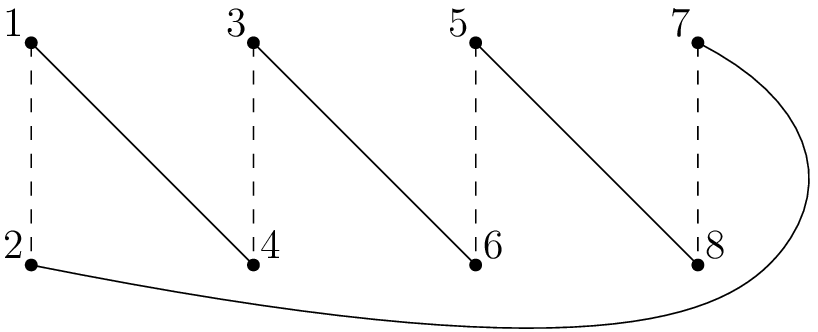}
%	\smallskip
	\caption{Graph corresponding to \cref{exm:decomp:dist-supp}	\ref{exm:c}}
	\label{fig:exm:c}
\end{subfigure}
\hfill
\begin{subfigure}[b]{0.45\textwidth}
	\includegraphics[width=\textwidth]{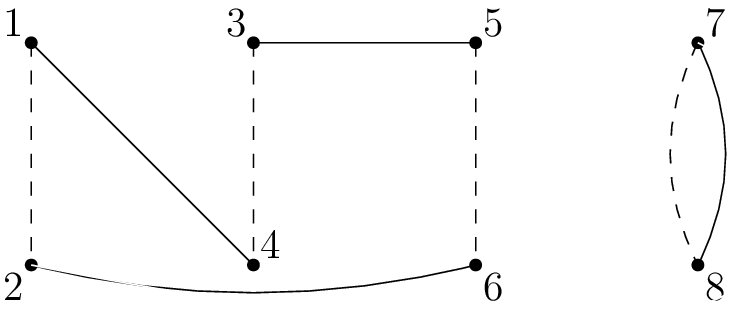}
%	\smallskip
	\caption{Graph corresponding to \cref{exm:decomp:dist-supp}	\ref{exm:d}}
	\label{fig:exm:d}
\end{subfigure}

\medskip

\caption{%
	Four examples from \cref{exm:decomp:dist-supp} \ref{exm:a}--\ref{exm:d}.
	The dashed lines correspond to the identity matching $\id_8$.
	The solid lines correspond to the matching $\eta$ from the respective examples.
}
\label{fig:exm}
\end{figure}

\begin{exmt}[Some Examples]
\label{exm:decomp:dist-supp}
Consider \PMs on $8$ elements, shown in
\cref{fig:exm}.
%\cref{fig:exm:a,fig:exm:b,fig:exm:c,fig:exm:d}.

\begin{enumerate}[label = \ensuremath{(\alph*)}]
\item \label{exm:a}
Take
\(
	\eta
\cq
	\bra{ \bra{1,3}, \: \bra{2,4}, \: \bra{5,7}, \: \bra{6,8} }.
\)
The graph is shown in \cref{fig:exm} $\subref{fig:exm:a}$.

Then $c \cq \mcc(\eta) = (0,2,0,0)$.
Hence $\# \eta = \# c = 4$ and $d(\eta) = d(c) = 2$.

\item \label{exm:b}
Take
\(
	\eta
\cq
	\bra{ \bra{1,4}, \: \bra{2,3}, \: \bra{5,8}, \: \bra{6,7} }.
\)
The graph is shown in \cref{fig:exm} $\subref{fig:exm:b}$.

Then $c \cq \mcc(\eta) = (0,2,0,0)$.
Hence $\# \eta = \# c = 4$ and $d(\eta) = d(c) = 2$.

\item \label{exm:c}
Take
\(
	\eta
\cq
	\bra{ \bra{1,4}, \: \bra{3,6}, \: \bra{5,8}, \: \bra{7,2} }.
\)
The graph is shown in \cref{fig:exm} $\subref{fig:exm:c}$.

Then $c \cq \mcc(\eta) = (0,0,0,1)$.
Hence $\# \eta = \# c = 4$ and $d(\eta) = d(c) = 3$.
%	Any `$4$-cycle', ie element of $\widebar \mfm_4$, satisfies this.

\item \label{exm:d}
Take
\(
	\eta
\cq
	\bra{ \bra{1,4}, \: \bra{2,6}, \: \bra{3,5}, \: \bra{7,8} }.
\)
The graph is shown in \cref{fig:exm} $\subref{fig:exm:d}$.

Then $c \cq \mcc(\eta) = (1,0,1,0)$.
Hence $\# \eta = \# c = 3$ and $d(\eta) = d(c) = 2$.
\qedhere
\end{enumerate}
\end{exmt}

The key information required from the law of the \CS, which is that of a uniform $\kk$-\PM, is its expected support and distance.
The mixing time actually only depends on the expected support; the expected distance is a tool used in the proof which need not be calculated explicitly.

\begin{defn}[Expected Support and Distance of Uniform $k$-Rematching]
\label{def:decomp:dist-supp:param-exp}
Let $\kk \in \mbn \setminus \bra{1}$.
Let $\Eta \sim \Unif( \mfm_\kk )$.
Define the \textit{expected support} $\kappa_\kk$ and \textit{expected distance} $\rho_\kk$ as follows:
\[
	\kappa_\kk
&
\cq
	\ex{ \# \Eta }
=
	\ex{
		\sumt[\infty]{\ell=2}
		\ell \mcc_\ell(\Eta)
	}
=
	\abs{\mfm_\kk}^{-1}
	\sumt{\eta \in \mfm_\kk}
	\sumt[\kk]{\ell=2}
	\ell \mcc_\ell(\eta);
\\
	\rho_\kk
&\cq
	\ex{ d(\Eta) }
=
	\ex{
		\sumt[\infty]{\ell=2}
		(\ell-1) \mcc_\ell(\Eta)
	}
=
	\abs{\mfm_\kk}^{-1}
	\sumt{\eta \in \mfm_\kk}
	\sumt[\kk]{\ell=2}
	(\ell-1) \mcc_\ell(\eta).
\]
\end{defn}

We are always interested in the expected support of a $\kk$-\PM.
Somewhat unusually, we abbreviate $\kappa \cq \kappa_\kk$.
Officially, $\kk$ is a function of $\nn$, so this is suppressing the $\nn$-dependence via abbreviating~$\kappa_{\kk_\nn}$.

\begin{lem}[Expected Support of Uniform $k$-Rematching]
\label{res:decomp:dist-supp:exp-supp}
	We have
	\[
		\kappa_\kk
	=
		\ex{ \# \Eta }
	=
		\kk - \tfrac{\kk}{2\kk-1}
	=
%		\kk - \tfrac12 \tfrac1{1 - 1/(2\kk)}
%	=
%		\kk - \tfrac12 \rbb{ 1 + \tfrac1{2\kk} + \Ohb{ \kk^{-2} } }
%	=
		\kk - \tfrac12 - \tfrac1{4\kk} + \Ohb{ \kk^{-2} }.
	\]
\end{lem}

%\begin{lem}[Expected Support and Distance of Uniform $k$-Rematching]
%\label{res:decomp:dist-supp:exp-supp}
%	Draw $\Eta \sim \Unif(\mfm_\kk)$, ie a uniform $\kk$-rematching.
%	Then
%	\[
%		\kappa_\kk
%	\cq
%		\ex{ \# \Eta }
%	=
%		\kk - \tfrac{\kk}{2\kk-1}
%	=
%		\kk - \tfrac12 + \Oh{1/\kk}.
%	\]
%\end{lem}

\begin{Proof}
This follows by a simple counting argument.
Indeed, $\# \Eta$ is simply $\kk$ minus the number of fixed points.
Thus, we just need to calculate the number of fixed points in expectation.
Start with the pair $\bra{1,2}$ matched. There are $2k-1$ other vertices to which $1$ can be matched. So the probability that it remains matched to $2$ is $1/(2\kk-1)$.
The expected number of fixed points is then $\kk/(2\kk-1)$, by linearity of expectation, as $\kk$ matches are made.
The lemma follows.
\end{Proof}

\begin{rmkt}[Support and Its Relation to Mixing]
\label{rmk:decomp:dist-supp:param-motivation}
One can think of the support as ``the number of random choices''.
We analogise with permutations:
	a $3$-cycle $(a, b, c)$ can be written as $(a, b) (b, c)$ and there are three choices, namely $a$, $b$ and $c$;
	a double-transposition $(a, b) (c, d)$ has four choices.

\bst analyse the mixing time of the \RW on the Cayley graph of the permutation group generated by a preset \CS $\Gamma$: a step comprises applying a uniform permutation given the \CS.
They show that the mixing time is inversely proportional to the support $\# \Gamma$.

The lower bound given by \bsc finds the time it takes for all cards to be touched.
Decomposing a permutation into its \CS, the number of cards touched equals the support $\# \Gamma$.
A coupon-collector argument can be applied when this support has size $\oh n$ to deduce a lower bound of $(\# \Gamma)^{-1} \nn \log \nn$ when there are $\nn$ cards.
It had long been conjectured that this, ie the time at which all cards have been touched, is indeed the correct mixing time. \bst establish this.

We adjust this heuristic to \PMs.
The number of pairs interacted with in a given rematching is the support of that rematching. Suppose $\kk_\tt$ are interacted with on the $\tt$-th step.
We wait until all the original pairs have been interacted with.
We want to apply a coupon-collector counting argument to estimate this time.
We have to be careful, though.
	Suppose that
		$\bra{ \bra{1, 2}, \bra{3, 4}, \bra{5, 6} }$
	becomes
		$\bra{ \bra{1, 3}, \bra{2, 4}, \bra{5, 6} }$
	and then
		$\bra{ \bra{1, 3}, \bra{2, 5}, \bra{4, 6} }$;
	the first two pairs are interacted with in the first step,
	but which are in the second step?
	Certainly the third pair is, but is the first or the second?
	It does not matter: the first and second pairs have already been `collected' in the first step. All that matters is that the second step included the third pair.
	Any as-yet `uncollected' pair is in its original position, by definition.
%This allows us to apply the coupon-collector argument in the usual way, except with a variable number of picks in each step.
This allows us to apply the coupon-collector argument it the usual way:
	the collection takes time approximately
	\(
		\inf\bra{ \tt \ge 0 \mid \kk_1 + \cdots + \kk_\tt \ge \nn \log \nn }.
	\)

The steps are independent, so the law of large numbers says that $\kk_1 + \cdots + \kk_\tt \approx \tt \kappa$.
It is thus natural to conjecture a mixing time of
\(
	\kappa^{-1} \nn \log \nn.
\)
This is what our main theorem verifies.
\end{rmkt}

\section{Analysis of Coalescence--Fragmentation Chain}
\label{sec:cf}

\subsection{Conditional Uniformity}
\label{sec:cf:cu}

Recall that a general $\kk$-\PM can be written as a product of disjoint, single-cycle \PM.
The different single-cycle \PMs correspond to different blocks in the partition.
The order of their application is thus inconsequential.
The single-cycle \PMs are broken down into swaps.
See \cref{alg:decomp:gen:gen-cycle} for more details.
General permutations of a given \CS have a similar independence property.

These properties lead us to the notion of a \textit{refresh time}, corresponds to the start of a new block.

%Following \bsc[\S4.2], we define \textit{refresh times}, \textit{conditional uniformity} and \textit{relaxed conditional uniformity}.

\begin{defn}[Refresh Time]
\label{def:cf:cu:refresh}
	Let $\cc \in \mfc_\infty$.
	Recall that $d(\cc) = \sumt[\infty]{\ell=2} (\ell-1) \cc_\ell$ is the swap distance from the identity.
	Sampling uniformly from $\mfm_\infty(\cc)$ involves $d(c)$ swaps, grouped together in different batches:
		the application of an $\ell$-cycle requires $\ell-1$ swaps;
		each batch corresponds to a block in the associated partition.
	Call $\ss \in \bra{1, ..., d(\cc)}$ a \textit{refresh time for $\cc$} if it is of the form $\ss = \sumt[m]{\ii=2} (\ii-1) \cc_\ii + 1$ for some $m \in \mbn$.
	In particular, $1$ is a refresh time, since the empty sum is $0$.
\end{defn}

\bst[Definition 4.1] define refresh times similarly.
They always apply a permutation with the same \CS. This means that all the refresh times can be defined in advance.
Our \CS is not preset, but varies from step to step. Thus, the refresh times vary according to the \CS of the \PM chosen in a given step. It will be enough, however, to condition in advance on the sequence of \CSs, from which we can define the refresh times.

We now describe \cref{alg:decomp:gen:gen-cycle,alg:decomp:gen:gen-match} in terms of these refresh times.
The reformulation is given as \cref{alg:cf:cu:original}.
Importantly, we distinguish the first and second markers in the choice of two matches for a swap:
	the first marker is $\ii$ and the second $\jj$ in \cref{alg:decomp:gen:gen-cycle}.

\begin{algt}[Conditional Uniformity]
\label{alg:cf:cu:original}
Let $c \in \mfc$ be a \CS.
Let $\ss \in \mbn$.
%If $\ss = 1$, then set
Set
\(
	\SS'_0 \cq \emptyset
\Qand
	\SS_0 \cq \emptyset.
\)
If $\ss > 1$, then 
let $\rbr{\ii_1, \jj_1}, ..., \rbr{\ii_{\ss-1}, \jj_{\ss-1}}$ denote the pairs chosen in the first $\ss-1$ steps
and set
\[
	\SS'_{\ss-1}
\cq
	\cup_{\rr=1}^{\ss_--1}
	\bra{\ii_\rr, \jj_\rr}
\Qand
	\SS_{\ss-1}
\cq
	\cup_{\rr=\ss_-}^{\ss-1}
	\bra{\ii_\rr, \jj_\rr},
\]
where
\(
	\ss_-
\cq
	\sup\bra{ \ss' \le \ss \mid \text{$\ss'$ is a refresh time} }
\)
is the most recent refresh time before $\ss$.
This way,
$\SS'_{\ss-1}$, respectively $\SS_{\ss-1}$,
is the set of indices used in the
previous blocks, respectively current block.\footnote{%
	This definition of $(\SS_\ss)_{\ss\ge0}$ is a natural extension of \cref{alg:decomp:gen:gen-cycle} where only a single cycle is considered.
	Performing a swap destroys the two original matches and two new ones are created; there is no real way of associating the old matches with the new ones.
%	, since the labels are somewhat arbitrary.
	However, the set union of the objects interacted with is always well-defined}
%If $\ss = 1$, then set
%\(
%	\SS'_0 \cq \emptyset
%\Qand
%	\SS_0 \cq \emptyset.
%\)

Perform the following steps sequentially for $\ss = 1, ..., d(\cc)$.
There are two cases according to whether $\ss$ is a refresh time for $\cc$ or not.
If $\ss$ is a refresh time, then $\ss_- = \ss$ and thus $\SS_{\ss-1} = \emptyset$.

\begin{itemize}
	\item 
	If $\ss$ is a refresh time for $\cc$,
		then it corresponds to the start of a new cycle.
	\begin{itemize}[noitemsep, topsep = 0pt]
		\item 
		The first marker $\ii_\ss$ is chosen uniformly on
			$[\nn] \setminus \SS'_{\ss-1}$.
		
		\item 
		The second marker $\jj_\ss$ is chosen uniformly on
			$[\nn] \setminus \rbr{ \SS_{\ss-1} \cup \bra{\ii_\ss} }$.
	\end{itemize}
	
	\item 
	If $\ss$ is not a refresh time for $\cc$,
		then it corresponds to the continuation of a cycle.
	\begin{itemize}[noitemsep, topsep = 0pt]
		\item 
		The first marker $\ii_\ss$ is chosen uniformly on
			$\SS_{\ss-1}$.
			
		\item 
		The second marker $\jj_\ss$ is chosen uniformly on
			$[\nn] \setminus \rbr{ \SS'_{\ss-1} \cup \SS_{\ss-1} \cup \bra{\ii_\ss} }$.
	\end{itemize}
	
	\item 
	Perform a uniform swap of the pairs $\ii_\ss$ and $\jj_\ss$.
\qedhere
\end{itemize}
\end{algt}

\begin{lem}[Conditional Uniformity Algorithm]
	Let $\cc \in \mfc$.
	The output of \cref{alg:cf:cu:original} is a uniform \PM with cycle structure $\cc$,
	ie is a uniform element of $\mfm(\cc)$.
\end{lem}

\begin{Proof}
This is an immediate consequence of
the algorithms and results of \S\ref{sec:decomp:gen}.
%\cref{alg:decomp:gen:gen-cycle,alg:decomp:gen:gen-match,res:decomp:gen:gen-cycle,res:decomp:gen:gen-match}.
	%
\end{Proof}

%	Perform a uniform swap with $\rbr{\ii_\ss, \jj_\ss}$ in each step $\ss$.
%	This generates a uniform \PM with \CS $\cc$.
%	%	, ie uniform element of $\mfm_\infty(\cc)$.

We next define a concept of \emph{relaxed} conditional uniformity:
	in essence, we relax the correlation between the different single-cycles in the decomposition of the \PM.
We want to pretend that we can sample the markers $\ii$ and $\jj$ as follows:
	$\ii$ completely uniform for refresh times and uniform on the indices used so far in the current cycle otherwise;
	$\jj$ uniform on everything except $\ii$.
Notationally,
	$\ii \sim \Unif([\nn])$ for refresh times
and
	$\ii \sim \Unif(\SS_{\ss-1})$ for non-refresh times;
	$\jj \sim \Unif([\nn] \setminus \bra{\ii})$ always.
This is, of course, not possible since the different single-cycles must correspond to different blocks of the corresponding partition. This will not always be the case in the relaxed version.
Were this possible, however, it would make analysing the chain considerably simpler.
We show that the relaxed version can be coupled with the original for a long enough period of time for us to couple.

\bst[Definition 4.2] define an analogous relaxation, although it is somewhat simpler in their case because their cycle structure is unchanging, unlike ours.
Further, they can always take $\ii_\ss \cq \jj_{\ss-1}$, ie the first marker for the current swap to be the second marker from the previous swap.
We cannot do this due to the previously-discussed lack of identifiability.

\begin{algt}[Relaxed Conditional Uniformity]
\label{alg:cf:cu:relaxed}
Let $\cc \in \mfc$ be a cycle structure.
Use the same notation for $( \SS_{\ss-1} )_{\ss \ge 0}$ as in \cref{alg:decomp:gen:gen-cycle,alg:cf:cu:original}.
Perform the following steps sequentially for $\ss = 1, ..., d(\cc)$.
There are two cases according to whether $\ss$ is a refresh time for $\cc$ or~not.

\begin{itemize}
	\item 
%	Suppose that $\ss$ is a refresh time.
	If $\ss$ is a refresh time for $\cc$, then
	sample
		$\ii_\ss \sim \Unif([\nn])$
	and
		$\jj_\ss \sim \Unif([\nn] \setminus \bra{\ii_\ss})$.
	
	\item 
%	Suppose that $\ss$ is not a refresh time.
	If $\ss$ is not a refresh time for $\cc$, then
	sample
		$\ii_\ss \sim \Unif(\SS_{\ss-1})$
	and
		$\jj_\ss \sim \Unif([\nn] \setminus \bra{\ii_\ss})$.
	
	\item 
	Perform a uniform swap with $\rbr{\ii_\ss, \jj_\ss}$ in each step $\ss$.
\qedhere
\end{itemize}
%This is using the same notation for $( \SS_{\ss-1} )_{\ss \ge 0}$ as in \cref{alg:decomp:gen:gen-cycle,alg:cf:cu:original}.
	%
\end{algt}

\begin{defn}[Relaxed Conditional Uniformity]
\label{def:cf:cu:relaxed}
	We call the evolution defined by \cref{alg:cf:cu:relaxed} the \textit{relaxed law} when \CS $\cc$ is chosen independently in each step and with distribution $\mcc(\Unif(\mfm))$.
\end{defn}

The choices of $(i,j)$ in the relaxed version (\cref{alg:cf:cu:relaxed}) clearly can violate the conditions in the original (\cref{alg:cf:cu:original}).
The next lemma shows that the relaxed version does not violate the conditions \whp when order $1$ steps are taken, since $\kk = \oh \nn$.

\begin{lem}[Relaxed Conditional Uniformity]
\label{res:cf:cu:relaxed}
	Suppose that $\Delta$ steps are taken under relaxed conditional uniformity (\cref{alg:cf:cu:relaxed}).
	Then the probability that any choice of $(i,j)$ violates the original conditional uniformity conditions (\cref{alg:cf:cu:original}) is at most $2 \Delta \kk / \nn$.
\end{lem}

\begin{Proof}
The condition is violated in a given step if the marker falls in the set of those already chosen.
This set has size at most $\kk$.
Two markers are chosen each time in a uniform manner.
%The lemma follows.
	%
\end{Proof}

We henceforth proceed using the relaxed version of conditional uniformity.
This is convenient for adapting a coupling which is based on an idea of \textcite[\S3]{S:random-trans}.
The original coupling of \cite[\S3]{S:random-trans} is designed for random transpositions, which naturally satisfies the relaxed version.

\subsection{Schramm's Coupling for the Coalescence--Fragmentation Chain}
\label{sec:cf:tiling}

Key to the analysis of \bsn is their use of a coupling between two realisations of the \cf chain; see \cite[\S 4.2.2]{BS:cutoff-conj-inv}.
The same coupling had already been used by \bszc[\S 3].
The original idea is due to \textcite[\S 3]{S:random-trans}.
%\footnote{%
%	\textcite{S:random-trans} did not control coupling time, but rather at weak convergence of the appropriately rescaled partitions to a Poisson--Dirichlet distribution.
%	The same goes for \bklmt}
There is a crucial difference in \cite{BSZ:k-cycle,BS:cutoff-conj-inv} compared with \cite{S:random-trans}:
	the introduction of the measure-preserving map $\Phi$, given in \cref{alg:cf:tiling:coupling} below.
%	\footnote{%
%		The map $\Phi$ was introduced in \bszc in order to get quantitative control on the coupling time, which is required to establish cutoff. \textcite{S:random-trans}, and \bklmt later, did not need this quantification; taking $\Phi$ to be the identity was sufficient for them in their establishment of weak convergence}
The introduction of this map was one of the main innovations of \bszc.
A version of \citeauthor{S:random-trans}'s coupling has been used recently by \bklmt[\S 5.2] in a set-up similar to ours, but without the adaptation of \cite{BSZ:k-cycle,BS:cutoff-conj-inv}.

Our description follows closely that of \cite{BSZ:k-cycle,BS:cutoff-conj-inv}.
Some changes are required to take into account the fact that, for us, a single block does not always split when both markers fall in it.
%; recall \cref{rmk:red:part:evolution}.

First we describe the marginal evolution of the partition.
We view this as a tiling of $(0,1]$.
%Further, there will be a \textit{distinguished} tile.

\begin{Proof}[Set-Up for Tiling]
\qedtriangle
We describe how to simulate a single \emph{round}, ie $\kk$-rematch, via individual \emph{steps}, ie swaps ($2$-rematches).
To extend to multiple rounds, the procedure is repeated independently.

Let $c \in \mfc'$.
%, ie a \CS for a \PM on $\nn$ pairs subject to having at most $\kk$ non-fixed points, ie $\sumt[\infty]{\ii=1} \ii \cc_\ii = \nn$ and $\cc_1 \ge \nn - \kk$.
This corresponds to choosing a $\kk$-rematching amongst $\nn$ objects.
%Set $d \cq d(\cc) = \sumt[\kk]{\ii=2} (\ii-1) \cc_\ii$.
We apply $d(\cc)$ swaps. This involves choosing markers for each $\ss \in \bra{1, ..., d(\cc)}$. This choice is performed differently according to whether or not $\ss$ is a refresh time; recall \cref{def:cf:cu:refresh,alg:cf:cu:relaxed}.

We use $\PP = (\PP_\tt)_{\tt\ge0}$ to denote the \cf process.
We denote it $\widebar \PP = (\widebar \PP_\ss)_{\ss\ge0}$ when looking at a single step, broken down by swaps indexed by $\ss$. This implicitly assumes that the \CSs have been conditioned on. The process on the swap-timescale is then well-defined.
Both $\ss$ and $\tt$ here indicate \emph{time}:
	$\tt$ in the sense of the number of \emph{rounds},
	whilst $\ss$ in the sense of \emph{swaps}.
%We emphasise that the subscript $\tt$ in $\PP_\tt$ indicates \emph{time} in the sense of number of \emph{rounds}, whilst the subscript $\ss$ in $\widebar \PP_\ss$ denotes \emph{time} in 

The set $\mfp$ comprises all $\tfrac1\nn$-integer partitions.
Given $\pp = (\pp_1, ..., \pp_\nn) \in \mfp$, we tile the interval $(0, 1]$ using the intervals $\bra{ (0, \pp_1], ..., (0, \pp_\nn] }$---the specific rule does not matter.
We choose markers $\uu$ and $\vv$ in $\bra{1/\nn, ..., \nn/\nn} \subseteq (0,1]$ below and use them, scaled by $\nn$, as markers in \cref{alg:cf:cu:relaxed}.
%This implicitly uses the natural bijection between $\bra{1/\nn, ..., \nn/\nn}$ and $\bra{1, ..., \nn}$.
	%
\end{Proof}

%\begin{subtheorem}{thm}
%	\label{alg:cf:tiling:marg}
	%
\begin{defn}[Distinguished Tile]
\label{def:cf:tiling:distinguished}
Suppose that $\ss \ge 1$ and that $\widebar P_0, ..., \widebar P_{\ss-1} \in \mfp$ have been defined.
%Recall the definition of a refresh time from \cref{def:cf:cu:refresh}.

\begin{itemize}
	\item 
	If $\ss$ is a refresh time, which includes $\ss = 1$,
		then select $u \sim \Unif\rbr{\bra{1/\nn, ..., \nn/\nn}}$ and \textit{distinguish} the tile containing $u$; use it as the first marker in \cref{alg:cf:cu:relaxed}.
	
	\item 
	If $\ss$ is not a refresh time, and hence $\ss \ge 2$,
%		then \textit{distinguish} the tile containing the second marker of the previous swap, ie $\jj_{\ss-1}$ in \cref{alg:cf:cu:relaxed}.
		then \textit{distinguish} the tile containing the first marker of the current swap, ie $\ii_\ss$ in \cref{alg:cf:cu:relaxed}.
\end{itemize}
There is a distinguished tile containing the first marker for the step $\ss$ in either case.

Define $\widebar P_\ss$ to be the new partition, written in non-increasing block size order.
\end{defn}

\begin{rmkt}[Permutations vs \PMs]
\bst[\S 4.2.2] use an analogous \emph{distinguished tile} construction.
It is variant on fundamental ideas introduced by \textcite{S:random-trans}.
There is a key difference, however, in theirs compared with ours:
\begin{itemize}[noitemsep]
	\item 
	they take the distinguished tile to be the second marker from the previous swap, ie $\ii_\ss = \jj_{\ss-1}$;
	
	\item 
	we must choose $\ii_\ss$ uniformly from the already-used indices $\SS_{\ss-1} \ni \jj_{\ss-1}$.
\end{itemize}
%	recall \cref{rmk:decomp:dist-supp:param-motivation} for a discussion of this difficulty.
The marginal evolution of the distinguished tile is thus simpler in their set-up.

%This fundamental difference only arises when $\kk > 2$:
%	if $\kk = 2$, then every time is a refresh time and so both markers are chosen independently of previous steps.
\bklmt[\S 5.2] use a variant on \citeauthor{S:random-trans}'s coupling for the \PM \RW. They do not face similar issues, however, because they only study the $2$-\PM \RW. This means that each round involves choosing only a single pair to swap and rounds are independent.

We thus need to extend the content of \bsc[\S 4.2.2] and \bklmc[\S 5.2] with new ideas.
\end{rmkt}

The next algorithm describes our marginal evolution under relaxed conditional uniformity.

\begin{algt}[Marginal Evolution of the Tiling]
\label{alg:cf:tiling:marg-evo}
We now describe a single step of the evolution of the tiling.
Multiple steps are obtained by repeating the single-step evolution.

The single-step evolution is given by a map $\TT : (\pp; \uu, \vv; \bb) \to \pp'$ with inputs as follows:
\begin{center-small}
	a tiling
		$\pp \in \mfp$;\quad%
	two markers
		$\uu, \vv \in \bra{1/\nn, ..., 1}$;\quad%
%	are the first and second markers in \cref{alg:cf:cu:relaxed}, respectively,
%and
	a coin toss
		$\bb \in \bra{0,1}$.
\end{center-small}

Let $\II$ and $\JJ$ be the tiles containing $\uu$ and $\vv$, respectively.
Reorder the tiles in $\pp$ so that $\II$ is at the left and $\uu = 1/\nn$.
We use $\uu$ and $\vv$ as the first and marker in \cref{alg:cf:cu:relaxed}, respectively.
%There are two possibilities.
	%
\begin{itemize}
	\item 
	If $\II \ne \II'$, then merge tiles $\II$ and $\II'$ into a single tile.
	The new tile has size $\abs \II + \abs \II'$.
%	Set $\II'$ to be this new, single tile.
	
	\item 
	If $\II = \II'$, then propose a split into two fragments at $\vv - 1/\nn$.
	Split if and only if $\bb = 1$:
		the new left-hand tile  has size $\vv - 1/\nn$
	and
		the new right-hand tile has size $\abs \II - (\vv - 1/\nn)$.
%	\begin{itemize}
%		\item 
%		If $\bb = 0$, then do not split.
%%		Set $\II' \cq \II$ to be this single, unsplit tile.
%		
%		\item 
%		If $\bb = 1$, then do split at location $\vv$, giving rise to two tiles:
%			the left  has size $\vv - 1/\nn$
%		and
%			the right has size $\abs \II - (\vv - 1/\nn)$.
%%		This results in two tiles, the left of size $\vv - 1/\nn$ and the right of size $\abs \II - (\vv - 1/\nn)$.
%%		Set $\II'$ to be the right-hand tile, which contains $\vv$.
%	\end{itemize}
\end{itemize}
	%
%The current tiling is not necessarily ordered in non-increasing size.
The output tiling $\pp'$ is the tiling after this change written in non-increasing order.
%This describes the marginal evolution of the tiling.

This is extended to multiple steps by drawing $(\uu, \vv)$ as in \cref{alg:cf:cu:relaxed}, scaled by $\nn$, and letting $\bb \sim \Bern(\tfrac12)$ independently at each step.
The tile $\II \ni \uu$ is distinguished at each step.
\end{algt}

%\begin{rmkt}[Evolution of the Distinguished Tile]
%	%
%The variable $\ww$ accounts for the fact that our first marker in the subsequent step can be any of the markers used since the last refresh time.
%\bst[\S 4.2.2] take the second marker in one swap to be the first marker in the subsequent step.
%They thus take $\ww \cq \vv$;
%so $\II'$ is the tile containing the second marker $\vv$.
%%Recall the discussion in \cref{rmk:decomp:dist-supp:param-motivation}.
%	%
%\end{rmkt}

%It remains to define the distinguished tile.
%We let the variable $\ww$ represent the first marker in the subsequent step; this is uniform amongst the markers used since the last refresh time.
%Define $\II'$ to be the tile containing $\ww$.
	%
%\end{subtheorem}

%\begin{lem}[Marginal Evolution of the Tiling]
%	The evolution described by \ref{alg:cf:tiling:marg-evo} has the relaxed law, in the sense of \cref{def:cf:cu:relaxed}.
%\end{lem}
%
%\begin{Proof}
%	%
%This is immediate from the definition of relaxed conditional uniformity, \cref{alg:cf:cu:relaxed}.
%	%
%\end{Proof}

A continuous version of these dynamics has also been studied.
There, $\uu, \vv \sim \Unif((0,1])$, $\bb \sim \Bern(\theta)$ and $\ww = \vv$.
$\theta \cq \tfrac12$ corresponds to our set-up, but general $\theta$ can be studied.
The following lemma was proved by \textcite{P:partitions-pd} for general $\theta$.
The case $\theta = 1$ was proved by \textcite{T:partitions-pd}.
See \cite[Theorem~7.1]{GUW:quantum-heisenberg} for another proof.
The main result of \textcite{S:random-trans} is that a corresponding tiling for random transpositions converges weakly to this continuous limit.
It is important for \bklmt too; see \bklmc[Lemma~5.2].

\begin{lem}[Invariant Distribution of Tiling]
\label{res:cf:tiling:pd}
	The $\PD(\theta)$ distribution is invariant for the continuous version of the tiling dynamics in \cref{alg:cf:tiling:marg-evo} when splits are accepted with probability $\theta$.
\end{lem}

\vspace*{-\smallskipamount}

We now describe a coupling of two tilings.
The history of this coupling was described at the start of this subsection.
Let $\PP$ and $\QQ$ denote two realisations which are to be coupled.
We describe a single round, as for the marginal evolution in \cref{alg:cf:tiling:marg-evo}.
Multiple rounds are obtained by repeating the single round procedure as described in \cref{alg:cf:tiling:marg-evo}, using \cref{alg:cf:cu:relaxed}.
We use the notation $\widebar \PP = (\widebar \PP_\ss)_{\ss\ge0}$ and $\widebar \QQ = (\widebar \QQ_\ss)_{\ss\ge0}$ for a single round, analogously to before.

Suppose that the current partitions are $\pp$ and $\qq$.
Create two tilings of $(0,1]$ using $\pp$ and using $\qq$.
We differentiate between the blocks that are \textit{matched} versus those that are \textit{unmatched}:
\begin{center-small}
	two blocks from $\pp$ and $\qq$ are \textit{matched} if they are of identical size.
\end{center-small}\noindent%
It may be possible to match the blocks in multiple ways; we choose an arbitrary matching.
Our goal is to match blocks as quickly as possible, but in a way that does not create small unmatched blocks.
Blocks are chosen at rate proportional to their size according to the marginal evolution.
Thus, large unmatched blocks are relatively easy to fix, but small blocks take longer to select.
%We do not use non-decreasing size-order to tile $(0, 1]$.
We place all matched parts to the right; the unmatched parts then occupy the left part.

\begin{figure}
	\centering
%\noindent%
\begin{minipage}{0.9\textwidth}
	\includegraphics[width=1\textwidth]{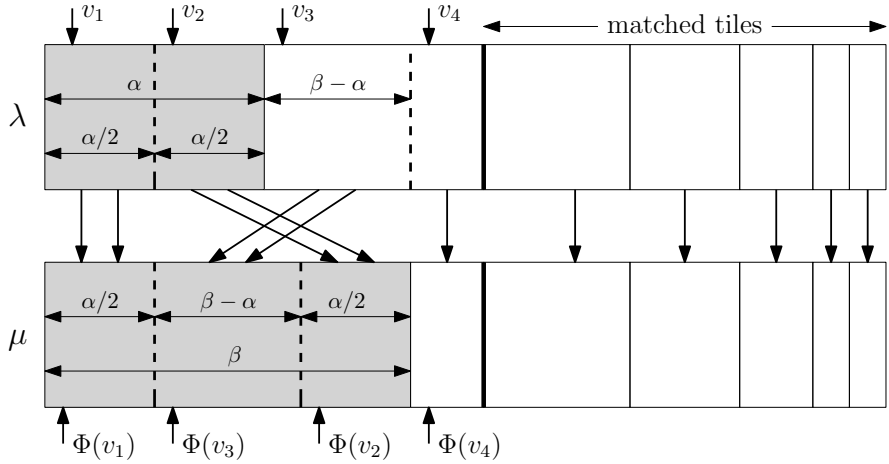}
	
	\caption{%
			Two tilings $\pp$ and $\qq$.
			The grey shaded tiles are the distinguished tiles, namely
				$I \in \pp$
			and
				$J \in \qq$.
			They have width $\abs I = \alpha$ and $\abs I = \beta$, respectively.
			The arrows represent the map $\Phi$ from \cref{alg:cf:tiling:coupling}.
			Four example pairs $(v_i, \Phi(v_i))_{i=1}^4$ are given.
		}
	\label{fig:schramm_transformation}
\end{minipage}
\vspace*{-3mm}
\end{figure}

\noindent%
\begin{algt}[Coupling of Tilings]
\label{alg:cf:tiling:coupling}
Let $\cc \in \mfc$ be a cycle structure.
Let $\ss \in \bra{1, ..., d(\cc)}$.

Suppose that $\widebar \PP_\ss = \pp$ and $\widebar \QQ_\ss = \qq$.
Suppose that the chosen markers are $(\uu_\pp, \vv_\pp)$ and $(\uu_\qq, \vv_\qq)$, respectively, for step $\ss$; these will be chosen in a coupled way.
Let $\II$ and $\II_\qq$ be the tiles containing the first markers $\uu_\pp$ and $\uu_\qq$, respectively.
Assume that either $\II_\pp$ and $\II_\qq$ are matched to each other or they are both unmatched.
We verify that this property is preserved by the coupling in \cref{res:cf:tiling:coupling:well-defined}.

We apply the transformation $\TT$ from the marginal evolution, \cref{alg:cf:tiling:marg-evo}.
We move the tiles $\II_\pp$ and $\II_\qq$ to the front of their respective tilings and assume that $\uu_\pp = 1/\nn = \uu_\qq$, as in \cref{alg:cf:tiling:marg-evo}.
We use the same coin toss $\bb$ for both $\pp$ and $\qq$.
We set
\(
	\pp' \cq \TT(\pp; \uu_\pp, \vv_\pp; \bb)
\Qand
	\qq' \cq \TT(\qq; \uu_\qq, \vv_\qq; \bb).
\)

It remains to construct $\vv_\qq$ as a function of $\vv_\pp$.
If $\II_\pp$ is matched, and hence is matched to $\II_\qq$, by assumption, then
	set $\vv_\qq \cq \vv_\pp$.
Now suppose that $\II_\pp$ is unmatched, and hence $\II_\qq$ is also unmatched, by assumption.
We obtain $\vv_\qq$ by applying a measure-preserving map $\Phi$ to $\vv_\pp$:
	set $\vv_\qq \cq \Phi(\vv_\pp)$.

It remains to define $\Phi$.
Let $\alpha \cq \abs{\II_\pp}$ and $\beta \cq \abs{\II_\qq}$ be the respective lengths of $\II_\pp$ and $\II_\qq$; assume that $\alpha \le \beta$, without loss of generality.
Denote $\gamma \cq \ceil{\alpha \nn/2 - 1}/\nn$.
Define
\[
	\Phi
:
	[0,1] \to [0,1]
:
	v
\mapsto
\begin{cases}
	v
	&\text{if}\quad v > \beta \text{ or } 2/\nn \le v \le \gamma + 1/\nn,
\\
	v - \gamma
	&\text{if}\quad \alpha < v \le \beta,
\\
	v + \beta - \alpha
	&\text{if}\quad \gamma + 1/\nn < v \le \alpha.
\end{cases}
\]
The map $\Phi$ is illustrated in \cref{fig:schramm_transformation};
cf \bsc[Figure~1].
% and \bklmc[Figures~5.2--5.7].
	%
\end{algt}

\vspace*{-\bigskipamount}

\begin{rmkt}[Potential Difficulties Arising from Rejecting Splits]
\label{rmk:cf:tiling:split-rejection}
Consider the scenario in which $\vv_\pp \notin \II_\pp$ but $\vv_\qq \in \II_\qq$.
The tile containing $\vv_\pp$ is always merged with $\II_\pp$ in $\PP$.
A split of $\II_\qq$ is proposed in $\QQ$, but may be rejected, namely if $\bb = 0$.
This scenario does not arise for random transpositions in \cite{S:random-trans} or conjugacy-invariant \RWs in \cite{BS:cutoff-conj-inv}; there, splits are always accepted.

It does not cause any issues for us, though. Indeed, larger tiles are easier to manage, since they are selected faster. Thus, not splitting is not be an issue.
\end{rmkt}

\begin{rmkt}[Weak Convergence vs Mixing]
\label{rmk:cf:tiling:weak-conv-vs-mixing}
One always takes $\vv_\qq \cq \vv_\pp$ in the original coupling of \textcite{S:random-trans}.
This is the case for \bklmt too, who adjust the coupling of \textcite{S:random-trans} to an application analogous to the $2$-\PM \RW.
Both of these articles study \emph{weak convergence}, rather than \emph{mixing}.
The adaptation to include the map $\Phi$ was one of the fundamental innocations introduced by \bszt; it was then used by \bst.
It is crucial when studying mixing, as we explain now.

If one takes $\Phi$ to be the identity, then this leads to the undesirable property that arbitrarily small, unmatched blocks may appear.
These small, unmatched blocks remain in the system for a long time---indeed, it takes a long time for them to even be selected. This prevents coalescence.
The map $\Phi$ rectifies this issue:
	the worst thing that can happen is for the smallest unmatched block to become smaller by a factor $2$ and this only happens with small probability;
	see \cref{res:cf:tiling:splitting}.
So, if the unmatched blocks are large initially, then they all remain relatively large.
The coupling quickly selects and matches large, unmatched blocks. Thus coalescence occurs quickly.
This is precisely why this map $\Phi$ was introduced by \bszt.
%, where it is one of the main innovations.

Another important property of the coupling is that the total number of unmatched blocks, ie the number in $\PP$ plus the number in $\QQ$, never increases; see \cref{res:cf:tiling:splitting} again.
\end{rmkt}

We now verify that the coupling is well-defined.

\begin{lem}[Coupling of Tilings]
\label{res:cf:tiling:coupling:well-defined}
	Suppose that the two distinguished tiles are either matched to each other or both unmatched at the start of a step of \cref{alg:cf:tiling:coupling}.
	Then, this is the also the case at the end of the step.
	
	The coupling is a genuine coupling, ie has the correct marginals, and is Markovian.
\end{lem}

\begin{Proof}
The first claim involves some routine case analysis; see \bszc[Lemma~4.3].
Merging two tiles in one system and \emph{not} splitting in the other preserves these conditions; cf \cref{rmk:cf:tiling:split-rejection}.

The second claim follows from the construction using \cref{alg:cf:cu:relaxed,alg:cf:tiling:marg-evo}.
\end{Proof}

%It will be convenient to define
%\[
%	\floor{x}_n
%\cq
%	\tfrac1n \floor{nx}
%\Qfor
%	x \in \mbr
%\Qand
%	n \in \mbn.
%\]

\begin{lem}[Evolution of Unmatched Blocks]
\label{res:cf:tiling:splitting}
	Let $\pp, \qq \in \mfp$ and let $\pp', \qq' \in \mfp$ be the corresponding integer partitions after one step of the coupling, ie of \cref{alg:cf:tiling:coupling}.
	Let $\UU$ and $\UU'$ be the sizes of the smallest unmatched block in the pair $(\pp, \qq)$ and $(\pp', \qq')$, respectively.
%	Choose $\rr \in \mbn$ so that $2^\rr \le \UU < 2^{\rr+1}$.
	The following hold:
	\begin{itemize}
		\item 
		\(
			\UU' \ge \tfrac1\nn \floor{\tfrac12 \UU \nn}
		\Qand
			\pr{ \UU' \le 2^{\floor{\log_2 U}} } \le 4 \UU / \nn;
		\)
		
		\item 
		the total number of unmatched partitions
		in $(\pp', \qq')$ vs $(\pp, \qq)$
		cannot increase.
	\end{itemize}
\end{lem}

\begin{Proof}
The proof of this lemma is almost the same as \bszc[Lemma~19]; see also \bsc[Lemma~4.4], where the details were omitted.
Merging two tiles in one system and \emph{not} splitting in the other preserves these conditions; cf \cref{rmk:cf:tiling:split-rejection}.
We omit the details here too.
\end{Proof}

We now compare properties of the coupling used in the current article with those of \bsc.
In essence, the key in \cite{BS:cutoff-conj-inv} is that the tiles do not get too small.
The \cf processes are the same except that we have the additional property of rejecting some splits.
Thus, it is simple to couple the two approaches so that blocks are larger in our process than in theirs.
This means that an identical proof as given in \bsc applies here, leading to \cref{res:cf:tiling:coupling:prob}~below.

Recall that $(\widebar \PP_\ss)_{\ss\ge0}$ and $(\widebar \QQ_\ss)_{\ss\ge0}$ denote two tilings, coupled in the above manner, on the swap-timescale. That is, incrementing $\ss$ to $\ss+1$ corresponds to applying a single \emph{swap}, not a full~\emph{round}.

For ease of presentation, assume that a single round involves choosing a divergent (in $\nn$) number of swaps.
This allows us to define easily $\widebar \PP_{\ceil{\delta^{-9}}}$ and $\widebar \QQ_{\ceil{\delta^{-9}}}$
with
%	$\Delta \cq \ceil{\delta^{-9}}$
%and
	$\delta > 0$ arbitrary but fixed.
Otherwise, simply concatenate sufficiently many rounds so that a least $\ceil{\delta^{-9}}$ swaps are made.
We are using the relaxed law and $\ceil{\delta^{-9}} = \Th1 = \oh \nn$, so this concatenation has no negative effects.

The next lemma shows that $\widebar \PP_{\ceil{\delta^{-9}}} = \widebar \QQ_{\ceil{\delta^{-9}}}$ \whp if the initial tilings $\widebar \PP_0$ and $\widebar \QQ_0$
	start with few unmatched blocks
and
	the smallest unmatched block is not very small.
%	specifically at least size $\delta$.
%We justify these assumptions on the initial tiling shortly.

\begin{lem}[Tiling Coupling Probability; cf {\bsc[Lemma~4.11]}]
\label{res:cf:tiling:coupling:prob}
	Let $\widebar \PP_0, \widebar \QQ_0 \in \mfp$ be two tilings.
	Assume that there are at most $3$ unmatched blocks between $\widebar \PP_0$ and $\widebar \QQ_0$.
	Write $A_\delta$ for the event that the size of the smallest unmatched block is at least $\delta$, for $\delta > 0$.
	Then,
	\[
		\LIMSUP{\tozero \delta}
		\LIMSUP{\toinf  \nn}
		\pr{
			\widebar \PP_{\ceil{\delta^{-9}}}
		\ne
			\widebar \QQ_{\ceil{\delta^{-9}}}
		}
		\ONE\rbr{ A_\delta }
	=
		0.
	\]
\end{lem}

\begin{Proof}
This lemma follows in a completely analogous way to how \bsc[Lemma~4.11] does for the conjugacy-invariant \RW on the symmetric group.
Rejecting splits has no ill effects.
%The fact that we reject half the proposed splits on average
%The limiting invariant distribution is $\PD(1)$ there,
%while it is $\PD(\tfrac12)$ here, due to the rejection of half the proposed splits.
%This change is inconsequential.
	%
\end{Proof}

We now briefly justify why we need only consider initial tilings with few unmatches blocks and smallest unmatched block not very small.
Rigorous analysis comes later in \cref{res:3coup:time:init}.

\begin{Proof}[Justification of Assumptions in \cref{res:cf:tiling:coupling:prob}]
\qedtriangle
We perform a path coupling approach.
We start at swap distance $1$.
It is easy to couple the tilings so that the swap distance remains $1$ via a simple relabelling.
The associated tilings then have at most $3$ unmatched blocks.

We use a `burn-in' phase.
This will be long enough so that the associated tilings look \emph{roughly} like they should in equilibrium. In particular, there will be few very small blocks.
We are able to deduce that the smallest unmatched block has size order $1$, ie not vanishing with $\nn$, \whp.
%
%The details are made rigorous in \cref{res:3coup:time:init}.
	%
\end{Proof}

\section{The Three-Stage Coupling of Two Systems}
\label{sec:3coup}

\subsection{Definition of Coupling}
\label{sec:3coup:def}

The overall coupling has three stages; cf\bsc[\S4.2].
It is trivial to couple two \PM systems, so that their relative swap distance remains constant, even on the swap-timescale, via a simple relabelling; see \cref{def:3coup:coup:const-swap-dist,res:3coup:coup:const-swap-dist}.
We call this the \textit{distance-preserving coupling}.
It is this simple coupling which is used in Stages~\ref{enm:3coup:def:stages:1} and~\ref{enm:3coup:def:stages:3}.
Stage~\ref{enm:3coup:def:stages:2} uses our adaptation of \textcite{S:random-trans}'s coupling.

\begin{enumerate}[itemsep = 0pt, topsep = \smallskipamount, leftmargin = 15mm]
\renewcommand{\labelenumi}{\upshape Stage \arabic{enumi}}
\renewcommand{\theenumi}{\upshape \arabic{enumi}}
\item \label{enm:3coup:def:stages:1}
The first stage is a burn-in period. It uses the distance-preserving coupling.
We wish the burn-in period to end in such a configuration that the two tilings have few unmatched blocks and any unmatched blocks are not vanishingly small; cf \cref{res:cf:tiling:coupling:prob}.
The length of the burn-in period asymptotically dominates the other two stages.

\item \label{enm:3coup:def:stages:2}
The next stage uses the adaptation of Schramm's coupling described in \S\ref{sec:cf:tiling}.
It will be run for time order $1$ on the swap-timescale; it does not necessarily involve an integer number of steps on the \PM-timescale, which involves approximately $\kk$ swaps.

\item \label{enm:3coup:def:stages:3}
The final stage simply finishes off \PM initiated in the second stage so that an integer number of \PMs have been applied.
It uses the distance-preserving coupling.
\end{enumerate}

We construct the distance-preserving coupling one swap at a time.
The informal idea is simple:
\begin{enumerate}[noitemsep]
	\item 
	relabel in the two \PMs so that each is at the identity;
	
	\item 
	draw a new matching and replace the identity with this new matching;
	
	\item 
	undo the relabelling in each \PM.
\end{enumerate}
%We give a precise and formal definition.

\begin{defn}[Swap Distance-Preserving Coupling]
\label{def:3coup:coup:const-swap-dist}
Suppose that the two \PM \RWs are at $\mu$ and $\nu$, respectively.
%Suppose that the current states of the two chains are $\mu$ and $\nu$, respectively.
Choose relabellings $\sigma$ and $\tau$ which translate $\mu$ and $\nu$ to the identity, respectively:
%\[
%	\cup_1^\nn
%	\bra{ \bra{ \mu_{\sigma(2\ell-1)}, \mu_{\sigma(2\ell)} } }
%=
%	\cup_1^\nn
%	\bra{ \bra{ 2\ell-1, 2\ell} }
%=
%	\cup_1^\nn
%	\bra{ \bra{ \nu_{\tau(2\ell-1)}, \nu_{\tau(2\ell)} } }.
%\]
\[
	\mu_{\sigma(\ii)}
=
	\ii
=
	\nu_{\tau(\ii)}
\Qforall
	\ii \in [2\nn].
\]
Now choose a \PM \uar,
say
\(
	\eta
=
	\cup_1^\nn
	\bra{ \bra{ \eta_{2\ell-1}, \eta_{2\ell} } }.
\)
Define $\mu'$ and $\nu'$ by `undoing' the relabelling of $\sigma$ and $\tau$, but starting from $\eta$ rather than the identity:
\[
	\mu'_\ii \cq \eta_{\sigma^{-1}(\ii)}
\Qand
	\nu'_\ii \cq \eta_{\tau^{-1}(\ii)}
\Qfor
	\ii \in [2\nn].
\]

A single swap, rather than the full \PM, is obtained by decomposing the new \PM $\eta$ into individual swaps and applying one at a time, choosing the relabellings $\sigma$ and $\tau$ anew each time.
\end{defn}

\begin{lem}[Swap Distance-Preserving Coupling]
\label{res:3coup:coup:const-swap-dist}
	The coupling of \cref{def:3coup:coup:const-swap-dist}
		is a genuine coupling of the \PM \RWs,
		is Markovian
	and
		preserves the swap distance,
		even on the swap-timescale.
\end{lem}

\begin{Proof}
This is immediate from the construction.
\end{Proof}

Recall the parameters $\kappa$ and $\rho$ representing, respectively, the expected support and expected distance of a uniformly chosen $\kk$-\PM from \cref{def:decomp:dist-supp:param-exp}.
This gives rise to an `average' \PM- and a swap-timescale.
The precise timescales are only well-defined if a sequence of \CSs are prescribed in advance and the \PMs or swaps are chosen conditional on this.

%\begin{rmkt}[Comparison Between \PM- and Swap-Timescales]
%\label{rmk:3coup:def:timescale-comparison}
%	%
%If $\eta \in \mfm'$---ie, an $\nn$-\PM with at most $\kk$ non-fixed points---is chosen,
%then we apply $d(\eta)$ swaps.
%Recall that $\rho = \ex{ d(\Unif(\mfm_\kk)) }$.
%So, the swap-timescale runs, on average, a factor $1/\rho$ slower than the \PM-timescale.
%That is, on average, $\rho$ steps on the swap-timescale corresponds to a single matching on the \PM-timescale.
%	%
%\end{rmkt}

We used $\ceil{\delta^{-9}}$ swaps used with our adaptation of Schramm's coupling.
%in \cref{res:cf:tiling:coupling:prob}.
Coalescence is achieved with probability tending to $1$ as $\delta \to 0$ on the event there are initially at most $3$ unmatched tiles and the unmatched tiles have size at least $\delta$; see \cref{res:cf:tiling:coupling:prob}.
If this event fails, then we use the distance-preserving coupling instead.
Importantly, the number of unmatched tiles is non-increasing under Schramm's coupling, so the relative distance of the \PMs remains at most $2$.

Recall from \cref{def:decomp:dist-supp:param-exp} that $\kappa = \ex{ \# H }$ where $H \sim \Unif(\mfm_\kk)$; it is the expected support of a uniform $\kk$-\PM, or equivalently of a uniform $\kk$-rematching.
Roughly, this is the number of uniform choices per round.
It is thus natural for our times, such as the mixing time, to scale inversely in $\kappa$.

\begin{defn}[Three-Stage Coupling]
\label{def:3coup:def:def}
We define the three coupling stages on the swap-timescale:
\begin{gather*}
	\text{\textup{Stage 1}} \ \text{ is } \
	[0,     \ss_1),
\quad
	\text{\textup{Stage 2}} \ \text{ is } \
	[\ss_1, \ss_2)
\Qand
	\text{\textup{Stage 3}} \ \text{ is } \
	[\ss_2, \ss_3),
\\
	\text{where}
\quad
	\ss_1 \cq \floor{ (\beta \nn - \delta^{-9}) / \kappa } \rho,
\quad
	\ss_2 \cq \ss_1 + \ceil{\delta^{-9}}
\Qand
	\ss_3 \cq \ceil{ \beta \nn / \kappa } \rho.
\end{gather*}
We use the distance-preserving coupling of \cref{def:3coup:coup:const-swap-dist} in Stages~\ref{enm:3coup:def:stages:1} and~\ref{enm:3coup:def:stages:3}.
If the smallest unmatched block has size at least $\delta$ at time $\ss_1$, then we use Schramm's coupling in Stage~\ref{enm:3coup:def:stages:2}, lifted to the \PM chain; otherwise, we use the distance-preserving coupling.
\end{defn}

\subsection{Coupling Time from Neighbouring Perfect Matchings}
\label{sec:3coup:time}

Suppose that we start with two neighbouring \PMs, ie ones which differ by a single swap.
Monotonicity of the number of unmatched blocks in the tilings implies that there are always at most $3$ unmatched blocks and thus that the relative distance of the \PMs is always at most $2$.

Recall the definition of $A_\delta$:
	the smallest unmatched block in the tiling has size at least $\delta$.
Our first aim is to estimate the probability that the chains jointly satisfy $A_\delta$ at time $\ss_1$, in the limit $\delta \to 0$.
%Our first aim is to show that the chains at time $\ss_1$ jointly satisfy the event $A_\delta$, which says that the smallest unmatched block in the tiling has size at least $\delta$, with some probability depending on $\beta$, then take the limit $\delta \to 0$.
This is the content of \cref{res:3coup:time:init}.
It requires an adaptation of a hyper-graph argument introduced by \bst[\S 3].
Given that $A_\delta$ is satisfied at time $\ss_1$, we use \cref{res:cf:tiling:coupling:prob} to coalesce the chains using Schramm's coupling with probability tending to $1$ as $\delta \to 0$.
The monotonicity in the number of unmatched blocks in Schramm's coupling means that the swap distance remains bounded by $2$, even if this coalescence fails.

The following result controls the contraction in the relative distance between two \PMs.
The definition of the contraction rate $\theta(\beta)$ and threshold $\beta_0 \in (0, \infty)$ in \cref{res:3coup:time:contr} below are given in \bsc[Lemma~2.1], but the precise definitions are not important.
What is important is that $\theta(\beta)$ is the asymptotic proportion of vertices in the giant component of an auxiliary graph process discussed in the next chapter. This is discussed more in the following two results and proofs.
%All that is important is that $1 - \theta(\beta)^2 \approx e^{-\beta}$---precisely,
%\(
%	\beta / \log\rbr{ 1 - \theta(\beta)^2 }
%\to
%	1
%\)
%\asinf \beta; see \bsc[Lemma~2.4].

\begin{prop}[Relative Distance Contraction]
\label{res:3coup:time:contr}
	Let $\MC {\widebar \Mu} \ss$ and $\MC {\widebar \Nu} \ss$ be two $\kk$-\PM \RWs on the swap-timescale.
	Suppose that $d(\widebar \Mu_0, \widebar \Nu_0) = 1$.
	Let $\beta \in (\beta_0, \infty)$.
	Define $\ss \cq \floor{ \beta \nn / \kappa } \rho$.
	Then
	\[
		\LIMSUP{\delta \to 0}
		\LIMSUP{\nn \to \infty}
		\ex{ d(\widebar \Mu_\ss, \widebar \Nu_\ss) }
	\le
		1 - \theta(\beta)^2.
	\]
\end{prop}

Key to proving this proposition is controlling the probability of the event $A_\delta$ at time $\ss_1$.
The skeleton argument above implies that we need to prove the following result; also, recall \cref{res:cf:tiling:coupling:prob}.

\begin{lem}[Properties at the Start of Stage~\ref{enm:3coup:def:stages:2}; {cf \bsc[Lemma~4.2]}]
\label{res:3coup:time:init}
	Consider two $\kk$-\PM \RWs on the swap-timescale, say
		$\widebar \Mu = \MC{\widebar \Mu}{s}$ and $\widebar \Nu = \MC{\widebar \Nu}{s}$
	with associated tilings
		$\widebar \PP = \MC{\widebar \PP}{s}$ and $\widebar \QQ = \MC{\widebar \QQ}{s}$, respectively.
	Use the distance-preserving coupling from \cref{def:3coup:coup:const-swap-dist}.
%	Couple $(\widebar \Mu, \widebar \Nu)$ using the coupling from \cref{def:3coup:coup:const-swap-dist}.
	Suppose that $d(\widebar \Mu_0, \widebar \Nu_0) = 1$.
	Assume that $\beta \in (\beta_0, \infty)$ and run for time $\ss_1$.
	The following hold.
	
	\begin{itemize}
		\item 
		There are at most $3$ unmatched tiles in $(\widebar \PP_{\ss_1}, \widebar \QQ_{\ss_1})$.
		
		\item 
		Write $A_\delta$ for the event that smallest unmatched block has size at least $\delta$ in $(\widebar \PP_{\ss_1}, \widebar \QQ_{\ss_1})$.
		Then,%
		\[
			\LIMINF{\tozero \delta}
			\LIMINF{\toinf  \nn}
			\prt{ A_\delta }
		\ge
			\theta(\beta)^2.
		\]
	\end{itemize}
\end{lem}

\cref{res:3coup:time:contr} follows relatively easily from \cref{res:3coup:time:init}, as we now show.

\begin{Proof}[Proof of \cref{res:3coup:time:contr}]
We start $\widebar \Mu$ and $\widebar \Nu$ at swap-distance $1$, ie $d(\widebar \Mu_0, \widebar \Nu_0) = 1$.
We use the distance-preserving coupling until $\ss_1$, on the swap-timescale.
If $A_\delta$ holds at this time, then we run for a further $\ceil{\delta^{-9}}$ units of time, now using Schramm's coupling.
Coalesce is achieved with probability tending to $1$ as $\delta \to 0$, by \cref{res:cf:tiling:coupling:prob}.
Further, the relative distance is always bounded $2$. Thus the expected relative distance tends to $0$ as $\delta \to 0$ on this event.
If $A_\delta$ does not hold, then we use the distance-preserving coupling.
We use the distance-preserving coupling in $[\ss_2, \ss_3 = \ss)$.
The result now follows from \cref{res:3coup:time:init} which controls the probability of $A_\delta$.
\end{Proof}

The proof of \cref{res:3coup:time:init} is an adaptation of that of \bsc[Lemma~4.2].
It requires the construction of an auxiliary graph process, similar to that in \bsc[\S 3], which we give in \S\ref{sec:graph}.
We now explain how to deduce \cref{res:3coup:time:init} assuming results on that graph process, referencing \S\ref{sec:graph}.

\begin{Proof}[Proof of \cref{res:3coup:time:init}]
\bst[\S 3] introduce an auxiliary graph process to control the sizes of unmatched tiles.
\cref{res:3coup:time:init} will follow analogously to \bsc[Lemma~4.2] once we have constructed an auxiliary graph process in a suitably analogous way to \bsc[\S 3].
We construct such an auxiliary graph process in \S\ref{sec:graph}; see, in particular, \cref{alg:graph:constr:pm}.
Precisely, we use the argument of \bsc[Lemma~4.2] along with \cref{res:graph:var:giant},
which is analogous to \bsc[Theorem~3.1],
and the relation between \cref{alg:graph:constr:pm} and that of \bsc[\S 3],
described in \cref{rmk:graph:equivalence}.

The key part of the proof of \bsc[Lemma~4.2] is the following (paraphrased).
\begin{quote}
	Let $A_1$ [no relation to $A_\delta$] be the event that the four points comprising the two transpositions fall within the largest component of the associated graph at time $\ss_1 \approx \beta n \rho / \kk$.
	The relative size of the giant component converges to $\theta(\beta)$.
	Thus $\pr{A_1} \to \theta(\beta)^4$.
\end{quote}
There are some minor parity constraints in \bsc. This is why they consider \emph{two} transpositions. We need only consider a single swap, which corresponds to a single transposition.
The symmetry of the problem implies that the two labels in this swap may be chosen uniformly at random without replacement.
This is why the limiting probability is the product of the limiting probabilities that the individual labels are in the giant.
We also replace $\kk$ with $\kappa$, as discussed in \cref{rmk:decomp:dist-supp:param-motivation}.

There is one further part in the proof of \bsc[Lemma~4.2] which does not obviously transfer to our set-up and proof.
It goes as follows (paraphrased).
\begin{quote}
	The rescaled cycles sizes at time $s_1$ converge in distribution to a $\PD(1)$ random variable, multiplied by $\theta(\beta)$; see \bsc[Theorem~3.6].
	This implies that, conditional on the event $A_1$ above, the relative size of the cycles containing the four points comprising the two transpositions can be thought of as the size of four independent samples from a $\PD(1)$ distribution, multiplied by $\theta(\beta)$.\footnote{%
		Formally, a $\PD(1)$ random variable corresponds to an tiling of $(0,1]$ broken up into infinitely many blocks.
		Draw $U \sim \Unif(0,1)$ and let $S$ be the size of the tile containing $U$.
		This $S$ is what \bsn mean by ``size''}
	The probability that any one of these four samples has a size smaller than $\delta / \theta(\beta)$ tends to $0$ \aszero \delta.
\end{quote}
The cycle sizes in our set-up are somewhat different to those in \bsc:
	splits are rejected half the time here,
	while they are never rejected there.
This means that the cycle sizes in our set-up are stochastically larger than those in theirs. Thus, the same lower bound holds.

We expect our cycle sizes to follow a $\PD(\tfrac12)$ limiting distribution---see \cref{rmk:3coup:time:PD(1/2)}---but such a refined statement is not required for the simple lower bound described above.
\end{Proof}

\begin{rmk}[Poisson--Dirichlet Convergence of Rescaled Cycle Sizes]
\label{rmk:3coup:time:PD(1/2)}
Convergence in distribution for the cycles of the $\kk$-\PM \RW to $\PD(\tfrac12)$ seems extremely likely to hold.
\bsc[Theorem~3.6] is described by the authors as a ``simple adaptation of the proof of \textcite{S:random-trans}''; they provide some, but not all, of the details in their appendix.
We have already seen how the ``interchange process with reversals'' studied by \bklmt is analogous to the $2$-\PM \RW.
Their main result \bklmc[Theorem~1.1] is that the appropriately rescaled cycle sizes converge to $\PD(\tfrac12)$.
See also \bklmc[Lemma~5.2]; cf \cref{res:cf:tiling:pd}.
The relevant proofs in \bklmc[\S 5.2, ``Schramm's coupling''], are ``identical or nearly identical to the corresponding proofs in \cite{S:random-trans}, so we omit the details, but give comments where there are differences related to the rejection of splits'' (paraphrased).
There are, of course, further arguments in \bklmc.

%In summary,
%	\bst use a ``simple adaptation of the proof of \citeauthor{S:random-trans}''
%and
%	\citeauthor{BKLM:interchange-rev}'s use arguments almost identical to those of \citeauthor{S:random-trans}.
%It is thus more than reasonable to expect that one should be able to adjust the arguments in \bklmc to extend from $\kk = 2$ to general $\kk$, using ideas from the proof of \bsc[Theorem~3.6] and our extension of \citeauthor{S:random-trans}'s coupling to $\kk > 2$.

We emphasise that our proof does not need this convergence. % in distribution.
Our cycles are at least as large as those used in \bsc.
Those cycles satisfied the required lower bounds.
Thus, ours do too.

We leave the question of convergence in distribution open.
We conjecture that a combination of
	the ideas from \bsc[Theorem~3.6], which are ``a simple adaptation of the proof of \citeauthor{S:random-trans}'',
	the ideas in \bklmc, particularly those which are ``nearly identical to the corresponding proofs by \citeauthor{S:random-trans}''
and
	our extension of \citeauthor{S:random-trans}'s coupling to $\kk > 2$
are sufficient to prove the claim.
\end{rmk}

\section{Auxiliary Graph Process}
\label{sec:graph}

Time has come to introduce and analyse the aforementioned auxiliary graph process, analogous to that introduced by \bst[\S 3].
There, the authors use a fixed \CS with support $\kk$ and consider a number $\tt$ of rounds with $\tt \eqsim \beta \nn / \kk$, for some $\beta \in (0, \infty)$.
Our set-up involves choosing the \CS randomly for each round.
The expected support $\kappa = \kappa_\kk$ satisfies
\(
	\kappa
=
	\kk - \tfrac12 + \Oh{\tfrac1\kk};
\)
recall \cref{def:decomp:dist-supp:param-exp,res:decomp:dist-supp:exp-supp}.
Thus, our number $\tt$ of rounds satisfies $\tt \eqsim \beta \nn / \kappa$.
%We thus consider a number $\tt$ of rounds satisfying $\tt \eqsim \beta \nn / \kappa$.

\subsection{Constructing the Auxiliary Graph Process}
\label{sec:graph:const}

We describe how to construct a hyper-graph in a way analogous to \bsc[\S3].
There, they apply a permutations $(\gamma_\tt)_{\tt\ge1}$, each with preset \CS. Such a permutation can be broken down into cycles, say $\gamma_\tt = \gamma_{\tt,1} \circ \cdots \circ \gamma_{\tt, \rr}$.
The hyper-edge $\bra{a_1, ..., a_\ell}$ is present in their hyper-graph at time $\TT$
	if and only if
$\gamma_{\tt, \ss} = \rbr{a_1, ..., a_\ell}$ for some $1 \le \ss \le \rr$ and $1 \le \tt \le \TT$.
The hyper-edge $\bra{a_1, ..., a_\ell}$ is independent of the order of its entries.
Thus, $a_1, ..., a_\ell$ may appear in order in the cycle $\gamma_{\tt, \ss}$.

This is equivalent to adding a clique with support $\bra{a_1, ..., a_\ell}$, ie adding all edges between $a_1, ..., a_\ell$, in a normal, non-hyper, graph.
The equivalence comes from the fact that we are only interested in the size of connected components.
We find this second viewpoint more natural.
%It is this second viewpoint which we feel is more natural.

Onto \PMs.
If we choose $\ell$ pairs to rematch into a cycle, say with labels $a_1, ..., a_\ell$, then we add a clique $\bra{a_1, ..., a_\ell}$.
We need a well-defined and consistent way of relabelling the matches after the rematching.
We explain precisely what we mean by this, since it is a key step.

Each \PM on $2\nn$ objects involves $\nn$ matches, or pairs, $\mm_1, ..., \mm_n$, which are labelled $1, ..., \nn$ in some manner. Suppose that we interact with the first two pairs, $\mm_1 = \bra{a, b}$ and $\mm_2 = \bra{c, d}$, giving rise to new matches $(\mm'_1, \mm'_2)$ satisfying $\mm'_1 \cup \mm'_2 = \bra{a, b, c, d} = \mm_1 \cup \mm_2$.
The quadruple $\bra{a, b, c, d}$ was initially matched as $\bra{\mm_1, \mm_2} = \bra{\bra{a, b}, \bra{c, d}}$.
Suppose that the interaction changes this to $\bra{\bra{a, c}, \bra{b, d}} = \bra{\mm'_1, \mm'_2}$.
There is no natural way of choosing $\mm'_1 \cq \bra{a, c}$ or $\mm'_1 \cq \bra{b, d}$.
%; recall \cref{rmk:decomp:gen:difficulty}.

We were in exactly the same quandary when sampling a uniform cycle via swaps.
%see \cref{rmk:decomp:gen:difficulty}.
%see \S\ref{sec:decomp:gen}.
We use here exactly the same solution as there: we choose the particular labelling in the new matching uniformly; cf \cref{alg:decomp:gen:gen-cycle}, which generates a uniform cycle via swaps.

Recall that
	if two objects in the same cycle are swapped,
	then we split components only half the time for the \PM \RW,
whereas
	splits always occur
	in this scenario
%	when two objects in the same cycle are transposed
	for the conjugacy-invariant \RWs of \bsc.
Importantly, this difference is irrelevant to the graph process since
hyper-edges/cliques are only \emph{added}, never \emph{removed}.
One can view this as merging two components of the graph when two cycles merge, unless they were already connected, but never splitting a component of the graph, even if a cycle splits.
Thus, this accept/reject of splits does not play~a~role.

\begin{algt}[Auxiliary Graph Process for the \PM \RW]
\label{alg:graph:constr:pm}
	%
%Fix $\nn \in \mbn$.
Let
$(\cc_\tt)_{\tt\ge1} \in \mfc^\mbn$
be a sequence of \CSs.
We construct a random graph process $(\GG_\tt)_{\tt\ge0}$.
%The graph $\GG_\tt$ has $\nn$ vertices for each $\tt \ge 0$.
We use an inductive construction.
Define $\GG_0 \cq ([\nn], \emptyset)$ to be the empty graph.
% on $\nn$ vertices.
Suppose that $\tt \ge 0$ and that $G_\tt$ has been defined.
We now define $G_{\tt+1}$.

\begin{itemize}
	\item 
	Choose the next \CS, ie $\cc_{\tt+1} \in \mfc$.
	
	\item 
	Choose a $\abs{\cc_{\tt+1}}$-subset of $[\nn]$ \uar, say $\bra{\bb_1, ..., \bb_{\abs{\cc_{\tt+1}}}} \subseteq [\nn]$.
	
	\item 
	Choose an associated partition\footnote{%
		the partition decides which elements of $\bra{\bb_1, ..., \bb_{\abs{\cc_{\tt+1}}}}$ go into which sub-cycle}
	$\pp \in \mfp_{\abs{\cc_{\tt+1}}}$ \uar.
	
	\item 
	Perform the following steps independently for each $\ell$-cycle in the decomposition $(\cc_{\tt+1}, \pp)$.\footnote{%
		The different $\ell$-cycles in the decomposition are disjoint.
		Thus, the order they are considered in is inconsequential}
	\begin{itemize}[noitemsep, topsep = 0pt]
		\item 
		Suppose that the labels of the $\ell$-cycle are $\bra{\aa_1, ..., \aa_\ell} \subseteq \bra{\bb_1, ..., \bb_{\abs{c_{\tt+1}}}} \subseteq [\nn]$.
		
		\item 
		Add the clique $\bra{\aa_1, ..., \aa_\ell}$,
		ie all edges between the vertices $\aa_1, ..., \aa_\ell$.\footnote{%
			Alternatively, if using the hyper-graph viewpoint,
			add the hyper-edge $\bra{\aa_1, ..., \aa_\ell}$}
		
		\item 
		Relabel the vertices $\aa_1, ..., \aa_\ell$ \uar.
%		\footnote{%
%			We could actually choose an arbitrary fixed relabelling or not relabel at all. This is because the partition $\pp$ is chosen uniformly.
%			Choosing uniformly from $(b_1, ..., b_\ell)$ is equivalent to choosing uniformly from a uniform permutation of these vertices}
	\qedhere
	\end{itemize}
\end{itemize}
\end{algt}

\begin{rmkt}[Comparison with {\bsc[\S 3]}]
\label{rmk:graph:equivalence}
We compare this graph process with the generalisation of that in \bsc[\S 3] for conjugacy-invariant \RWs, where we allow different \CSs to be picked at stage in a quenched sense.
This generalisation makes the algorithm for constructing the conjugacy-invariant graph process \emph{identical} to that used for \PMs, ie \cref{alg:graph:constr:pm} above, with one exception:
	the labels in the subset $\bra{\aa_1, ..., \aa_\ell}$ are randomised for the \PM version,
	but not for conjugacy-invariant version.

This relabelling is inconsequential.
% to the law of the graph process.
Indeed, the partition is chosen uniformly and independently each time.
The relabelling is only needed in order to couple with the \PM \RW.
This immediately gives a natural coupling between the \PM and conjugacy-invariant versions of the graph process.
	%
%\end{rmkt}

%\begin{rmkt}[Inconsequentiality of Relabelling]
	%
%The relabelling is inconsequential in a further sense.
%We fundamentally are only interested in sizes of connected components of the graph process.
Furthermore, we are only interested in the sizes of components later; see \cref{res:graph:var:giant}.
If two vertices are connected, then it does not matter to which of these two vertices other vertices are connected: the same connected component will be formed.

The relabelling is there only to circumnavigate the identifiability issue mentioned before.
\end{rmkt}

We have thus reduced the problem to a situation similar to that in \bsc.
There, a \CS is fixed and used forever: a $\cc \in \mfc$ is chosen and $\cc_\tt \cq \cc$ for all $\tt \ge 1$.
Now, the sequence $(\cc_\tt)_{\tt \ge 1}$ need not be constant.
The particular application that we are interested in is when each \CS is chosen independently and according to a uniform $\kk$-rematching in an $\nn$-\PM, ie $(\cc_\tt)_{\tt \ge 1} \sim \mcc\rbr{\Unif(\mfm')}^\mbn$.

\subsection{Approximating Variable Cycle Structures by a Fixed One}
\label{sec:graph:approx}

Our desire is to show that which particular \CS is used is irrelevant:
	in essence,
	all that matters is the rate at which an $\ell$-cycle is applied for each $\ell$.
We think of the growth of the graph process though an independent approximation.
The process involves breaking a $\kk$-\PM into single-cycle \PMs which are (weakly) correlated to previously applied single-cycle \PMs.
Ignore the correlations for the moment and just determine the law of the choice of single-cycle \PMs when the $\kk$-\PM~is~chosen~\uar.

If \CS $\cc$ is chosen, then $\cc_\ell$ is the number of $\ell$-cycles which are applied, for each $\ell$.
We can view this as a ``drawing balls from an urn'' problem in the following sense.
	Place $\BB \cq \sumt[\infty]{\ell=2} \cc_\ell$ balls in an urn:
		$\cc_\ell$ of colour $\ell$ for each $\ell \ge 2$.
	Set $\bb \cq 0$ and $\SS \cq \emptyset$.
Repeat the following steps until $b = B$.
\begin{itemize}[noitemsep]
	\item 
	If $\bb < \BB$, then draw a ball uniformly.
	Suppose that it is of colour $\ell$.
	
%	\begin{itemize}[noitemsep, topsep = 0pt]
		\item 
		Choose an ordered collection of $\ell$ elements \uar from $[\nn] \setminus \SS$.
		
		\item 
		Apply an $\ell$-cycle with this ordered collection.
%	\end{itemize}
	
	\item 
	Add these elements to $\SS$.
	Do not return the ball to the urn.
	
	\item 
	Increment $\bb$ by $1$.
%	Stop if $\bb = \BB$.
\end{itemize}
This perfectly simulates the application of a uniform \PM with \CS \cc.
We approximate by returning the ball to the urn and not updating the set $\SS$.
Each of the $\BB$ steps then has the same description.
\begin{itemize}[noitemsep]
	\item 
	Draw $\ell \in \mbn \setminus \bra{1}$ proportional to $(\cc_\ell)_{\ell=2}^\infty$.
	
	\item 
	Apply a uniformly chosen $\ell$-cycle.
\end{itemize}
A random number of $\ell$-cycles are applied in a single round; the expected number is $\cc_\ell$.
%We thus wish to determine the average number of $\ell$-cycles which are applied in a single round.

We extend this from always using the same \CS to define the law of this random number to choosing a random \CS for each round.
%Fix $\kk : \nn \mapsto \kk_\nn : \mbn \to \mbn$.
%Assume that $\kk_\nn/\nn \to 0$ \asinf \nn.
%It is standard to suppress the $\nn$-superscript from the notation.
%However, there are multiple limits going on here, so we leave it in for the sake of precision in some places.
Let $\Gamma \sim \mcc\rbr{\Unif(\mfm')}$, ie the \CS of a uniform $\nn$-\PM with at most $\kk$ non-fixed points.
Draw $\cc_\tt = (\cc_{\tt,\ell})_{\ell=1}^\infty \sim \mcc\rbr{\Unif(\mfm')}$ independently for each $\tt \ge 1$.
% with the same distribution as $\Gamma$.
A single round now involves applying $\gamma_\ell \cq \ex{\Gamma_\ell}$ $\ell$-cycles independently on average.
Then,
\[
	\# \gamma
=
	\sumt[\infty]{\ell=2}
	\ell \gamma_\ell
=
	\ex{ \# \Gamma }
=
	\kappa
\Qwhere
	\gamma
\cq
	(\gamma_\ell)_{\ell=1}^\infty.
\]
That is, $\gamma$ is \emph{almost} a \CS with support $\kappa$, ie the average support of a uniformly chosen $\kk$-\PM.
It is not quite, though, as each $\gamma_\ell$ need not be a non-negative \emph{integer}.

We would like to be able to say,
%\begin{quote} \slshape
	``Instead of choosing a random \CS in each step, just use the average $\gamma$, then apply some concentration results.
	This is legitimate since the order in which the $\ell$-cycles are applied is irrelevant \emph{for the random graph process}.''
%\end{quote}
The fact that $\gamma \notin \mbn_0^{\mbn_0}$ prohibits this.
It turns out to be unimportant, though.
We group together multiple steps and approximate those by a genuine \CS:
	roughly, we replace $\gamma_\ell$ with $\gamma_\ell' \cq \floor{\gamma_\ell / \eps} \in \mbn$, corresponding to $1/\eps$ steps.

%The only slight difference is that the values $\gamma_\ell$ may not be non-negative integers; they are merely restricted to being non-negative reals.
%Instead, $\gamma_\ell$ is the long-run average of the number of times an $\ell$-cycle is applied
%This distinction will turn out to be unimportant.

\subsection{Size of the Largest Component of the Graph}

The purpose of this section is to determine the proportion of vertices in the largest component of the auxiliary graph process, asymptotically \asinf \nn.
%This will be denoted $\theta$ below.
%
The following theorem is an adaptation of \bsc[Theorem~3.1] to our set-up.
In it, there is a critical threshold $\beta_0$ which $\beta$ must be above and a proportion $\theta(\beta)$, which will be the asymptotic proportion of vertices in the giant.
The particular values and definitions of these parameters is unimportant, but can be found in \bsc[Lemma~2.1].
%What matters is that $\theta(\beta)^2 \approx 1 - e^{-\beta}$ \asinf \beta.

Recall that we consider a number $\TT$ of round satisfying $\TT \eqsim \beta \nn / \kappa$. We make this precise now.
We evaluate the graph process of \cref{alg:graph:constr:pm} after this many rounds.

\begin{thm}[cf {\bsc[Theorem~3.1]}]
\label{res:graph:var:giant}
	There exists a critical threshold $\beta_0 \in (0, \infty)$ and a function $\theta : \mbr_+ \to (0, 1)$ with the following properties.
	Fix $\beta \in (\beta_0, \infty)$ arbitrarily.
	Suppose that $\TT_\beta$ satisfies $\TT_\beta \kappa / \nn \to \beta$ \asinf \nn.
	Consider the random graph process $(\GG_\tt)_{\tt \ge 0}$ evaluated at $\TT_\beta$.
	Then, the proportion of vertices which lie in the largest component
%	ie the size of the largest component divided by $\nn$,
	converges to $\theta(\beta)$ in probability \asinf \nn.
\end{thm}

%We also point out another fact of fundamental importance.
%The order in which the individual cycles are applied is inconsequential for the graph process $\GG_\tt$ at time $\tt$:
%	applying a single $3$-cycle in a round and then three in the next
%is equivalent to
%	applying two in each round,
%or even to
%	applying none in the first round an all four in the final round.
%The number of $\ell$-cycles applied across the $\tt$ rounds in total is $\tt \gamma_\ell$ in expectation.
%The number applied in each step is independent between steps and negatively correlated between different $\ell$.
%Thus laws of large numbers imply that these all concentrate around their means.
%So, while the intermediary evolution may differ slightly, the final result really is equivalent to what one would obtain with a fixed \CS $\gamma$.

We sketch the ideas behind \cref{res:graph:var:giant}.
% concerning the size of the giant component.
Even just the sketch proof is relatively technical.
We include the majority of the details, but suppress the explicit description of $1 \pm \oh1$ terms. Controlling these efficiently is more of a notational challenge than a mathematical one.
We trust that the details provided are sufficient for a masochistic\footnote{%
		\emph{masochist}: a person who enjoys an activity that appears to be painful or tedious}
reader to construct a rigorous proof.

We expect that the sketch is more complicated than it needs to be, but we have not found a simplification.
Indeed, we even conjecture that a `quenched' version of the theorem holds; see \S\ref{sec:graph:quenched}.
% for details.
%We expect that some more intelligent, less technical proof exists to establish this claim.

%Consider first a fixed \CS $\gamma \in \mfm'$,
%%ie not changing with time,
%as analysed by \bst[\S 3].
%It turns out that $\theta$ depends only on the proportions of the mass in cycles of a given size in $\gamma$, namely on
%\[
%	\widebar \alpha
%\cq
%	(\widebar \alpha_\ell)_{\ell=1}^\infty
%\cq
%	\LIM{\toinf \nn}
%	(\alpha_\ell)_{\ell=1}^\infty
%\Qwhere
%	\alpha_\ell
%\cq
%	\ell \gamma_\ell / \kk
%\in
%	[0,1]
%\Qfor
%	\ell \in \mbn.
%\]

The proof involves comparing our graph process with that of \bsc[\S 3] and applying \bsc[Theorem~3.1].
There are two key reductions.
We describe these two independently, then conclude.

\subsubsection*{Truncating the Cycle Sizes and Applying a Law of Large Numbers}

Let $\Gamma_\tt = (\Gamma_{\tt, \ell})_{\ell=1}^\kk \sim^\iid \mcc\rbr{ \Unif(\mfm') }$ for $\tt \in \mbn_0$, ie \iid $\kk$-rematchings in the space of $\nn$-\PMs.
Use \CS $\Gamma_\tt$ in round $\tt \in \mbn$.
Recall that $\gamma_\ell = \ex{\Gamma_{0, \ell}}$ is the number of $\ell$-cycles applied on average per round.

If $\gamma_\ell$ were an integer for each $\ell$, then we could simply use the \CS $\gamma = (\gamma_\ell)_{\ell=1}^\kk$ for each round and then conclude via a Law of Large Number (\textit{\LLN}). But alas, it is not.
In fact, $\aa_\ell \cq \gamma_\ell \ell \to 1$ as $\ell \to \infty$.
This is known for uniformly random permutations with $\aa_\ell \cq 1$ for all $\ell$. An analogous proof holds for \PMs; we omit the details.
If $\kk \asymp 1$, for example, then we can apply a \LLN to say that each $\ell$-cycle ($\ell \in \bra{1, ..., \kk}$) is applied a typically number of times.
However, if $\kk$ is sufficiently large, ie $\kk \gg \sqrt{\nn \log \nn}$, then the number of times that a $\kk$-cycle is applied is actually $\oh1$.
For such a large $\kk$, though, there is not significant difference between applying a $\kk$-cycle, a $(\kk-1)$-cycle, etc. We thus group together indices and assume that each group is applied a typical number of times.

We now proceed more formally.
Assume first that $\kk \to \infty$ \asinf \nn. We explain the easier $\kk$-bounded case after.
We group together indices $\ell$ which are `approximately equal'.
Let $\xi > 0$ with $\xi \to 0$ \asinf \nn, but vanishing as slowly as we desire.
Asymptotically, all the mass of the support comes from $\ell$-cycles with $\ell > \xi \kk \gg 1$. Indeed, this follows simply from the expectation $\gamma_\ell \asymp \ell$:
	\[
		\sumt[\floor{\xi \kk}]{\ell = 1}
		\ell \gamma_\ell
	\asymp
		\sumt{\ell \le \xi \kk}
		1
	=
		\oh{\kk},
	\quad
		\sumt[\kk]{\ell = \ceil{\xi \kk}}
		\ell \gamma_\ell
	\approx
		\sumt{\ell \ge \xi \kk}
		1
	=
		\kk (1 - \xi)
	\approx
		\kk
	\Qand
		\kk
	\approx
		\kappa,
	\]
	where the ``$\approx$'' signs hide $1 \pm \oh1$ factors,
	including $1 \pm \xi = 1 \pm \oh1$ factors.
%	which change from appearance to appearance.
We use the following grouping.
Let $\eps > 0$ with $\eps \to 0$ \asinf \nn, again vanishingly slowly.
Let
\[
	\II_\ii^-
\cq
	\xi \kk (1 + \eps)^\ii \wedge \kk,
\quad
	\II_\ii^+
\cq
	\xi \kk (1 + \eps)^{\ii+1} \wedge \kk
\Qand
	\II_\ii
\cq
	\oci{ \II_\ii^-, \: \II_\ii^+ }
%	\ocb{ \xi \kk (1 + \eps)^\ii, \: \xi \kk (1 + \eps)^{\ii+1} }
\cap
	\mbn
\Qfor
	\ii \in \mbn_0.
\]
Let
\(
	\ii_{\max}
\cq
	\inf\bra{ \ii \in \mbn \mid \xi \kk (1 + \eps)^{\ii+1} \ge \kk }.
\)
Then, for $0 \le \ii < \ii_{\max}$, we have
\[
	\sumt{\ell \in \II_\ii}
	\gamma_\ell
=
%	\sumt[\floor{\xi \kk (1 + \eps)^{\ii+1}}]{\ell = \ceil{\xi \kk (1 + \eps)^\ii}}
	\sumt[\floor{\II_\ii^+}]{\ell = \ceil{\II_\ii^-}}
	\aa_\ell / \ell
\approx
	\rbb{ {\xi \kk (1 + \eps)^{\ii+1}} - {\xi \kk (1 + \eps)^\ii} }
\big/
	\rbb{ \xi \kk (1 + \eps)^\ii }
=
	\eps.
\]
The \LLN along with a union bound over $\ii \in \bra{0, ..., \ii_{\max} - 1}$ gives
\[
	\pr{
		\sumt[\TT]{\tt=1}
		\sumt{\ell \in \II_\ii}
		\Gamma_{\tt, \ell}
	\approx
		\eps \TT
	\
		\text{uniformly}
	\
		\forall
	\:
		\ii \in \bra{0, ..., \ii_{\max} - 1}
	}
\approx
	1,
\]
where $\TT \approx \beta \nn / \kk$ is the number of rounds.
This requires the ``$\approx$'' sign inside the probability to be sufficiently weak compared with the decay of $\eps$ and $\xi$.
A similar bound holds jointly for $\ii_{\max}$, but taking into account the fact that $\abs{\II_{\ii_{\max}}}$ has a slightly different form, due to the truncation at $\kk$.

We use the following approximation $\GG'$ to the original graph process $\GG$:
if
	an $\ell$-cycle is applied in $\GG$ with $\ell > \xi \kk$,
then
	find $\ii$ with $\ell \in \II_\ii$ and apply an $\II_\ii^-$-cycle in $\GG'$.
The processes $\GG$ and $\GG'$ can easily be coupled so that $\GG_\tt \supseteq \GG'_\tt$ for all $\tt \ge 0$.
These $\ell$ and $\II_\ii^-$ satisfy $\ell \approx \II_\ii^-$ uniformly.
Thus, by continuity of $\theta$, it is still the case that $\GG'_\TT$ has a giant containing a proportion $\theta'(\beta) \approx \theta(\beta)$ of the vertices asymptotically.
This allows us to analyse $\GG'$ instead of $\GG$.

Analysis of $\GG'$ is still not trivial. We cannot apply ``an $\eps$-proportion of an $\II_\ii^-$-cycle'' in a single step.
We would like to simply `enlarge' the \CS by a factor $1/\eps$ and multiply the number of rounds by $\eps$.
We explain this concept via the following analogous situation.
\begin{itemize}[noitemsep]
	\item 
	Alternate between applying a $2$- and $3$-cycle; thus each is applied half the time.
	
	\item 
	`Enlarge' this by a factor $2$:
		apply a $2$- and $3$-cycle every round.
	
	\item 
	Divide the number of rounds by $2$:
		replace $\TT$ by $\tfrac12 \TT$.
\end{itemize}
This does not give rise to the same graph:
	choosing a $2$- and $3$-cycle in the same round conditions them to be disjoint;
	this is not the case when they are chosen in different rounds.
It is reasonable to suspect that this difference is minor, however.
Indeed, \bsc[Theorem~3.1] implies that this is the case when a fixed \CS is `doubled',
	ie there are twice as many $\ell$-cycles for each $\ell \ge 2$,
and the number of rounds is halved.
We show below that an analogous result holds for our $\eps$-application.

It remains to comment on the $\kk$-bounded case.
We do not need any rounding for this case since $\gamma_\ell \asymp 1$ uniformly and thus all $\ell$-cycles are applied a constant proportion of the time.
We simply condition that the number of $\ell$-cycles applied is typical for each $\ell$, of which there are $\kk \asymp 1$ different values.
We then rescale time by common denominator of $\gamma_1 \cdots \gamma_\kk$, which is order $1$.

We show next that these adjusted processes give rise to giants of the same size asymptotically.

\subsubsection*{Approximating Sampling without Replacement by Sampling with Replacement for Fixed \CS}

Suppose that an $\ell$-cycle is being applied and that indices $\aa = \bra{\aa_1, ..., \aa_\rr}$ have already been chosen this round, by the application of previous cycles.
$\ell$ indices $\bb = \bra{\bb_1, ..., \bb_\ell} \subseteq \bra{\aa_1, ..., \aa_\rr}$ are chosen uniformly without replacement from the restricted set $\bra{1, ..., \nn} \setminus \aa$.
The clique $\bb = \bra{\bb_1, ..., \bb_\ell}$ is added to the graph.
Suppose, instead, that we draw the indices with replacement and from the entirety of $\bra{1, ..., \nn}$:
	$\bb''_1, ..., \bb''_\ell \sim^\iid \Unif(\bra{1, ..., \nn})$;
	set $\bb'' \cq \bra{\bb''_1, ..., \bb''_\ell}$.
Certainly
\(
	\bb'' \setminus \aa
%	\bra{\bb''_1, ..., \bb''_\ell} \setminus \bra{\aa_1, ..., \aa_\rr}
\subseteq
	\bb
%	\bra{\bb_1, ..., \bb_\ell}
\)
stochastically.
Define the graph $\GG''$ via the indices $\bb'' \setminus \aa$ at the application of each cycle.
%, with $\GG$ defined via $\bb$ in the usual~way.

The fact that $\ell \le \kk \ll \nn$ implies that
\(
	\abs{\bb'' \setminus \aa''}
\approx
	\abs{\bb}
=
	\ell
\)
\whp.
In particular, for every $\ell \in \bra{1, ..., \kk}$, we can find an $\ell''$
such that
	$\ell'' \approx \ell$ uniformly
and
	at least $\ell''$ distinct elements are chosen \whp when an $\ell$-cycle is applied.
We think of this as ``rounding $\ell$ down to account for double counting''.

We can couple $\GG$ and $\GG''$ by adding an $\ell''$-clique to $\GG''$ whenever an $\ell$-clique is added to $\GG$.
There is some small probability that the inequality fails, but only a uniformly $\oh1$ probability.
%We compensate for this by increasing the number of rounds applied in $\GG''$:
	set $\TT'' \cq \beta \nn / \kappa''$,
	where
		$\kappa'' \approx \kappa$ is the support of this slightly reduced-size \CS.
%	and
%		$\beta'' \approx \beta$ is slightly larger than $\beta$ in such a way that at least $\TT = \beta \nn / \kappa$ rounds are successfully applied in $\GG''$ \whp.
Then, $\GG''_{\TT''}$ has a giant containing an asymptotic proportion $\theta''(\beta) \approx \theta(\beta)$ of the vertices,
using continuity and uniformity.
%where the final relation uses continuity of $\theta$ and uniformity of $\ell \approx \ell''$.

Finally, we release the restriction of
	applying exactly $\ell''$ to account for double counting in $\bb''$
and
	removing the previously-considered indices of $\aa''$:
	we simply choose $\ell$ uniformly with replacement and add this clique.
This only increases the size of the giant.

All in all, we have shown that
%if the number $\TT$ of rounds and support $\kappa$ of the \CS are increased by a $1 + \oh1$ factor, then
the giant of the graph in which the indices are sampled with replacement, rather than without replacement, contains a proportion $\theta(\beta)$ in probability.

\subsubsection*{Concluding Given the Above Reductions}

We conclude the sketch by combining the two reductions just established.
%It is more convenient to apply the second reduction first, even though they arose naturally in the other order.

\begin{enumerate}
	\item 
	Replace $(\GG, \TT)$ with $(\GG', \TT')$, where $\TT' \approx \TT$.
	This is the `rounding down' process, in which we apply an $\ell'$-cycle in $\GG'$ whenever an $\ell$-cycle is applied in $\GG$,
	where $\ell' \cq \II_\ii^-$ with $\ell \in \II_\ii$.
	We also condition that a typical number of each $\ell'$-cycles are applied and group these together.
%	This grouping together is what affects the time.
	
	\item 
	Replace $(\GG', \TT')$ with $(\GG'', \TT'')$, where $\TT'' \approx \TT'$.
	This replaces the ``sampling without replacement'' in each round with ``sampling with replacement''.
\end{enumerate}
Importantly, there is no longer a concept of ``multiple disjoint cycles in a single round'' when sampling \emph{with} replacement.
This means that the `enlargement' described at the end of the first part does not actually change the process at all.
Thus our \emph{random choice} graph process does indeed correspond, asymptotically, to the \emph{average choice}, encoded by $\gamma$ and $\theta$.

\subsubsection*{Alternative Proof: Copying \bsn's Argument from \bsc[\S 3]}

We believe that our \cref{res:graph:var:giant} can also be proved by following closely \bsn's proof of \bsc[Theorem~3.1] in \bsc[\S 3].
Doing so, one sees that the particular structure of $\Gamma$ is unimportant for their proof.
Indeed, this almost has to be the case since their argument works when $\Gamma$ comprises $\tfrac12 \kk$ disjoint transpositions, a single $\kk$-cycle or anything in-between.
\bsn give a helpful verbal summary of this lemma, which we lightly paraphrase.

\begin{quote} \slshape
	It is perhaps surprising that \bsc[Lemma~3.2] is sufficient for the proof of \bsc[Theorem~3.1].
	The lemma essentially only records whether a cycle is microscopic (finite) or ``more than microscopic''.
	In particular, whether the mass of the \CS comes from many small mesoscopic or fewer big cycles makes no difference.
\end{quote}

We have not checked carefully every detail in this argument.
Indeed, the reductions that we described above are sufficient for our annealed set-up, so there was no need.
However, the \LLNs we used would not be so amenable to the quenched set-up, described below.
The best way to prove a quenched statement may be to simply go through \bsc[\S 3], making the appropriate adjustments.
These are no doubt relatively easy conceptually, but likely challenging technically.

\subsubsection*{Convergence of Cycle Structure to Independent Poisson Process}

We remark for the sake of interest, rather than the proof, that the full vector of cycle lengths for a uniform permutation converges to that of an independent Poisson process in \TV if $\kk = \oh \nn$; see \cite{AT:cycle-structure-random-perm,B:comment-chen-stein,DP:permutations-unpublished}.
\TV analyses the entire vector:
	it is stronger than the more common weak convergence,
	which only analyses finite-dimensional marginals.
\textcite{B:comment-chen-stein} uses the Chen--Stein method, which approximates certain (weakly) dependent variables by independent Poisson random variables.
We have not checked carefully all the details, but we strongly suspect that the same argument can be used to establish convergence in \TV for a uniform \PM too.

\subsection{Conjectured Extension to General `Quenched' Cycle Structures}
\label{sec:graph:quenched}

We have done our best to leave the above description as general as possible.
In particular, we could estimate the law of the \CS of a uniform \PM.
The reasons for our not doing this are twofold.

First and most important, we do not need to.
The important term to control is $1 - \theta(\beta)^2$;
see \cref{res:3coup:time:contr}.
This is always approximately $e^{-\beta}$ in the limit \toinf \beta, regardless of the law;
see \bsc[Lemma~2.4] or \cref{res:cutoff:upper:beta-theta}.
In particular, if we draw the \CSs according to a different law, then this approximation still holds.
\cref{res:3coup:time:contr} is evaluated at $(\beta \nn / \kappa) \rho$ on the swap-timescale, which is equivalent to $\beta \nn / \kappa$ on the \PM-timescale.
The $\beta$ in the contraction $1 - \theta(\beta)^2 \approx e^{-\beta}$ and the $\beta$ in the time $\beta \nn / \kappa$ end up cancelling.
This is all made clear and rigorous in \S\ref{sec:cutoff:upper} below.

Second and more abstractly, the current formulation leads itself more naturally towards extension.
We do not really need anywhere the randomness in the choice of the \CS at each round.
For example, suppose that $\Gamma_0$ and $\Gamma_1$ are two fixed \CSs---say all transpositions ($2$-cycles) and all $3$-cycles, respectively.
Use $\Gamma_b$ in the $\tt$-th step if $\tt \equiv b$ mod $2$.
All our arguments would go through outputting the same results as if one of $\Gamma_0$ and $\Gamma_1$ were chosen uniformly and independently at each round.
The former is a `quenched' statement and the latter an `annealed'.

We believe that this can be extended even further.
If there is some `average behaviour' of the quenched sequence which manifests itself on the \PM-timescale order $\nn / \kappa$, then we expect that this `average behaviour' can be used to define $\theta$ appropriately.
Indeed, the graph process is insensitive to the order in which the different $\ell$-cycles are applied; it is `Abelian' in this sense. The coupling decomposes cycles into products of transpositions; it does not care what order these are applied or whether the transposition came from an $\ell$-cycle or an $\ell'$-cycle.
All that needs controlling carefully is the size of the small cycles and of the giant component after order $n$ swaps have been applied, however those swaps may arise; recall the proofs of \cref{res:3coup:time:contr,res:3coup:time:init}.

A quenched version of the lower bound actually holds easily.
We elaborate
in \cref{res:cutoff:lower:quenched}.
%in \S\ref{sec:cutoff:lower}.

%\newpage
\section{Upper Bound for Cutoff}
\label{sec:cutoff:upper}

The ideas in this concluding section follow closely those employed by \bst, differing only very slightly.
%The upper bound follows a path coupling argument.
Nevertheless, we include almost all the details for concreteness.

We are going to use the path coupling technique of \textcite{BD:path-coupling};
see \cite[Theorem~14.6]{LPW:markov-mixing} for a modern description.
The following proposition is a rephrasing of \cref{res:3coup:time:contr}, which is on the swap-timescale; the proposition below is given on the \PM-timescale.

\begin{prop}[Relative Distance Contraction]
\label{res:cutoff:upper:contr}
	Let $(\Mu_\tt)_{\tt\ge0}$ and $(\Nu_\tt)_{\tt\ge0}$ be two \PMs chains on the \PM-timescale.
	Suppose that $d(\Mu_0, \Nu_0) = 1$.
	Fix $\beta \in (\beta_0, \infty)$.
	Let $\TT_\beta \cq \floor{ \beta \nn / \kappa }$.
	Recall the contraction rate $\theta$ from \cref{res:3coup:time:contr}.
	There exists a coupling of $(\Mu_\tt)_{\tt\ge0}$ and $(\Nu_\tt)_{\tt\ge0}$ such that
	\[
		\lambda_\beta
	\cq
		\LIMSUP{\toinf \nn}
		\ex{ d(\Mu_{\TT_\beta}, \Nu_{\TT_\beta}) }
	\Quad{satisfies}
		\LIMSUP{\toinf \nn}
		\lambda_\beta
	\le
		1 - \theta(\beta)^2.
	\]
\end{prop}

We first informally justify the upper bound of $\tfrac1\kappa \nn \log \nn$ on the mixing time.
The standard path coupling bound says that the \TV distance after time $\mm \TT$ decays exponentially as $\rbr{ 1 - \theta(\beta)^2 }^\mm$.
There is a diameter pre-factor which is $\nn - 1 \approx n$.
Thus, to get \TV distance $\oh1$, we need
\[
	\mm
\approx
	- \log \nn / \log\rbb{1 - \theta(\beta)^2}.
\]

The function $\theta$ depends on the law of the \CS of a uniform \PM.
Somewhat surprisingly, however, we do not need to control this.
We use the following lemma which holds \emph{regardless} of the law.

\begin{lem}[cf {\bsc[Lemma~2.4]}]
\label{res:cutoff:upper:beta-theta}
	We have
	\[
		\LIM{\toinf \beta}
%		\frac{\beta}{\log\rbr{1 - \theta(\beta)^2}}
		\beta \big/ \log\rbb{1 - \theta(\beta)^2}
	=
		- 1.
	\]
\end{lem}

\begin{Proof}
The proof is elementary analysis.
See \bsc[Lemma~2.4] for analogous details.
\end{Proof}

This lemma then tells us, for the above $\mm$, that
\[
	\mm \TT_\beta
\approx
	- \tfrac1\kappa \nn \log \nn \cdot \beta / \log\rbb{1 - \theta(\beta)^2}
\to
	\tfrac1\kappa \nn \log \nn
\quad
	\asinf \beta.
\]
This informally justifies the upper bound of $\tfrac1\kappa \nn \log \nn$.
We now proceed formally and rigorously.

\begin{Proof}[Proof of Upper Bound in \cref{res:intro:main}]
Let $(\Mu_\tt)_{\tt\ge0}$ and $(\Nu_\tt)_{\tt\ge0}$ be two \PMs chains.
Recall the $d$ denotes the swap-distance; in particular, $d(\mu, \nu) \ge \one{\mu \ne \nu}$.
Thus,
\[
	\tv{ \Mu_\tt - \Nu_\tt }
\le
	\ex{ d(\Mu_\tt, \Nu_\tt) }
\Qforall
	\tt \ge 0,
\]
for any coupling of $\Mu_\tt$ and $\Nu_\tt$.
Let $\beta \in (\beta_0, \infty)$ and $\mm \in \mbn$.
Recall that $\TT_\beta = \floor{ \beta \nn / \kappa }$.
Iterating as in the path coupling method and applying \cref{res:cutoff:upper:contr} at each iteration, we obtain
\[
	\tv{ \Mu_{\mm \TT_\beta} - \Nu_{\mm \TT_\beta} }
\le
	\nn \lambda_\beta^\mm
=
	\nn \exp{ \mm \log \lambda_\beta }.
%	\nn \rbb{ 1 - \theta(\beta)^2 }^\mm.
\]
This uses the fact that
\(
	d(\Mu_0, \Nu_0)
\le
	\max_{\mu, \nu \in \mfm_\nn} d(\mu, \nu)
=
	\nn - 1
\le
	\nn.
\)

Let $\eps > 0$. We want the \TV distance to be at most $\eps$.
It thus suffices for $\mm$ to satisfy
\[
	\mm
\ge
	\mm_{\beta, \eps}
\cq
	\rbb{ \log n + \log(1/\eps) }
\big/
	\log\rbr{ 1 / \lambda_\beta }.
%	\frac
%		{\log \nn + \log(1/\eps)}
%		{\log\rbr{ 1 / \lambda^\nn_\beta }}.
%%		{- \log\rbr{ 1 - \theta(\beta)^2 }}.
\]
It thus suffices to consider $\tt$ with $\tt \ge \mm_{\beta, \eps} \TT_\beta$.
Let $\delta > 0$ be arbitrarily small but constant.
Set
\[
	\tt
\cq
	\rbr{ 1 + \delta } \rbr{ \nn \log \nn / \kappa }.
\]
\cref{res:cutoff:upper:contr,res:cutoff:upper:beta-theta}
%together
imply that we can choose $\beta_\delta$ and $\nn_{\delta, \eps}$
large enough
so that
\[
	\tt
\ge
	\mm_{\beta_\delta, \eps} \TT_{\beta_\delta}
\Qforall
	\nn \ge \nn_{\delta, \eps}.
\]
%as required.

This completes the upper bound
in \cref{res:intro:main}
as $\eps$ and $\delta$ were arbitrary.
%in the proof of cutoff as $\eps, \delta > 0$ were arbitrary.
	%
\end{Proof}

\section{Lower Bound for Cutoff}
\label{sec:cutoff:lower}

The lower bound is just a simple coupon-collector argument, using the number of fixed points as a distinguishing statistic.
We omit the details of this calculation, referencing to analogous ones.
%We include the concluding details.

A uniform $\kk$-\PM has $\kk/(2\kk-1) = \Th1$ fixed points \wrt the identity in expectation.
The application of an $\ell$-cycle involves choosing $\ell$ elements of $[\nn]$ \uar without replacement. This is approximately the same as choosing with replacement since $\ell \le \kk = \oh \nn$. In fact, if one does draw \uar with replacement, then the number $\ell'$ of draws required to get $\ell$ distinct elements satisfies $\ell' / \ell = 1 + \oh1$ \whp.
A coupon-collector argument shows that if only $(1 - \delta) \nn \log \nn$ uniform choices are made, then divergently many elements of $[\nn]$ will not have been selected.
The resulting \PM then has a divergent number of fixed points.
The number of fixed points thus acts as a distinguishing~statistics.
%Again, the order in which the $\ell$-cycles are applied does not matter. All that matters is the number of times an $\ell$-cycle has been applied for each $\ell$.

A formal and rigorous proof in the case of a fixed \CS is given by \bsn in their appendix, specifically \bsc[Appendix~A].
It can be adapted to prove the following result.

\begin{prop}[Fixed Points]
\label{res:cutoff:lower:fp}
	Let $\kk \in \mbn$ satisfy $\kk / \nn \to 0$ \asinf \nn.
	Let $(\cc_\tt)_{\tt \ge 1} \in (\mfm')^\mbn$ be an arbitrary sequence of \CSs, each corresponding to a $\kk$-rematching in an $\nn$-\PM.
	Let $(\Mu_\tt)_{\tt \ge 0}$ be the `quenched' \PM \RW in which \CS $\cc_\tt$ is used in round $\tt \ge 1$ with $\Mu_0 = \id$, the identity.
	Let
	\[
		\tt_\lambda
	\cq
		\inf\brb{ \tt \ge 0 \midb \sumt[\tt]{\tt'=0} \abs{\cc_{\tt'}} \ge \lambda \nn \log \nn }
	\Qfor
		\lambda \in (0, \infty).
	\]
	Fix $\lambda \in (0, 1)$ and $\KK \in \mbn$.
	Then the number of fixed points in $\MM_\tt$ is at least $\KK$ \whp if $\tt \le \tt_\lambda$.
\end{prop}

The lower bound on mixing follows easily from this.

\begin{Proof}[Proof of Lower Bound in \cref{res:intro:main}]
Suppose that the sequence $(\cc_\tt)_{\tt \ge 1}$ of \CSs is drawn.
%The expected support of a uniform $\kk_\nn$-\PM is $\kappa^\nn$, by definition.
The corresponding $\kk$-\PMs are chosen independently.
Let $\delta \in (0, 1)$, independent of $\nn$; set
\[
	\tt
\cq
	(1 - \delta) \rbr{ \nn \log \nn / \kappa }.
\]
Recall that $\nn / \kappa \ge \nn / \kk \to \infty$ \asinf \nn.
The Law of Large Numbers thus implies that
\[
	\sumt[\tt]{\tt'=0} \abs{\cc_{\tt'}}
\le
	(1 - \tfrac12 \delta) \nn \log \nn
\quad
	\whp,
\Quad{ie}
	\tt \le \tt_{1 - \delta / 2}
\quad
	\whp.
\]
There are thus divergently many fixed points in the \PM at time $\tt$ \whp, by \cref{res:cutoff:lower:fp}.
Contrastingly, the expected number of fixed points in a uniform \PM is at most $\kk/(2\kk-1) \le 1$, for any $\kk$; see \cref{res:decomp:dist-supp:exp-supp}.
The number of fixed points in the \PM is thus a distinguishing statistic.

This completes the upper bound
in \cref{res:intro:main}
as $\delta$ was arbitrary.
\end{Proof}

\begin{rmkt}[Extension to `Quenched' Cycle Structures]
\label{res:cutoff:lower:quenched}
This argument extends easily to quenched cycle structures, where the sequence $(\cc_\tt)_{\tt \ge 1}$ of \CSs is prescribed in advance, provided the support $\kappa_\tt \cq \# \cc_\tt$ is uniformly $\oh \nn$.
Define $\TT$ to be the natural coupon-collector threshold, ie
\[
	\TT
\cq
	\inf\brb{ \tt \ge 0 \midb \kappa_1 + \cdots + \kappa_\tt \ge \nn \log \nn }.
\]
Then, there are divergently many fixed points at $\tt \cq (1 - \delta) \TT$ \whp if $\delta > 0$ is independent of $\nn$.
\end{rmkt}

\renewcommand{\bibfont}{\sffamily\small}
\renewcommand{\bibfont}{\sffamily}
\printbibliography[heading = bibintoc, title = {Bibliography}]

\end{document}